\UseRawInputEncoding
\documentclass[12pt]{article}
\usepackage[colorlinks=true,linkcolor=blue,urlcolor=blue,citecolor=blue]{hyperref}
\usepackage{pgfplots}
\usepackage{mathrsfs}
\usepackage{amssymb}
\usepackage{amsmath}
\usepackage{amsfonts}
\usepackage{amstext}
\usepackage{amsthm}
\usepackage{tikz}
\usetikzlibrary{arrows,patterns,positioning}
\usepackage{graphicx}
\usepackage{subfig}
\usepackage{color}
\usepackage{indentfirst,latexsym,bm}
\usepackage{mathrsfs}
\usepackage{cite}
\usepackage[all]{xy}
\usepackage{newlfont}
\usepackage{enumerate}
\usepackage{array,longtable}
\numberwithin{equation}{section}
\theoremstyle{remark}\newtheorem*{remark*}{Remark}\theoremstyle{plain}
\newtheorem*{corollary*}{Corollary}
\newtheorem*{lemma*}{Lemma}
\newtheorem*{proposition*}{Proposition}
\newtheorem*{theorem*}{Theorem}
\allowdisplaybreaks[2]
\newtheorem{theorem}{Theorem}[section]
\newtheorem{lemma}[theorem]{Lemma}
\newtheorem{proposition}[theorem]{Proposition}

\theoremstyle{remark}
\newcommand{\R}{\mathbb R}
\newcommand{\norm}[1]{\left\lVert #1\right\rVert}
\newcommand{\dd}{\mathop{}\!\mathrm{d}}

\DeclareMathOperator{\diver}{div}

\begin{document}
\setlength{\parindent}{2em}
\setlength{\emergencystretch}{1em}
\title{Rayleigh's conjecture for clamped plates in dimensions four and five
\thanks{Research supported by NNSF of China (No. 12371110 and No. 12301133).}}

\author{
Guowei Dai\thanks{School of Mathematical Sciences,
Dalian University of Technology, Dalian 116024, P.R. China.
\newline E-mail: daiguowei@dlut.edu.cn.}
\and
Yingxin Sun\thanks{School of Mathematical Sciences,
Dalian University of Technology, Dalian 116024, P.R. China.
\newline E-mail: sunyingxin2023@mail.dlut.edu.cn.}
\and
Yong Zhang\thanks{School of Mathematical Sciences,
Jiangsu University, Zhenjiang 212013, P.R. China.
\newline E-mail: zhangyong@ujs.edu.cn.}
}
\date{}

\maketitle
\renewcommand{\abstractname}{Abstract}

\begin{abstract}
In 1877, Lord Rayleigh conjectured that the disk minimizes the
fundamental frequency of a clamped plate of prescribed area.
We prove the higher-dimensional conjecture in dimensions four and
five by establishing a sharp comparison between the first clamped-plate
eigenvalue and the square of the first Dirichlet eigenvalue on the same
domain. The proof uses a nonlinear transformation built from the
ball eigenfunctions. A local identity keeps all terms of the Hessian.
Exact coefficient inequalities and bounds for the infinite remainders
establish positivity, while a weak formulation handles sign changes
and critical sets. We characterize equality in the volume inequality
and identify an obstruction to this transformation in dimension six.
\end{abstract}



\section{Introduction}
In 1877, Lord Rayleigh \cite{Rayleigh} conjectured that, among
homogeneous clamped plates of prescribed area, the circular plate has
the least fundamental frequency. Its higher-dimensional formulation
asks whether a ball minimizes the first eigenvalue of the clamped
biharmonic operator under a volume constraint. We establish this
inequality in dimensions four and five through a sharp comparison with
the first Dirichlet eigenvalue on the same open set.

For an open set $\Omega\subset\R^n$ with $0<\left|\Omega\right|<+\infty$, define
\[
\begin{aligned}
 \lambda_1\left(\Omega\right)
 &=\inf_{v\not\equiv0,\ v\in H^1_0\left(\Omega\right)}
 \frac{\int_\Omega\left|\nabla v\right|^2\dd x}{\int_\Omega v^2\dd x},\\
 \Gamma_1\left(\Omega\right)
 &=\inf_{u\not\equiv0,\ u\in H^2_0\left(\Omega\right)}
 \frac{\int_\Omega\left(\Delta u\right)^2\dd x}{\int_\Omega u^2\dd x}.
\end{aligned}
\]
Here $\left|\Omega\right|$ is Lebesgue measure, $\nabla$ is the Euclidean gradient,
and $\Delta$ is the Laplacian. The Sobolev space $H^k\left(\Omega\right)$
consists of square-integrable functions whose weak derivatives through
order $k$ are square-integrable; $H^k_0\left(\Omega\right)$ is the closure
of $C_c^\infty\left(\Omega\right)$ in this space. The latter denotes smooth
functions with compact support in $\Omega$. All spaces are real.
On a smooth bounded domain, the plate problem is
\[
 \Delta^2u=\Gamma u\quad\text{in }\Omega,\qquad
 u=\partial_{\mathbf n}u=0\quad\text{on }\partial\Omega,
\]
where $\mathbf n$ is the outward unit normal and
$\partial_{\mathbf n}u=\nabla u\cdot\mathbf n$.
The boundary conditions fix both displacement and normal slope.
For a homogeneous isotropic thin plate with fixed material parameters
and thickness, the square of each vibration frequency is proportional
to the corresponding eigenvalue $\Gamma$ \cite{GGS}.
Thus a lower bound for $\Gamma_1$ gives a lower bound for the fundamental
frequency using only geometric data. Unlike an upper bound obtained
by inserting a single trial displacement into the Rayleigh quotient, such
a bound must control every admissible displacement. In higher
dimensions the same variational problem belongs to the theory of
fourth-order elliptic operators and spectral shape optimization
\cite{GGS,Henrot}.
Write $B_\rho=\left\{x\in\R^n:\left|x\right|<\rho\right\}$,
$\omega_n=\left|B_1\right|$, and
$B_\Omega=B_{\left(\left|\Omega\right|/\omega_n\right)^{1/n}}$.
The conjectured inequality is
$\Gamma_1\left(\Omega\right)\ge\Gamma_1\left(B_\Omega\right)$.

The analogous membrane problem was resolved by Faber and Krahn in the 1920s; see
\cite[Chapter 3]{Henrot}. Their inequality follows from symmetric
decreasing rearrangement, which replaces superlevel sets by concentric
balls of the same measure and decreases the Dirichlet energy. For a clamped plate, however, the same argument runs into two
obstacles. Rearrangement does
not in general decrease the squared-Laplacian energy, and the first
clamped eigenfunction need not be sign-definite. This second difficulty comes from the absence of a general positivity principle for the clamped
biharmonic operator \cite{GGS}. Taking an absolute value does not remove
it: across a regular nodal hypersurface, $\left|u\right|$ has a jump in its first
normal derivative and need not belong to $H^2$.

Szeg\H{o} proved a fixed-sign result in 1950 \cite{Szego}. Talenti's
comparison theorem for Poisson equations \cite{Talenti76} supplied a way
to rearrange second-order data, and his 1981 work on the clamped plate
\cite{Talenti81} produced a universal, though nonsharp, lower bound. Its
constant is about $0.98$ in the plane and tends to $1/2$ as the dimension
increases. Nadirashvili's planar proof, announced in 1992 and published
in 1995 \cite{Nadirashvili}, incorporated further constraints into this
approach. Ashbaugh and Benguria \cite{AshbaughBenguria} then gave a
Bessel-function proof in dimensions two and three. These works overcame
the sign-changing ground-state obstruction by reducing the problem to
coupled radial functions on two balls. The survey of Ashbaugh, Benguria,
and Laugesen \cite{AshbaughBenguriaLaugesen} explains both this reduction
and its relationship with the distinct buckling problem.

The obstruction in dimension four can be stated quantitatively. Let
$J_\nu$ and $I_\nu$ denote the Bessel and modified Bessel functions of
the first kind, with $\nu=n/2-1$. Let $q=j_{\nu,1}$ be the first
positive zero of $J_\nu$, and let $R$ be the first positive zero of
\[
 J_\nu\left(r\right)I_{\nu+1}\left(r\right)
 +J_{\nu+1}\left(r\right)I_\nu\left(r\right).
\]
The ball formulas are
$\lambda_1\left(B_1\right)=q^2$ and
$\Gamma_1\left(B_1\right)=R^4$
\cite{AshbaughBenguria,GGS}.
Normalize the two radii by $a^n+b^n=1$, with $0\le a\le b$.
For $a,b>0$, the relaxed quotient is
\[
 k^4\left(a\right)=
 \min_{\left(v,w\right)\ne\left(0,0\right)}
 \frac{\int_{B_a}\left(\Delta v\right)^2\dd x
       +\int_{B_b}\left(\Delta w\right)^2\dd x}
      {\int_{B_a}v^2\dd x+\int_{B_b}w^2\dd x},
\]
where $v,w$ are radial functions in
$H^2\left(B_a\right)\cap H^1_0\left(B_a\right)$ and
$H^2\left(B_b\right)\cap H^1_0\left(B_b\right)$, respectively,
subject to the single flux constraint
$a^{n-1}v'\left(a\right)=b^{n-1}w'\left(b\right)$.
Unlike a clamped function on each ball, their boundary slopes need not
vanish separately. The fourth-root parameter has the following limiting values
\cite[Section 5, p.~390]{AshbaughLaugesen}
\[
 k\left(0\right)=R,\qquad
 k\left(2^{-1/n}\right)=2^{1/n}q.
\]
For $n=2,3$, Ashbaugh and Benguria prove that the minimum occurs when
the first ball shrinks to radius $0$. In contrast, Ashbaugh and Laugesen
\cite[Section 5]{AshbaughLaugesen} prove that for $n\ge4$ it occurs
at two equal balls, and
\[
 2^{1/n}q<R,\qquad
 \Gamma_1\left(\Omega\right)
 \ge c_n\Gamma_1\left(B_\Omega\right),\qquad
 c_n=2^{4/n}\left(\frac qR\right)^4<1.
\]
Thus $c_4\simeq0.953797$ and $c_5\simeq0.921845$, whereas
$c_n\to1$ as $n\to+\infty$. This is the optimal bound obtained from the two-ball relaxation.

The relaxed class contains radial pairs that satisfy the boundary
and flux constraints but need not come from a single clamped
function. Its minimum is $c_nR^4<R^4$. Thus the relaxed problem
cannot give the conjectured constant. A proof must use additional
constraints satisfied by clamped functions, or find a different comparison.

\begin{figure}[!htbp]
\centering
\begin{tikzpicture}[x=1cm,y=1cm,font=\small,
  every path/.style={line width=.6pt}]
\node at (0,2.15) {$n=2,3$};
\node at (7.0,2.15) {$n\ge4$};
\draw[black!35] (3.45,-1.45)--(3.45,2.35);
\fill (-1.15,.65) circle (.035);
\node[below=4pt] at (-1.15,.65) {$a\downarrow0$};
\draw[fill=black!8] (.8,.65) circle (.9);
\node at (.8,.65) {$B_1$};
\node at (.8,-.5) {$b\uparrow1$};
\draw[fill=black!8] (5.9,.65) circle (.76);
\draw[fill=black!8] (8.1,.65) circle (.76);
\node at (5.9,.65) {$B_a$};
\node at (8.1,.65) {$B_b$};
\node at (7,-.5) {$a=b=2^{-1/n}$};
\node at (0,-1.1) {$\min_a k^4\left(a\right)=R^4$};
\node at (7,-1.1) {$\min_a k^4\left(a\right)=2^{4/n}q^4<R^4$};
\node[text badly centered,text width=12cm] at (3.45,-2.0)
 {$a^n+b^n=1,\qquad a^{n-1}v'\left(a\right)=b^{n-1}w'\left(b\right)$};
\end{tikzpicture}
\caption{The two-ball relaxation: a degenerate pair minimizes the
auxiliary quotient in dimensions two and three, whereas an equal pair
minimizes it for $n\ge4$ \cite{AshbaughBenguria,AshbaughLaugesen}.
The circles represent balls schematically; the point at left represents
a disappearing ball. The strict gap at right is a limitation of the
relaxation, not a counterexample to the Rayleigh conjecture.}
\label{fig:two-balls}
\end{figure}
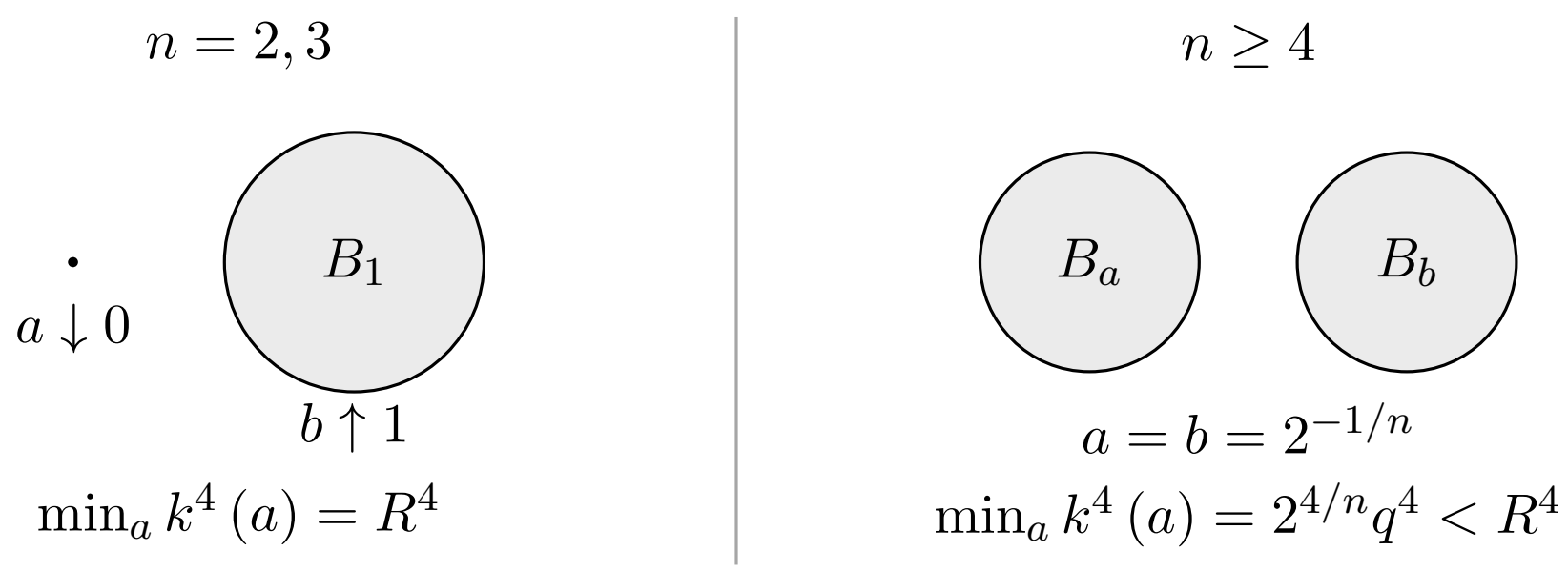

Figure~\ref{fig:two-balls} displays the minimizing configurations of
the auxiliary problem. The change at dimension four concerns this
relaxed class, whose constraints permit boundary slopes excluded by
the clamped boundary conditions.

Several later developments clarify the scope of this difficulty.
Chasman and Langford \cite{ChasmanLangford} extend the Talenti--Ashbaugh--Benguria
approach to the bi-Hermite operator in Gaussian space, where the
comparison sets are half-spaces. They prove a volume-dependent-factor
bound and report numerical evidence for a factor at least $0.91$;
this concerns the Gaussian problem. Leylekian obtains conditional ball
optimality from information about the critical values of an optimal
eigenfunction \cite{LeylekianCritical}, and from an order-reduction
principle with additional mean-value or boundary hypotheses
\cite{Leylekian}. Freitas and Leylekian extend the fixed-sign theorem
to plates under nonnegative tension \cite[Theorem A]{FreitasLeylekian}.
For a related functional, Ashbaugh, Bucur, Laugesen, and Leylekian
\cite{AshbaughBucurLaugesenLeylekian} prove the all-dimensional
fourth-order Saint--Venant inequality, which maximizes mean deflection
under a constant transverse load. More recently, Bucur and Laugesen
\cite{BucurLaugesen} obtain a simpler planar proof and nonlinear
extensions from a paired Poisson comparison. The unrestricted Euclidean $L^2$ plate inequality in $n\ge4$
is not established in these works.

Our main result is the following. Set $\alpha=q^2/R^2$.
\begin{theorem}\label{thm:main}
Let $n\in\left\{4,5\right\}$ and let $\Omega\subset\R^n$ be open with
$0<\left|\Omega\right|<+\infty$. Then
\begin{equation}\label{eq:main}
 \Gamma_1\left(\Omega\right)
 \ge\frac{R^4}{q^4}\lambda_1\left(\Omega\right)^2
 \ge R^4\left(\frac{\omega_n}{\left|\Omega\right|}\right)^{4/n}
 =\Gamma_1\left(B_\Omega\right).
\end{equation}
The first constant is attained on every ball. Equality in the final
volume comparison implies that there is a ball $B$ with
$\left|B\right|=\left|\Omega\right|$ and $\left|\Omega\triangle B\right|=0$, where $\triangle$
denotes symmetric difference. If $\Omega$ has locally Lipschitz
boundary, equality holds precisely for balls.
\end{theorem}

The proof of Theorem~\ref{thm:main}, including the equality statements,
is given in Section~\ref{sec:main-proof}.
The first inequality compares the two spectra on the same open set;
the volume inequality then follows from Faber--Krahn.
The proof applies to the variational test class and requires no
assumption on the existence or regularity of a shape optimizer. Equality on
rough open sets is described by the clamped Sobolev space in
Proposition~\ref{prop:equality}; agreement with a ball almost everywhere
alone is not sufficient.

\begin{figure}[ht]
\centering
\begin{tikzpicture}[>=stealth, node distance=1.0cm,
  every node/.style={draw, rounded corners, text centered,
    inner sep=6pt, text width=3.0cm}]
\node[draw,rounded corners,text centered,inner sep=6pt,text width=3.0cm] (plate)
  {$U(r)$ and $\phi(r)$\\radial plate and membrane profiles};
\node[draw,rounded corners, text centered,inner sep=6pt,text width=3.0cm,
  right=of plate] (ratio)
  {$s(r)=U(r)/(-U'(r))$\\strictly decreasing};
\node[draw,rounded corners, text centered,inner sep=6pt,text width=3.0cm,
  right=of ratio] (profile)
  {$f(s)=\phi(r)/(-U'(r))$\\radial comparison profile};
\node[draw,rounded corners, text centered,inner sep=6pt,text width=4.0cm,
  below=0.9cm of profile] (trial)
  {$F_u=|\nabla u|\,f\!\left(|u|/|\nabla u|\right)$\\comparison function};
\draw[->] (plate)--node[above,draw=none,text width=1.2cm] {$r\mapsto s$}(ratio);
\draw[->] (ratio)--node[above,draw=none,text width=1.2cm] {invert}(profile);
\draw[->] (profile)--node[right,draw=none,text width=1.0cm] {apply to $u$}(trial);
\end{tikzpicture}
\caption{From the radial ball profiles to the comparison function for a general trial function.}
\label{fig:radial-map}
\end{figure}
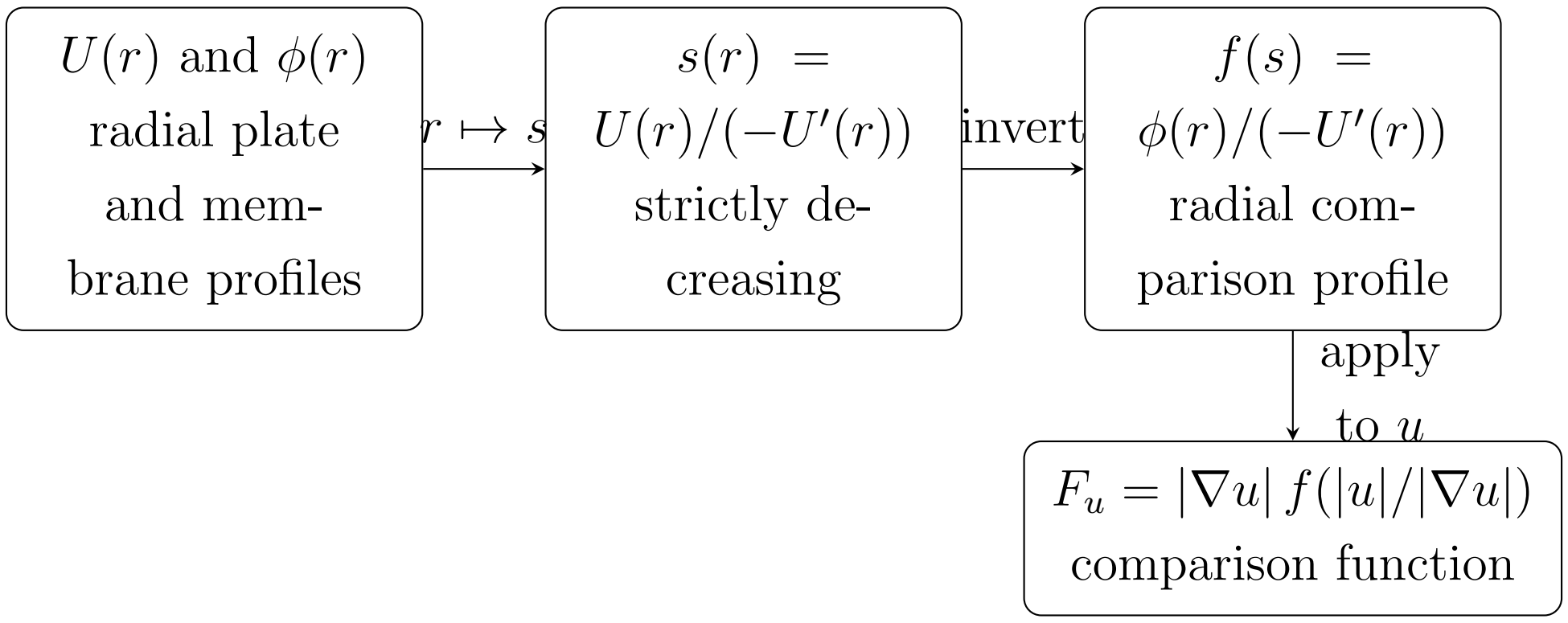

The proof constructs, for every $u\in C_c^\infty\left(\Omega\right)$,
$F_u\not\equiv0$, $F_u\in H^1_0\left(\Omega\right)$ when $u\ne0$, such that
\begin{equation}\label{eq:comparison-intro}
 \int_\Omega\left(\left|\nabla F_u\right|^2-\alpha F_u^2\right)\dd x
 \le\int_\Omega\left(\left(\Delta u\right)^2-u^2\right)\dd x.
\end{equation}
To choose the map, let $U$ and $\phi$ be the positive radial plate and
membrane profiles on $B_R$, normalized by $U''\left(R\right)=1$ and
$\phi'\left(R\right)=-1$. The ratio $s=U/\left(-U'\right)$ is a strictly monotone
function of the radius. Defining $f\left(s\right)=\phi/\left(-U'\right)$
gives
\[
 F_u=\left|\nabla u\right|f\left(\frac{\left|u\right|}{\left|\nabla u\right|}\right)
 \quad\text{where }\left|\nabla u\right|>0.
\]
The construction is summarized schematically in Figure~\ref{fig:radial-map}.
Its continuous extension is specified at critical points. For $u=U$, this definition gives $F_U=\phi$. By a calibration we mean this choice of comparison map
and divergence correction whose energy defect is nonnegative and
vanishes on the reference profile.

Radial equality alone does not imply an inequality for general
functions. The local identity uses a divergence correction involving
$\left(\Delta u\right)I_n-D^2u$, where $I_n$ is the identity matrix and
$D^2u$ is the Hessian. It completes the energy defect into a quadratic
form in the normal and tangential traces, together with the mixed and
trace-free tangential components. These terms are needed to control
arbitrary, nonradial perturbations of the reference profile. The remaining two-by-two
matrix condition reduces to positivity of three entire numerators
through its Schur complement: for a symmetric matrix
$\left(\begin{smallmatrix}a&d/2\\d/2&e\end{smallmatrix}\right)$ with $e>0$,
positive definiteness is equivalent to $a-d^2/\left(4e\right)>0$.
We cancel the factors that vanish at the ends of the radial interval.
Thirty-nine rational coefficient inequalities, together with bounds
for the full remainders, then prove positivity on the whole interval.

A second step extends the local comparison across nodal and
critical sets. Lipschitz extensions control the critical set, including
when it has positive measure. Kato's inequality supplies the nodal
measure, and the weighted contribution is nonnegative. The
coefficient inequalities and these weak extensions together yield
\eqref{eq:comparison-intro} for every $u\in C_c^\infty\left(\Omega\right)$,
including sign-changing functions. Scaling and the Dirichlet
variational principle then give the first inequality in \eqref{eq:main};
the Faber--Krahn theorem gives the second.

The restriction to dimensions four and five comes from the signs in this construction.
The coefficient inequalities for the local identity hold in dimensions
four and five. In dimension
six, the centre coefficient
$\alpha\left[\phi\left(0\right)/U\left(0\right)\right]^2-1$ is negative.
A smooth function constant on a large ball then makes the volume
contribution to the comparison defect negative, while its transition
shell contributes only a lower-order term. Consequently
\eqref{eq:comparison-intro} fails in dimension six for this particular
map, independently of the divergence correction. The obstruction
concerns this transformation; the Rayleigh inequality in dimensions
six and above is not decided by this argument.

When the present paper was in the final stages of preparation,
an independent proof of Theorem~\ref{thm:main} appeared in
\cite{Suragan2026}. While both papers use a similar collection of
ideas, their technical implementations are different.

The rest of this paper is organized as follows.
Section~\ref{sec:profiles} introduces the ball profiles and derives the
local Hessian identity. Section~\ref{sec:certificate} establishes the
coefficient inequalities. Section~\ref{sec:proof-main} proves
Theorem~\ref{thm:main}: Section~\ref{sec:weak-extension} supplies the
weak extensions, and Section~\ref{sec:main-proof} combines them with
the coefficient inequalities to obtain the integral comparison,
the spectral inequalities, and the stated equality cases. The section
ends with a further characterization of equality on arbitrary open sets.
Section~\ref{sec:six} treats the
six-dimensional obstruction. Appendix~\ref{sec:paper-input} is devoted
to the finite rational certification: its first subsection encloses the
spectral parameters, the next gives the profile-coefficient intervals,
and the subsequent subsections derive the numerator coefficients,
certify their lower bounds, and estimate the products of the finite
coefficient blocks. Appendix~\ref{sec:tail} treats the infinite parts:
it constructs positive majorants for the complete profile series,
proves the tails of the quadratic numerators, and then proves the
three determinant-factor tail bounds.

\section{Ball profiles and the local comparison}\label{sec:profiles}
We use the reference ball to build the comparison map and to determine its sharp constant. The coefficients in the Hessian identity come from the radial plate and membrane profiles.

\subsection{Reference profiles and the radial coordinate}
Standard
separation of variables on a ball gives the positive radial first
profiles below, with radial variable $r=\left|x\right|$
\cite{AshbaughBenguria,GGS}. Their values at $r=0$ are
understood as limits. We use $U''\left(R\right)=1$ and $\phi'\left(R\right)=-1$:
\[
 U\left(r\right)=\frac12\left(\frac Rr\right)^\nu
 \left(\frac{I_\nu\left(r\right)}{I_\nu\left(R\right)}-\frac{J_\nu\left(r\right)}{J_\nu\left(R\right)}\right),
 \quad
 \phi\left(r\right)=\frac Rq\left(\frac Rr\right)^\nu
 \frac{J_\nu\left(qr/R\right)}{J_{\nu+1}\left(q\right)}.
\]
To check the normalizations, denote the two terms in parentheses
in the formula for $2U$ by $F_I$ and $F_J$, respectively, including
the common factor $\left(R/r\right)^\nu$. The radial Bessel equations give
$\Delta F_I=F_I$ and $\Delta F_J=-F_J$. Hence
\[
 \Delta U=\frac{F_I+F_J}{2},\qquad \Delta^2U=U.
\]
At $r=R$ both $F_I$ and $F_J$ equal $1$. The differentiation formulas
for $r^{-\nu}I_\nu\left(r\right)$ and $r^{-\nu}J_\nu\left(r\right)$,
together with the equation defining $R$, give
\[
 F_I'\left(R\right)=\frac{I_{\nu+1}\left(R\right)}{I_\nu\left(R\right)}
 =-\frac{J_{\nu+1}\left(R\right)}{J_\nu\left(R\right)}
 =F_J'\left(R\right).
\]
Thus $U\left(R\right)=U'\left(R\right)=0$ and
$\Delta U\left(R\right)=U''\left(R\right)=1$.
Similarly, the equation $J_\nu\left(q\right)=0$ and the same
differentiation formula give $\phi\left(R\right)=0$,
$\phi'\left(R\right)=-1$, and $-\Delta\phi=\alpha\phi$.
In particular,
\[
 \eta=\left(\Delta U\right)'\left(R\right)
      =\frac{I_{\nu+1}\left(R\right)}{I_\nu\left(R\right)}>0.
\]
Here and below a prime denotes differentiation in $r$ unless another
variable is specified. We use $O\left(t^k\right)$ for a quantity
bounded in absolute value by a constant times $\left|t\right|^k$ in the stated
limit. Set
$m=n-1$, $T=R^2$, $Q=q^2$, $Z=R\eta$, and $K=T/16$.
The rational enclosures of $R,q,T,Q,Z$ used below are listed in
Appendix~\ref{subsec:appA-parameters}. Their endpoint signs are
certified there by the normalized Bessel-series calculation and its
explicit remainder bound.

Put $p=-U'>0$ and
\[
\begin{gathered}
 s=U/p,\qquad A=-p'/p,\qquad C=-1/r,\qquad
 V=\left(\Delta U\right)'/p,\\
 f=\phi/p,\qquad w=\phi'/p,\qquad z=f_s,\qquad g=f-sz.
\end{gathered}
\]
All quantities on the right initially depend on $r\in\left(0,R\right)$. A subscript
$s$ will mean differentiation after changing to the coordinate $s$.
The map $r\mapsto s\left(r\right)$ is a strictly decreasing bijection from
$\left(0,R\right)$ onto $\left(0,+\infty\right)$, with $s\left(0+\right)=+\infty$ and $s\left(R-\right)=0$.
For the strict monotonicity, $p$ is a positive constant times
\[
 r^{-\nu}J_{\nu+1}\left(r\right)
 \left(1-\left|c\right|\frac{I_{\nu+1}\left(r\right)}{J_{\nu+1}\left(r\right)}\right),
 \qquad c=J_\nu\left(R\right)/I_\nu\left(R\right)\in\left(-1,0\right).
\]
Here $c$ is a fixed scalar, and $j_{\mu,k}$ denotes the $k$th
positive zero of $J_\mu$ for the orders $\mu$ used below. The canonical
product for $J_\mu$ and the
order monotonicity of its positive zeros are standard
\cite[Section 10.21, especially (10.21.15)]{DLMF}.
Recall that a positive function is log concave when its logarithm is
concave. The first positive zero $R$ of the Bessel cross-product satisfies
$R<j_{\nu+1,1}$; this comparison is recorded in
Ashbaugh--Laugesen~\cite{AshbaughLaugesen}, based on the monotonicity
theorem of Lorch~\cite{Lorch}. The canonical product and the ordering
of the zeros of $J_{\nu+1}$ are given in
\cite[Section~10.21]{DLMF}. These facts show that
$r^{-\nu}J_{\nu+1}\left(r\right)$ is strictly log concave: its logarithm is a
constant plus $\log r+\sum_{k\ge1}\log\left(1-r^2/j_{\nu+1,k}^2\right)$.
For each positive zero $j=j_{\nu+1,k}$, the derivatives of the
corresponding logarithmic factor satisfy
\[
 \frac{\mathrm d^2}{\mathrm d r^2}\log\left(1-r^2/j^2\right)
 =-\frac{2\left(j^2+r^2\right)}{\left(j^2-r^2\right)^2}<0.
\]
The quotient $I_{\nu+1}/J_{\nu+1}$ has the product representation
\[
 \frac{I_{\nu+1}\left(r\right)}{J_{\nu+1}\left(r\right)}
 =\prod_{k\ge1}\frac{1+r^2/j_{\nu+1,k}^2}
                         {1-r^2/j_{\nu+1,k}^2}.
\]
Every factor is positive, increasing, and convex for $0<r<R$;
products preserve these properties, as does their locally uniform
limit. The root relation at $R$ makes
$1-\left|c\right|I_{\nu+1}/J_{\nu+1}$ vanish at $R$.
This factor is therefore positive and concave on $\left(0,R\right)$.
The logarithm of a positive concave function is concave, because
$\left(\log h\right)''=h''/h-\left(h'/h\right)^2\le0$.
Adding it to the strictly concave logarithm of the first factor
proves that $\log p$ is strictly concave.
Since $U\left(r\right)=\int_r^R p\left(t\right)\dd t$, we have
$s'=-1-Up'/p^2$. If $p'\ge0$, this is negative. If $p'<0$, strict
log concavity gives
\[
 U\left(r\right)<p\left(r\right)\int_r^{+\infty}
 \exp\left(\left(t-r\right)p'\left(r\right)/p\left(r\right)\right)\dd t
 =-\frac{p\left(r\right)^2}{p'\left(r\right)},
\]
Indeed, concavity bounds $\log p\left(t\right)$ above by its tangent
line at $r$ for $r<t<R$; extending the resulting exponential integral
to infinity gives the displayed strict bound. Multiplication by
$-p'\left(r\right)/p\left(r\right)^2>0$ gives
$-Up'/p^2<1$, and hence again $s'<0$.
At the centre the even profile series gives
$U\left(r\right)=U\left(0\right)+O\left(r^2\right)$ and
$p\left(r\right)=-U''\left(0\right)r+O\left(r^3\right)$, with
$U\left(0\right)>0$ and $U''\left(0\right)<0$.
At the boundary, Taylor's formula and the normalization give
\[
 U\left(r\right)=\frac12\left(R-r\right)^2+O\left(\left(R-r\right)^3\right),
 \qquad p\left(r\right)=R-r+O\left(\left(R-r\right)^2\right).
\]
Consequently $s\left(r\right)$ tends to infinity at the centre and
equals $\left(R-r\right)/2+O\left(\left(R-r\right)^2\right)$ near the
boundary. Continuity and strict monotonicity now give the asserted
bijection and justify differentiation with respect to $s$.

\subsection{The full Hessian identity}
For $P=\left|\nabla u\right|>0$, the map $F=P f\left(u/P\right)$ is fixed by
equality on the ball.
To control its energy for a general function, we introduce divergence
terms that vanish after integration. For a vector field
$J=\left(J^1,\ldots,J^n\right)$, its divergence is
$\diver J=\sum_{i=1}^n\partial_iJ^i$.
The coefficients below are chosen so that the radial Hessian of $U$
is a zero of the resulting quadratic form.

We next define the two correction coefficients $G$ and $b$:
\[
\begin{gathered}
 V_0=w^2-mCfw+\alpha f^2,\qquad
 \Psi=A'+s^2+\left(1+msC\right)V,\\
 G=\frac{V_0-\Psi}{mC^2},\qquad
 b=-\left(A+sV\right)-mGC-fw.
\end{gathered}
\]
Define also $a=1-g^2$, $d=2G-sG_s$, and
$e=\left[1+\left(m-1\right)G\right]/m$. The letter $d$ here is a scalar coefficient;
differentials in integrals are written $\dd$.
The differential relations used below can all be read from the
two radial equations. Indeed,
\[
 U''=Ap,\qquad \Delta U=\left(A+mC\right)p,\qquad
 p'=-Ap,\qquad \left(\Delta U\right)'=Vp.
\]
Differentiating the second identity and dividing by $p$ gives
$V=A'+mC'-A\left(A+mC\right)$.
The equation $\Delta\left(\Delta U\right)=U$ gives
$V'-AV-mCV=s$. Similarly, differentiating
$\phi=p f$ and $\phi'=pw$ and using $\Delta\phi=-\alpha\phi$
gives the two equations for $f'$ and $w'$ below. Thus
\[
\begin{gathered}
 s'=sA-1,\quad C'=C^2,\quad
 A'=A^2+mAC-mC^2+V,\quad V'=\left(A+mC\right)V+s,\\
 f'=w+Af,\quad w'=\left(A+mC\right)w-\alpha f,\quad w=-z-Ag.
\end{gathered}
\]
The last identity follows by writing
$f'=z s'=z\left(sA-1\right)=w+Af$ and using $g=f-sz$.
For clarity, the derivatives of the correction data before collecting
terms are
\[
\begin{aligned}
 V_0'&=2ww'-mC'fw-mCf'w-mCfw'+2\alpha ff',\\
 \Psi'&=A''+2ss'+m\left(s'C+sC'\right)V
                    +\left(1+msC\right)V',\\
 A''&=2AA'+mA'C+mAC'-2mCC'+V'.
\end{aligned}
\]
Replacing the derivatives on the right by the preceding equations,
and grouping the coefficients of $w^2,fw,f^2$ in the first line,
yields
\[
 V_0'=\left(2A+mC\right)V_0-mC^2fw,\qquad
 \Psi'-\left(2A+mC\right)\Psi=mC^2\left(A-2C+sV\right).
\]
In differentiating $G$, its denominator also contributes:
$G'=\left(V_0'-\Psi'\right)/\left(mC^2\right)-2CG$.
Differentiating the definition of $b$ gives
\[
 b'=-A'-s'V-sV'-mCG'-mC'G-f'w-fw'.
\]
Substitution in these two formulas yields
\[
 G'=\left(2A+\left(m-2\right)C\right)G-fw-\left(A-2C+sV\right),\quad
 b'+mCG'=-2f'w-2s'V.
\]
These relations yield the three identities needed to complete the
quadratic expression:
\[
\begin{split}
2b-sb_s-2zg&=-2aA-mCd,\\
b+G_s&=-dA-2meC,\\
-s^2-z^2+\alpha f^2+b_s&=aA^2+mdAC+m^2eC^2.
\end{split}
\]
Here division by $s'$ is legitimate because $s'<0$.
In particular, the preceding identity for $b'+mCG'$ becomes
\[
 b_s=-2zw-2V-mCG_s.
\]
For the first completion identity, use this expression and then
$w+z=-Ag$ to obtain
\[
\begin{aligned}
 2b-sb_s-2zg
 &=-2A-mC\left(2G-sG_s\right)-2g\left(w+z\right)\\
 &=-2A-mCd+2Ag^2=-2aA-mCd.
\end{aligned}
\]
For the second, $s'G_s=G'$ gives
\[
\begin{aligned}
 b+G_s+dA+2meC
 &=b-G'+2AG+2C+2\left(m-1\right)CG=0,
\end{aligned}
\]
where the last equality is the displayed formula for $G'$.
Finally, $-z^2-2zw=w^2-\left(w+z\right)^2=w^2-A^2g^2$.
Using this equality and the second completion identity, the difference
between the two sides of the third is
\[
\begin{aligned}
 &-s^2-z^2+\alpha f^2+b_s
          -\left(aA^2+mdAC+m^2eC^2\right)\\
 &\quad=w^2+\alpha f^2-s^2-A^2-2V+mCb+m^2eC^2\\
 &\quad=V_0-\left[A^2+mAC-mC^2+s^2
                         +\left(2+msC\right)V\right]-mC^2G\\
 &\quad=V_0-\Psi-mC^2G=0.
\end{aligned}
\]
Thus the completion uses the radial equations and the definitions of
$G,b$, without imposing any differential equation on the test function.

The following identity measures the defect from the ball profile in
every Hessian component; its sign will be established in the next section.
\begin{proposition}\label{prop:identity}
For a positive smooth function $u$ at a point where
$P=\left|\nabla u\right|>0$, set
\[
 F=P f\left(u/P\right),\quad J_1=P b\left(u/P\right)\nabla u,
 \quad J_2=G\left(u/P\right)\left[\left(\Delta u\right)I_n-D^2u\right]\nabla u.
\]
All scalar coefficients below are evaluated at $s=u/P$. In an orthonormal
frame whose first vector is $\nabla u/P$, write
\[
 D^2u=\begin{pmatrix}X&\zeta^{\mathsf T}\\
 \zeta&T_0+\tau I_m/m\end{pmatrix},\qquad
\operatorname{tr}T_0=0.
\]
Thus $X$ is the normal entry, $\zeta\in\R^m$ the mixed vector, and
$\tau$ the trace (sum of diagonal entries) of the tangential block.
For matrices, $\left|\cdot\right|$
denotes the square root of the sum of the squares of their entries.
Then
\begin{equation}\label{eq:hessian}
\begin{split}
&\left|D^2u\right|^2-u^2-\left|\nabla F\right|^2+\alpha F^2+\operatorname{div}\left(J_1+J_2\right)\\
={}&a\left(X-AP\right)^2+d\left(X-AP\right)\left(\tau-mCP\right)+e\left(\tau-mCP\right)^2\\
&+\left(2-g^2-d\right)\left|\zeta\right|^2+\left(1-G\right)\left|T_0\right|^2.
\end{split}
\end{equation}
\end{proposition}
\begin{proof}
Write $u_{ij}=\partial_i\partial_j u$. The matrix field
$N=\left(\Delta u\right)I_n-D^2u$, sometimes called the first
Newton tensor of the Hessian, has zero divergence, since
$\sum_i\partial_iN_{ij}=\partial_j\Delta u-\sum_i\partial_i u_{ij}=0$.
The orthogonal splitting of the tangential block gives
\[
 \left|D^2u\right|^2=X^2+2\left|\zeta\right|^2
                       +\left|T_0\right|^2+\tau^2/m,
 \qquad \Delta u=X+\tau.
\]
Since $\partial_iP=P^{-1}\sum_j u_j u_{ji}$,
$\nabla P=\left(X,\zeta\right)$ in the chosen frame. The chain rule gives
\[
 \nabla\left(u/P\right)=\frac{\nabla u-s\nabla P}{P},\qquad
 \nabla F=z\nabla u+g\nabla P,
\]
and hence
\[
 \left|\nabla F\right|^2=\left(zP+gX\right)^2+g^2\left|\zeta\right|^2.
\]
For $J_1$, the derivative of its scalar factor is
$\nabla\left(Pb\right)=b_s\nabla u+\left(b-sb_s\right)\nabla P$.
For $J_2$, the zero divergence of $N$ implies
$\diver\left(N\nabla u\right)=\operatorname{tr}\left(ND^2u\right)$,
and $N\nabla u=P\left(\tau,-\zeta\right)$.
These identities give
\[
\begin{aligned}
 \diver J_1={}&bP\left(X+\tau\right)
 +\left(b-sb_s\right)PX+b_sP^2,\\
 \diver J_2={}&G\left(2X\tau+\left(1-1/m\right)\tau^2
                  -2\left|\zeta\right|^2-\left|T_0\right|^2\right)\\
 &+G_s\left(P\tau-sX\tau+s\left|\zeta\right|^2\right).
\end{aligned}
\]
Here the second identity uses
\[
\begin{aligned}
 \operatorname{tr}\left(ND^2u\right)
   &=\left(\Delta u\right)^2-\left|D^2u\right|^2,\\
 \nabla\left(u/P\right)
   &=\left(1-sX/P,-s\zeta/P\right).
\end{aligned}
\]
The left-hand side of \eqref{eq:hessian}, before completing squares,
is therefore
\[
\begin{aligned}
 &aX^2+dX\tau+e\tau^2
       +\left(2-g^2-d\right)\left|\zeta\right|^2
       +\left(1-G\right)\left|T_0\right|^2\\
 &\quad+\left(2b-sb_s-2zg\right)PX
       +\left(b+G_s\right)P\tau
       +\left(-s^2-z^2+\alpha f^2+b_s\right)P^2.
\end{aligned}
\]
The three completion identities replace the last line by
\[
 -\left(2aA+mCd\right)PX-\left(dA+2meC\right)P\tau
       +\left(aA^2+mdAC+m^2eC^2\right)P^2.
\]
Expanding the first three terms on the right-hand side of
\eqref{eq:hessian} gives exactly this polynomial. The remaining two
Hessian terms already have the stated coefficients, which proves the
identity.
\end{proof}

The equality profile explains the geometry of this completion.
For $u=U$, the normal and tangential Hessian entries are
$X=U''=AP$ and $\tau=mU'/r=mCP$, while $\zeta=0$ and $T_0=0$.
Thus every term on the right of \eqref{eq:hessian} vanishes.
For a general function, the mixed vector and the trace-free tangential
block measure departures from this radial structure; see
Figure~\ref{fig:hessian}. Controlling only $X$ and $\tau$ would leave
these departures untreated.

\begin{figure}[!htbp]
\centering
\begin{tikzpicture}[font=\small,>=latex,line width=.65pt]
\draw[black!35] (-2,-1.1).. controls (-.8,.15) and (.3,.15).. (1.7,-1.1);
\draw (-2,-.55).. controls (-.8,.7) and (.3,.7).. (1.7,-.55);
\draw[black!35] (-2,0).. controls (-.8,1.25) and (.3,1.25).. (1.7,0);
\fill (-.15,.385) circle (.035);
\draw[->] (-.15,.385)--(-.15,1.9) node[above] {$e_1=\nabla u/P$};
\draw[->] (-.15,.385)--(1.7,.385) node[right] {$e_2$};
\node[anchor=north,text badly centered,text width=5cm] at (-.15,-1.35)
 {Level sets of $u$\\normal and tangent};
\node[draw,rounded corners=2pt,fill=black!3,inner sep=10pt,text badly centered,text width=5.5cm]
 at (6.4,.45)
 {$D^2u=\begin{pmatrix}X&\zeta^{\mathsf T}\\
 \zeta&T_0+\dfrac{\tau}{m}I_m\end{pmatrix}$\\[9pt]
 $\operatorname{tr}T_0=0,\qquad m=n-1$};
\node[text badly centered,text width=5.5cm] at (6.4,-1.9)
 {Ball profile: $X=AP$, $\tau=mCP$\\
 $\zeta=0$, $T_0=0$};
\end{tikzpicture}
\caption{The Hessian decomposition at a regular point, shown in a
two-dimensional section. The vector $e_2$ represents a tangential direction. The mixed entries $\zeta$ and the trace-free
tangential block $T_0$ are included in the local identity, so its sign
controls nonradial as well as radial test functions.}
\label{fig:hessian}
\end{figure}

It remains to establish, throughout $0<r<R$,
\begin{equation}\label{eq:positivity}
 G<0,\qquad e>0,\qquad ae-d^2/4>0.
\end{equation}
Since $ae=\left(ae-d^2/4\right)+d^2/4>0$ and $e>0$, these conditions
also give $a>0$. Indeed, for real numbers $\xi,\theta$, the exact completion
\[
 a\xi^2+d\xi\theta+e\theta^2
 =e\left(\theta+\frac{d}{2e}\xi\right)^2
       +\frac{ae-d^2/4}{e}\xi^2
\]
shows that the matrix
$\begin{pmatrix}a&d/2\\d/2&e\end{pmatrix}$ is positive definite:
its quadratic form is positive for every nonzero vector.
Apply this formula to $\xi=X-AP$ and $\theta=\tau-mCP$.
They also control the other Hessian terms: $G<0$ implies $e<1$;
$d<2\sqrt{ae}<2\sqrt a$ then gives
$2-g^2-d=1+a-d>\left(1-\sqrt a\right)^2\ge0$, and $1-G>0$.
The right-hand side in Proposition~\ref{prop:identity} is therefore nonnegative.

\section{Positivity of the comparison}\label{sec:certificate}
The next step replaces the apparently singular coefficients in
\eqref{eq:positivity} by numerators with absolutely convergent power
series. Here an entire function means a complex-valued function that
is holomorphic on the whole complex plane. It therefore has a Taylor
series about $y=0$ with infinite radius of convergence. We use this
property for the numerators below: after the removable powers of $y$
are divided out, their Taylor series still converge absolutely on
every bounded interval, in particular on $0\le y\le1$.

\subsection{Entire numerators and coefficient recurrences}
Put $y=1-r^2/T$, $x=1-y$, $u\left(y\right)=8U\left(r\right)/T$, and $v\left(y\right)=2\phi\left(r\right)/R$.
Within this section and the appendices,
$u,v$ denote only these auxiliary profiles, not
a variational trial function. A subscript $x$ means $-\partial_y$.
Thus $x=r^2/T$ increases from zero at the centre to $1$ at the
boundary, whereas $y$ increases in the opposite direction. The even
Bessel series show that $u$ and $v$ are entire functions of $y$.
Their boundary conditions give double and simple zeros at $y=0$,
respectively. Define the removable quotients
\[
 h=u/y^2,\ k=v/y,\ H=2h-yh_x.
\]
In particular $u_x=-yH$, since $\partial_x y=-1$. The remaining
numerators used to remove denominators are
\[
 D=2x\left\{y\left(h_xH-hH_x\right)-hH\right\}-yhH,\quad
 N_g=y\left(kh_x-hk_x\right)-kh,
\]
\[
\begin{split}
B={}&u_x^2-\left(2n+6\right)xu_xu_{xx}+4x^2u_{xx}^2-8x^2u_xu_{xxx}\\
&-m\left(n+2\right)uu_{xx}-2mxuu_{xxx},\\
N_G={}&4xTv_x^2+2mTvv_x+QT v^2-K\left(T^2u^2+4B\right),\\
L={}&N_G/y^3,\qquad N_a=\left(KD^2-TxN_g^2\right)/y^2,\\
N_d={}&LD-xh\left\{y\left(L_xH-2LH_x\right)-LH\right\},\\
N_e={}&4mKH^2+\left(m-1\right)yL,\qquad
N_{\det}=4H^2N_aN_e-N_d^2.
\end{split}
\]
The auxiliary profile amplitudes differ from the physical amplitudes.
Correspondingly, set $\widehat f=\sqrt K f$, $\widehat g=\sqrt K g$,
$\widehat G=KG$, $\widehat a=Ka$, $\widehat d=Kd$, and
$\widehat e=Ke$. In particular $\widehat a=K-\widehat g^2$.
We claim that direct substitution into the definitions of the normalized
coefficients yields the following identities:
\[
\widehat G=\frac{yL}{4mH^2},\quad
\widehat a=\frac{y^2N_a}{D^2},\quad
\widehat d=\frac{yN_d}{2mH^2D},\quad
\widehat e=\frac{N_e}{4m^2H^2},
\]
\[
\widehat a\widehat e-\widehat d^2/4
=\frac{y^2N_{\det}}{16m^2H^4D^2}.
\]
We now verify these formulas from the definitions. This verification is
included because the cancellations in the numerators are essential; the
displayed identities are not being used as an unexplained consequence of
substitution.
The relation $u_x=-yH$ and the normalization $U=Tu/8$ give
\[
 p=-\frac{T}{8}\frac{2r}{T}u_x=\frac{ryH}{4}.
\]
Thus $H>0$ for $0<x<1$. The formula for $s_x$ derived below,
together with $\mathrm d x/\mathrm d r=2r/T>0$ and $s'<0$,
then gives $D<0$ on the same interval. It suffices to prove $-L,N_e,N_{\det}>0$. Indeed,
\begin{equation}\label{eq:schur-numerators}
 \widehat a-\frac{\widehat d^{\,2}}{4\widehat e}
 =\frac{y^2N_{\det}}{4H^2D^2N_e}.
\end{equation}
Thus $N_e>0$ and $N_{\det}>0$ imply $\widehat a>0$ and hence $a>0$.
Equivalently, $4H^2N_aN_e=N_{\det}+N_d^2>0$.
The physical determinant is the displayed auxiliary determinant
divided by $K^2>0$.
For the reconstruction formulas, $s=Ryh/\left(2\sqrt{x}H\right)$ gives
\[
s_x=\frac{RD}{4x^{3/2}H^2},\qquad
\widehat g=\frac{R\sqrt{x}N_g}{D},\qquad
\widehat G_x=\frac{y\left(L_xH-2LH_x\right)-LH}{4mH^3}.
\]
To see the cancellation in $\widehat g$, write
$\widehat f=Rk/\left(2\sqrt{x}H\right)$. Differentiation gives
\[
 \widehat f_x=\frac{R}{4x^{3/2}H^2}
               \left[2x\left(k_xH-kH_x\right)-kH\right].
\]
Using $\widehat g=\widehat f-s\widehat f_x/s_x$, its numerator simplifies by
\[
 kD-yh\left[2x\left(k_xH-kH_x\right)-kH\right]=2xH N_g.
\]
This proves the displayed formula for $\widehat g$ and consequently
that for $\widehat a=K-\widehat g^2$.
Moreover, $\widehat G_s=\widehat G_x/s_x$, so that
\[
 s\widehat G_s
 =\frac{xyh}{2mH^2D}
           \left[y\left(L_xH-2LH_x\right)-LH\right].
\]
Subtracting this expression from $2\widehat G$ proves the formula for
$\widehat d$. The remaining coefficient is
\[
 \widehat e=\frac{K+\left(m-1\right)\widehat G}{m}.
\]
Inserting these three fractions over a common denominator gives the
determinant formula.
To express $\Psi$ in the variable $y$, put $\widetilde p=-\partial_r u$, with $u$ here regarded as
the composite profile $u\left(1-r^2/T\right)$. The quantities $s,A,V$
are unchanged by the constant rescaling from $U$ to $u$, so
\[
 \widetilde p^{\,2}\Psi
 =u_{rr}^{\,2}+\widetilde p\,u_{rrr}+u^2
   +\left(\widetilde p-mu/r\right)\left(\Delta u\right)'.
\]
Use
\[
\begin{gathered}
 u_r=\frac{2r}{T}u_x,\qquad
 u_{rr}=\frac2T\left(u_x+2xu_{xx}\right),\\
 u_{rrr}=\frac{4r}{T^2}\left(3u_{xx}+2xu_{xxx}\right),\qquad
 \left(\Delta u\right)'=
 \frac{4r}{T^2}\left(\left(n+2\right)u_{xx}+2xu_{xxx}\right).
\end{gathered}
\]
Substitution gives the polynomial $B$ above. For
$\widehat V_0=K V_0$, the analogous computation gives
\[
4u_x^2r^2\Psi=T^2u^2+4B,\quad
4u_x^2r^2\widehat V_0=4xTv_x^2+2mTvv_x+QT v^2.
\]
Multiply the first equality by $K$ and subtract it from the second.
The result is
\[
 N_G=4u_x^2r^2 K\left(V_0-\Psi\right)
     =4m u_x^2\widehat G=4my^2H^2\widehat G,
\]
where $C^2=1/r^2$ was used in the second equality. This proves
$\widehat G=yL/\left(4mH^2\right)$, including its normalization factor.

For a power series $A\left(y\right)$, write $\left[A\right]_j$ for the coefficient of $y^j$.
Let $u_j=\left[u\right]_j$, $v_j=\left[v\right]_j$, and write $\gamma=n/2$.
The profile equations in this variable are
\[
 \mathcal L^2u=\frac{T^2}{16}u,\qquad
 \mathcal Lv=-\frac Q4v,\qquad
 \mathcal L=\left(1-y\right)\partial_y^2-\gamma\partial_y.
\]
For example, the coefficient of $y^j$ in $\mathcal Lu$ is
$\left(j+1\right)\left(j+2\right)u_{j+2}-\left(j+1\right)\left(j+\gamma\right)u_{j+1}$.
Applying the same coefficient formula a second time gives
\[
\begin{aligned}
 \left[\mathcal L^2u\right]_j
 ={}&\left(j+1\right)\left(j+2\right)\left(j+3\right)\left(j+4\right)u_{j+4}\\
 &-2\left(j+1\right)\left(j+2\right)\left(j+3\right)
                                      \left(j+\gamma+1\right)u_{j+3}\\
 &+\left(j+1\right)\left(j+2\right)\left(j+\gamma\right)
                                      \left(j+\gamma+1\right)u_{j+2}.
\end{aligned}
\]
Equating this expression to $T^2u_j/16$ and dividing by the
positive coefficient of $u_{j+4}$ gives, for $j\ge0$,
\[
u_{j+4}=\frac{2\left(j+\gamma+1\right)}{j+4}u_{j+3}
-\frac{\left(j+\gamma\right)\left(j+\gamma+1\right)}{\left(j+3\right)\left(j+4\right)}u_{j+2}
+\frac{T^2u_j}{16\left(j+1\right)\left(j+2\right)\left(j+3\right)\left(j+4\right)},
\]
\[
 v_{j+2}=\frac{j+\gamma}{j+2}v_{j+1}-\frac{Qv_j}{4\left(j+1\right)\left(j+2\right)},
\]
The boundary data give $u_0=u_1=v_0=0$ and $u_2=v_1=1$.
For the next plate coefficient, differentiate
$\Delta U=U''+mU'/r$ at $R$:
\[
 U'''\left(R\right)=\eta-m/R=\left(Z-m\right)/R.
\]
On the other hand, $U=Tu/8$ and the boundary series gives
$\left.u_{xx}\right|_{y=0}=2$ and
$\left.u_{xxx}\right|_{y=0}=-6u_3$.
The radial chain rule displayed above then gives
$U'''\left(R\right)=\left(3-6u_3\right)/R$.
Equating these expressions yields $u_3=\left(n+2-Z\right)/6$.
The two recurrences now determine every later coefficient. In particular
\[
u_4=\frac{\left(n+2\right)\left(n+4-2Z\right)}{48},\quad v_2=n/4,
\quad v_3=\frac{n\left(n+2\right)-Q}{24}.
\]
These give $B=16+8my+2nmy^2+O\left(y^3\right)$ and
$4xv_x^2+2mvv_x+Qv^2=4+2my+\left(nm/2\right)y^2+O\left(y^3\right)$.
Since $u^2=O\left(y^4\right)$ and $K=T/16$, these two expansions give
\[
\begin{aligned}
 N_G={}&T\left[4+2my+\frac{nm}{2}y^2\right]
       -\frac{T}{16}\left[64+32my+8nmy^2\right]+O\left(y^3\right)\\
     ={}&O\left(y^3\right).
\end{aligned}
\]
Thus its coefficients of degrees $0$, $1$, and $2$ vanish exactly.
For the other division, put $h_1=\left[h\right]_1=u_3$.
Direct differentiation of $D,N_g$ gives
$D_0=-4$, $\left(N_g\right)_0=-1$, $D_1=2-8h_1$, and
$\left(N_g\right)_1=-2h_1$. Hence
\[
\begin{aligned}
 KD^2&=T\left[1+\left(4h_1-1\right)y\right]+O\left(y^2\right),\\
 TxN_g^2&=T\left[1+\left(4h_1-1\right)y\right]+O\left(y^2\right).
\end{aligned}
\]
The constant and linear coefficients therefore cancel as well.
Division by $y^3$ in $L$ and by $y^2$ in $N_a$ removes zeros of entire
functions; it introduces no singularity at $y=0$. All the numerators
defined above are consequently entire, and their series converge
absolutely throughout $0\le y\le1$.

An equivalent realization follows directly from the two equations
$\Delta U=W$ and $\Delta W=U$. With $W\left(y\right)=\Delta U\left(r\right)$,
write $W_j=\left[W\right]_j$. Then
\[
\begin{gathered}
 \mathcal Lu=2W,\qquad \mathcal LW=\frac{T^2}{32}u,\qquad
 W_0=1,\quad W_1=-Z/2,\\
 u_{j+2}=\frac{j+\gamma}{j+2}u_{j+1}
             +\frac{2W_j}{\left(j+1\right)\left(j+2\right)},\\
 W_{j+2}=\frac{j+\gamma}{j+2}W_{j+1}
             +\frac{T^2u_j}{32\left(j+1\right)\left(j+2\right)}.
\end{gathered}
\]
Eliminating $W$ gives the fourth-order scalar recurrence.

\subsection{Finite coefficients and complete-series positivity}
The following table is a convenient summary of the finite coefficient
bounds used in the positivity argument. Its entries are rounded down
from the rational enclosures in Appendix~\ref{subsec:appA-positive}:
the bounds for $-L$ and $N_e$ come from the first coefficient table
there, while the determinant entries of degrees zero through five come
from the direct-convolution table. For degrees six through twenty,
the lower bound $1/20$ is supplied by
Proposition~\ref{prop:finite-determinant-coefficients}. Thus the table
collects previously certified estimates; it is not an independent
numerical computation.
For
$A\in\left\{-L,N_e,N_{\det}\right\}$ the required common lower bounds are
\begingroup
\setlength{\tabcolsep}{\arraycolsep}
\setlength{\LTpre}{\medskipamount}
\setlength{\LTpost}{\medskipamount}
\begin{longtable}{>{$}c<{$}|>{$}r<{$}>{$}r<{$}>{$}r<{$}}
j&-L&N_e&N_{\det}\\
\hline
\endfirsthead
j&-L&N_e&N_{\det}\\
\hline
\endhead
\hline
\endfoot
\hline
\endlastfoot
0&12&63&3500\\1&21&62&22000\\2&12&37&62000\\
3&5&20&109000\\4&1&8&145000\\5&1/3&2&155000\\
6&1/20&7/10&1/20\\7&1/100&3/20&1/20\\8&1/1000&1/50&1/20\\
\end{longtable}
\endgroup
For $9\le j\le20$ we need only $\left[N_{\det}\right]_j>1/20$.
For $-L$ and $N_e$, the entries in the main-text table are the
downward-rounded versions of the first table in
Appendix~\ref{subsec:appA-positive}. For $N_{\det}$, degrees zero
through five are obtained by the direct convolution formula
\[
 [N_{\det}]_j=\sum_{i=0}^j(c_i e_{j-i}-z_i z_{j-i}),
 \]
whose certified values are listed in the next table of that appendix;
degrees six through twenty are covered by
Proposition~\ref{prop:finite-determinant-coefficients}. The profile
intervals and coefficient formulas used in all these computations are
given in Appendix~\ref{subsec:appA-profiles} and
Appendix~\ref{subsec:appA-formulas}, respectively.

For example, in dimension four at $j=20$, the selected positive
convolution term is $c_{15}^-e_{4,5}^-=0.138\cdot2.8$, where the
two lower bounds are entries in the coefficient tables of
Appendix~\ref{subsec:appA-positive}. The subtracted square is bounded
by $\mathcal Q_{20}=0.24477811<0.245$ using the displayed $d_i^+$
entries, and the omitted high-index contribution is bounded by
$\delta_4U_4=10^{-5}\cdot6000=0.06$. Hence
\[
 [N_{\det}]_{20}>
 0.138\cdot2.8-0.245-0.06=0.0814>\frac1{20}.
\]
The full substitution for $\mathcal Q_{20}$ is recorded in
Appendix~\ref{subsec:appA-positive}.

For the tail estimates write
\[
 \mathcal C=4H^2N_a,\qquad \mathcal E=N_e,\qquad
 \mathcal Z=N_d,\qquad N_{\det}=\mathcal C\mathcal E-\mathcal Z^2.
\]
These symbols denote entire functions of $y$; in particular
$\mathcal C$ is distinct from the radial coefficient $C=-1/r$.
For a series $A=\sum_{j\ge0}A_jy^j$ with absolutely summable
coefficients define
\[
 A^\#\left(y\right)=\sum_{j\ge0}\left|A_j\right|y^j,\qquad
 \norm{A}_{\mathcal A}=\sum_{j\ge0}\left|A_j\right|,\qquad
 \mathcal T_k A^\#=\sum_{j\ge k}\left|A_j\right|.
\]
Here $\norm{A}_{\mathcal A}$ is the coefficient-sum norm.
The absolute convergence permits rearrangement of the nonnegative
double sum, so the convolution formula yields
\[
 \norm{AB}_{\mathcal A}
 \le\sum_{j\ge0}\sum_{i=0}^j\left|A_i\right|\left|B_{j-i}\right|
 =\left(\sum_{i\ge0}\left|A_i\right|\right)
    \left(\sum_{k\ge0}\left|B_k\right|\right)
 =\norm{A}_{\mathcal A}\norm{B}_{\mathcal A}.
\]
Proposition~\ref{prop:quadratic-complete-tails} in
Appendix~\ref{subsec:appB-quadratic} gives the two quadratic
remainder bounds, while Proposition~\ref{prop:determinant-complete-tails}
in Appendix~\ref{subsec:appB-determinant} gives the three
higher-order bounds:
\[
\begin{gathered}
 \mathcal T_9L^\#<13/4,\qquad \mathcal T_9N_e^\#<61,\\
 \mathcal T_{21}\mathcal C^\#<7/4,\qquad
 \mathcal T_{21}\mathcal E^\#<10^{-9},\qquad
 \mathcal T_{21}\mathcal Z^\#<1/100.
\end{gathered}
\]

Let $\mathcal C_0,\mathcal E_0,\mathcal Z_0$ be the polynomials
obtained by retaining coefficients zero through twenty in
$\mathcal C,\mathcal E,\mathcal Z$, respectively. The upper-bound columns $c_i^+,e_i^+,d_i^+$ in the two
dimension tables of Appendix~\ref{subsec:appA-positive} give the
finite coefficient sums. The five-dimensional columns dominate the
four-dimensional ones. Their sums are
\[
 \begin{aligned}
 \sum_{i=0}^{20}c_i^+&=21318.20648<22000,\\
 \sum_{i=0}^{20}e_i^+&=305.9385,\\
 \sum_{i=0}^{20}d_i^+&=1675.432427.
 \end{aligned}
\]
For $\mathcal C_0$, only these coefficients through degree twenty
occur, so
\[
 \norm{\mathcal C_0}_{\mathcal A}
 \le\sum_{i=0}^{20}c_i^+<22000.
\]
For the complete series $\mathcal E$ and $\mathcal Z$, we must also
include the coefficients of degree at least twenty-one. These are
bounded by Proposition~\ref{prop:determinant-complete-tails} in
Appendix~\ref{subsec:appB-determinant}. Consequently,
\[
 \begin{aligned}
 \norm{\mathcal E}_{\mathcal A}
 &<305.9385+10^{-9}<310,\\
 \norm{\mathcal Z}_{\mathcal A}
 &<1675.432427+\frac1{100}<1690.
 \end{aligned}
\]
This gives the three norm bounds used below, with each finite sum
coming from Appendix~\ref{subsec:appA-positive} and each infinite
remainder coming from Appendix~\ref{subsec:appB-determinant}.
Although $\mathcal C_0,\mathcal E_0,\mathcal Z_0$ each have degree at
most twenty, their products can have degree forty. Thus their terms
of degrees twenty-one through forty must be included in the
determinant tail. The bound below is the finite-block estimate proved in
Appendix~\ref{subsec:appA-finite-tail}; its row-by-row table records
the contributions from every ordered pair of degrees whose sum is at
least twenty-one:
\[
\begin{aligned}
 \sum_{\substack{0\le i,j\le20\\i+j\ge21}}
       \left|\left[\mathcal C\right]_i\right|
       \left|\left[\mathcal E\right]_j\right|
 +
 \sum_{\substack{0\le i,j\le20\\i+j\ge21}}
       \left|\left[\mathcal Z\right]_i\right|
       \left|\left[\mathcal Z\right]_j\right|<100.
\end{aligned}
\]
Put $R_C=\mathcal C-\mathcal C_0$, and similarly define $R_E,R_Z$.
The identities
\[
 \mathcal C\mathcal E-\mathcal C_0\mathcal E_0
       =R_C\mathcal E+\mathcal C_0R_E,\qquad
 \mathcal Z^2-\mathcal Z_0^2=R_Z\left(\mathcal Z+\mathcal Z_0\right)
\]
and submultiplicativity imply
\[
\begin{aligned}
 \norm{N_{\det}
       -\left(\mathcal C_0\mathcal E_0-\mathcal Z_0^2\right)}_{\mathcal A}
 &\le\norm{R_C}_{\mathcal A}\norm{\mathcal E}_{\mathcal A}
       +\norm{\mathcal C_0}_{\mathcal A}\norm{R_E}_{\mathcal A}
       +2\norm{\mathcal Z}_{\mathcal A}\norm{R_Z}_{\mathcal A}\\
 &<\frac74\cdot310+22000\cdot10^{-9}
                      +\frac{2\cdot1690}{100}
       =576.300022<577.
\end{aligned}
\]
The difference estimated here starts at degree twenty-one, because
each of $R_C,R_E,R_Z$ does. The tail of the polynomial
$\mathcal C_0\mathcal E_0-\mathcal Z_0^2$ contributes less than $100$
by the preceding finite convolution estimate. The triangle inequality
therefore gives the complete determinant tail bound
\[
 \mathcal T_{21}\left(N_{\det}^\#\right)<100+577=677.
\]

The finite coefficient bounds show that the coefficients through
degree eight of $-L,N_e$, and through degree twenty of $N_{\det}$,
are positive, with constant terms greater than $12,63,3500$,
respectively. For every $0\le y\le1$, a remaining coefficient can
decrease the value of a series by no more than its absolute value.
More explicitly, if $A_j>0$ for $0\le j<k$, then for $0\le y\le1$
\[
 A\left(y\right)
 =A_0+\sum_{j=1}^{k-1}A_jy^j+\sum_{j\ge k}A_jy^j
 \ge A_0-\sum_{j\ge k}\left|A_j\right|.
\]
Apply this with $k=9$ to $-L,N_e$ and with $k=21$ to $N_{\det}$.
The first two bounds are Proposition~\ref{prop:quadratic-complete-tails};
the determinant bound combines Proposition~\ref{prop:determinant-complete-tails}
with the finite-block contribution in Appendix~\ref{subsec:appA-finite-tail}.
The three full-series tail bounds give
\[
 -L\left(y\right)>12-\frac{13}{4}>8,\qquad
 N_e\left(y\right)>63-61=2,\qquad
 N_{\det}\left(y\right)>3500-677>2500.
\]
In particular these bounds hold at $y=1$. The reconstruction formulas
and the exact Schur identity~\eqref{eq:schur-numerators} establish
\eqref{eq:positivity} throughout $0<r<R$.

The same positivity also controls points where the gradient vanishes.
The next lemma identifies the limiting defect there, which will be
needed when the local identity is integrated.
\begin{lemma}\label{lem:center}
For $n=4,5$ put $c_\infty=\phi\left(0\right)/U\left(0\right)$ and
$\chi=\alpha c_\infty^2-1$. Write
$U_c=u\left(1\right)=8U\left(0\right)/T$ and
$V_c=v\left(1\right)=2\phi\left(0\right)/R$ for the normalized
centre amplitudes. Then $\chi>0$.
\end{lemma}
\begin{proof}
The explicit reference profiles give
\[
 a_0=-\frac{U''\left(0\right)}{U\left(0\right)}
     =\frac{1+c}{n\left(1-c\right)}>0.
\]
To compute the limits, use the even centre expansions
\[
\begin{aligned}
 U\left(r\right)&=U\left(0\right)
                    \left(1-\frac{a_0r^2}{2}+O\left(r^4\right)\right),\\
 p\left(r\right)&=a_0U\left(0\right)r+O\left(r^3\right),\\
 \phi\left(r\right)&=\phi\left(0\right)
                    \left(1-\frac{\alpha r^2}{2n}+O\left(r^4\right)\right).
\end{aligned}
\]
The coefficient in the last line follows from
$\Delta\phi\left(0\right)=n\phi''\left(0\right)=-\alpha\phi\left(0\right)$.
Similarly, for $W=\Delta U$, the equation $\Delta W=U$ gives
$W'\left(r\right)=U\left(0\right)r/n+O\left(r^3\right)$.
Dividing these expansions by $p$ gives
\[
\begin{gathered}
 s=\frac1{a_0r}+O\left(r\right),\quad
 A=-\frac1r+O\left(r\right),\quad
 V=\frac1{na_0}+O\left(r^2\right),\\
 f=\frac{c_\infty}{a_0r}+O\left(r\right),\qquad
 w=-\frac{\alpha c_\infty}{na_0}+O\left(r^2\right).
\end{gathered}
\]
For $V_0$, the term $w^2$ is bounded, whereas its other two terms
have leading sum
\[
 -mCfw+\alpha f^2
 =\frac{\alpha c_\infty^2}{a_0^2r^2}
       \left(1-\frac mn\right)+O\left(1\right).
\]
For $\Psi$, the leading terms are
$A'=r^{-2}+O\left(1\right)$,
$s^2=a_0^{-2}r^{-2}+O\left(1\right)$, and
$msCV=-m/\left(na_0^2r^2\right)+O\left(1\right)$;
the remaining term $V$ is bounded. Since $n=m+1$, this gives
\[
 V_0=\frac{\alpha c_\infty^2}{na_0^2r^2}+O\left(1\right),\qquad
 \Psi=\frac{1+1/\left(na_0^2\right)}{r^2}+O\left(1\right).
\]
Consequently $G_\infty=\lim_{r\downarrow0}G\left(r\right)$ exists, and
\begin{equation}\label{eq:chi}
 \chi=\frac{16Q V_c^2}{T^2U_c^2}-1
      =na_0^2\left(1+mG_\infty\right),\qquad
 G_\infty=\frac{\chi/\left(na_0^2\right)-1}{m}.
\end{equation}
The same convergent series give $f=c_\infty s+O\left(s^{-1}\right)$,
so $z=c_\infty+O\left(s^{-2}\right)$ and
$g=f-sz=O\left(s^{-1}\right)$.
Also $G=G_\infty+O\left(r^2\right)$ and
$s'=-1/\left(a_0r^2\right)+O\left(1\right)$ imply
$sG_s=sG'/s'=O\left(r^2\right)$.
Consequently $a\to1$, $d\to2G_\infty$, and the limiting coefficient matrix is
\[
 \mathsf M_\infty=
 \begin{pmatrix}
  1&G_\infty\\
  G_\infty&\left(1+\left(m-1\right)G_\infty\right)/m
 \end{pmatrix}.
\]
The reconstruction remains nonsingular at $y=1$.
Indeed $h\left(1\right)=u\left(1\right)=U_c>0$.
The identity $p=ryH/4$ and the centre expansion of $p$ give
$H\left(1\right)=4a_0U\left(0\right)=Ta_0U_c/2>0$.
At $x=0$, the definition of $D$ reduces to
$D\left(1\right)=-U_cH\left(1\right)\ne0$.
The strict numerator bounds at $y=1$ and the exact Schur identity
therefore make $\mathsf M_\infty$ positive definite. Testing it on the
nonzero vector $\left(1,m\right)$ yields
\[
 \begin{aligned}
 a_0^2\begin{pmatrix}1&m\end{pmatrix}
              \mathsf M_\infty\begin{pmatrix}1\\m\end{pmatrix}
 &=a_0^2\left[1+2mG_\infty
                   +m\left(1+\left(m-1\right)G_\infty\right)\right]\\
 &=na_0^2\left(1+mG_\infty\right)=\chi>0.
\end{aligned}
\]

\end{proof}

\section{Proof of Theorem~\ref{thm:main}}\label{sec:proof-main}
We first establish the weak extensions needed to integrate the local
identity. Section~\ref{sec:main-proof} then proves the integral comparison,
the two spectral inequalities, and the equality statements of the theorem.

We return to a real trial function $u\in C_c^\infty\left(\Omega\right)$.
The pointwise identity is already proved wherever $u$ and $\nabla u$
are nonzero. To integrate it, we must treat the nodal set $\left\{u=0\right\}$
and the critical set $\left\{\nabla u=0\right\}$. The latter may have positive
Lebesgue measure. The Sobolev chain rule and the fact that a weak
gradient vanishes almost everywhere on each level set will be used in
their standard forms \cite[Chapter 9]{Brezis}.

\subsection{Weak extensions and exceptional sets}\label{sec:weak-extension}
The expansions near $r=0$ and $r=R$ determine the extensions of the
comparison map and the two divergence fields. The jump of the odd coefficient
$b$ is separated explicitly, because it produces the nodal measure.

\begin{lemma}\label{lem:lipschitz}
Extend $f,G$ evenly and $b$ oddly to the real $s$ axis. Write
$b\left(s\right)=\eta\operatorname{sign}\left(s\right)+b_{\mathrm{reg}}\left(s\right)$.
Here $\operatorname{sign}\left(s\right)$ equals $1,-1,0$ according as $s>0$,
$s<0$, or $s=0$, and $b_{\mathrm{reg}}\left(0\right)=0$. For
$\left(t,\xi\right)\in\R\times\R^n$ with $\xi\ne0$, the maps
\[
\Phi\left(t,\xi\right)=\left|\xi\right|f\left(\left|t\right|/\left|\xi\right|\right),\quad
\mathcal B\left(t,\xi\right)=\left|\xi\right|b_{\mathrm{reg}}\left(t/\left|\xi\right|\right),\quad
\mathcal V\left(t,\xi\right)=G\left(\left|t\right|/\left|\xi\right|\right)\xi
\]
have globally Lipschitz extensions. At $\xi=0$ their values are
$c_\infty\left|t\right|$, $B_\infty t$, and $0$, respectively. Here
$c_\infty=\phi\left(0\right)/U\left(0\right)>0$ and
$B_\infty=\lim_{s\to+\infty}b\left(s\right)/s$ is finite.
\end{lemma}
\begin{proof}
Write $\delta=R-r$. Since
$U'''\left(R\right)=\eta-m/R$ and
$\phi''\left(R\right)=m/R$, the boundary expansions are
\[
\begin{aligned}
 U\left(r\right)&=\tfrac12\delta^2
     -\tfrac16\left(\eta-m/R\right)\delta^3+O\left(\delta^4\right),\\
 p\left(r\right)&=\delta
     -\tfrac12\left(\eta-m/R\right)\delta^2+O\left(\delta^3\right),\\
 \phi\left(r\right)&=\delta+\tfrac{m}{2R}\delta^2
       +O\left(\delta^3\right).
\end{aligned}
\]
Dividing the first and third series by the second gives
\[
 s=\frac\delta2+\frac{\eta-m/R}{12}\delta^2
       +O\left(\delta^3\right),\qquad
 f=1+\frac\eta2\delta+O\left(\delta^2\right)
   =1+\eta s+O\left(s^2\right).
\]
In the expression for $b$, the apparently singular terms cancel:
\[
 A+sV=\frac1\delta+\frac m{2R}+O\left(\delta\right),\qquad
 fw=-\frac1\delta-\eta-\frac m{2R}+O\left(\delta\right).
\]
The cancelled numerator for $G$, followed by the definition of $b$, gives
\[
 G=-\frac{2\left(Z-2\right)}{3R}s+O\left(s^2\right),
 \qquad b=\eta+O\left(s\right).
\]
At the centre, write $U_0=U\left(0\right)>0$ and
$a_0=-U''\left(0\right)/U_0>0$, as in Lemma~\ref{lem:center}.
The even profile series give
\[
 U=U_0+O\left(r^2\right),\quad
 p=a_0U_0r+O\left(r^3\right),\quad
 \phi=\phi\left(0\right)+O\left(r^2\right),\quad
 s=\frac1{a_0r}+O\left(r\right).
\]
Thus $r=1/\left(a_0s\right)+O\left(s^{-3}\right)$ and
$f/s=\phi/U=c_\infty+O\left(s^{-2}\right)$.
For clarity, we spell out the change of variables in these expansions.
The relation
\[
 s=\frac{1}{a_0r}+O(r)
\]
first gives
\[
 r=\frac{1}{a_0s}+O(s^{-3}),\qquad
 r^2=\frac{1}{a_0^2s^2}+O(s^{-4}).
\]
Since $f/s=\phi/U=c_\infty+O(r^2)$, it follows that
\[
 f=c_\infty s+O(s^{-1}).
\]
The centre expansion already obtained for $G$ is
$G(r)=G_\infty+O(r^2)$. Substituting
$r^2=O(s^{-2})$ gives
\[
 G(s)=G_\infty+O(s^{-2}).
\]
We can also compute the coefficient of the linear term in $b$
directly. The centre expansions in Lemma~\ref{lem:center} give
\[
 A=-\frac1r+O(r),\quad
 C=-\frac1r,\quad
 V=\frac1{na_0}+O(r^2),\quad
 f=\frac{c_\infty}{a_0r}+O(r),\quad
 w=-\frac{\alpha c_\infty}{na_0}+O(r^2),
\]
while $G=G_\infty+O(r^2)$ and
$s=(a_0r)^{-1}+O(r)$. Substitution into
$b=-(A+sV)-mGC-fw$ gives
\[
\begin{aligned}
 b
 &=\frac1r-\frac{1}{na_0^2r}
   +\frac{mG_\infty}{r}
   +\frac{\alpha c_\infty^2}{na_0^2r}
   +O(r)\\
 &=\frac{1+mG_\infty}{r}
   +\frac{\alpha c_\infty^2-1}{na_0^2r}
   +O(r).
\end{aligned}
\]
Since $\chi=\alpha c_\infty^2-1
=na_0^2(1+mG_\infty)$ by \eqref{eq:chi}, this becomes
\[
 b=\frac{2(1+mG_\infty)}{r}+O(r)
  =2a_0(1+mG_\infty)s+O(s^{-1}).
\]
Consequently, with
\[
 B_\infty=2a_0(1+mG_\infty)
 =a_0(1+mG_\infty)+\frac{\alpha c_\infty^2-1}{na_0},
\]
we have
\[
 f=c_\infty s+O(s^{-1}),\qquad
 G=G_\infty+O(s^{-2}),\qquad
 b=B_\infty s+O(s^{-1})
 \quad\text{as }s\to+\infty.
\]
In particular, $B_\infty=\lim_{s\to+\infty}b(s)/s$.

These expansions also hold after
termwise differentiation to the orders used below.
The change of variable is analytic in the boundary distance and in
$1/s$ as $r$ approaches $0$ or $R$. The differentiated estimates
therefore follow from the corresponding convergent series.
In particular, subtracting the odd jump gives
\[
 b_{\mathrm{reg}}\left(s\right)
 =B_\infty s-\eta\operatorname{sign}\left(s\right)
       +O\left(\left|s\right|^{-1}\right)\quad\text{as }\left|s\right|\to+\infty.
\]
For example, at infinity the linear term in $f$ cancels in
$f-sf_s=O\left(s^{-1}\right)$, whereas
$b_{\mathrm{reg}}-sb'_{\mathrm{reg}}
=-\eta\operatorname{sign}\left(s\right)+O\left(\left|s\right|^{-1}\right)$.
Also $G_s=O\left(s^{-3}\right)$ and $sG_s=O\left(s^{-2}\right)$.
Near zero the preceding Taylor series give bounded derivatives from either side. On each remaining compact $s$ interval the coefficients
are smooth, since $s'\ne0$.
Consequently all of
$f_s,f-sf_s,b'_{\mathrm{reg}},b_{\mathrm{reg}}-sb'_{\mathrm{reg}}$
and $G,G_s,sG_s$ are bounded for $s\ne0$.

We now check the derivatives of the maps themselves. With
$\rho=\left|\xi\right|>0$ and $\sigma=\left|t\right|/\rho$, at points with $t\ne0$,
\[
\begin{gathered}
 \partial_t\Phi=\operatorname{sign}\left(t\right)f_s\left(\sigma\right),
 \qquad
 \nabla_\xi\Phi=
 \left(f\left(\sigma\right)-\sigma f_s\left(\sigma\right)\right)\frac\xi\rho,\\
 \partial_t\mathcal V
 =\operatorname{sign}\left(t\right)G_s\left(\sigma\right)\frac\xi\rho,\qquad
 D_\xi\mathcal V=G\left(\sigma\right)I_n
       -\sigma G_s\left(\sigma\right)\frac{\left(\xi_i\xi_j\right)_{i,j=1}^n}{\rho^2}.
\end{gathered}
\]
Here $D_\xi$ denotes the Jacobian matrix with respect to $\xi$,
and the displayed matrix has entries $\xi_i\xi_j$.
For the signed variable $s=t/\rho$, the remaining derivatives are
\[
 \partial_t\mathcal B=b'_{\mathrm{reg}}\left(s\right),\qquad
 \nabla_\xi\mathcal B
 =\left(b_{\mathrm{reg}}\left(s\right)
               -s b'_{\mathrm{reg}}\left(s\right)\right)\frac\xi\rho.
\]
Every displayed derivative is uniformly bounded. The series near $r=0$ and $r=R$ give the asserted limits as
$\rho\downarrow0$ with $t\ne0$ fixed.
At $\left(t,\xi\right)=\left(0,0\right)$, continuity follows from
\[
 \left|\Phi\left(t,\xi\right)\right|
   +\left|\mathcal B\left(t,\xi\right)\right|
   \le C_0\left(\left|t\right|+\left|\xi\right|\right),\qquad
 \left|\mathcal V\left(t,\xi\right)\right|\le C_0\left|\xi\right|,
\]
for a constant $C_0$ independent of $t,\xi$.
Integrating the derivative bounds along a line segment now proves the
Lipschitz estimates. Unless a segment is contained in $t=0$ or $\xi=0$,
each of these sets meets it at at most $1$ point; split the segment at
these points and pass to the limiting values by continuity.
A segment contained in $\xi=0$ is
handled directly by the three prescribed limiting values. On $t=0$,
the maps reduce to $\left|\xi\right|$, $0$, and $0$, respectively, so segments
contained in that hyperplane are handled directly as well.
\end{proof}

We use the following standard property of Sobolev functions. If
$v\in W^{1,1}_{\mathrm{loc}}(\Omega)$
and $c\in\mathbb R$, then
\[
 \nabla v=0\quad\text{for almost every }x\in\{v=c\}.
\]
The assertion is understood componentwise and holds for every fixed
level $c$. Applying it to $v=\partial_i u$ and $c=0$ gives
\[
 \partial_j\partial_i u=0
 \quad\text{for almost every }x\in\{\nabla u=0\},
 \qquad 1\le i,j\le n.
\]
This is the standard Sobolev level-set property; see
Brezis~\cite[Chapter~9]{Brezis} or Ziemer~\cite[Chapter~2]{Ziemer}.
We will use precisely this form, rather than any regularity assumption
on the level set.

Put $P=\left|\nabla u\right|$ and $F=\Phi\left(u,\nabla u\right)$, and define
\[
 J_{1,\mathrm{reg}}=\mathcal B\left(u,\nabla u\right)\nabla u,\qquad
 J_2=\left[\left(\Delta u\right)I_n-D^2u\right]\mathcal V\left(u,\nabla u\right).
\]
In particular, $F$ is Lipschitz. Let $L_\Phi$ be a Lipschitz
constant of $\Phi$ with respect to the Euclidean norm on
$\mathbb R\times\mathbb R^n$. Set
\[
 \mathcal U(x)=\bigl(u(x),\partial_1u(x),\ldots,\partial_nu(x)\bigr).
\]
At almost every point, the Sobolev chain rule and the operator-norm
bound for the derivative of a Lipschitz map give
\[
 |D F(x)|=|D(\Phi\circ\mathcal U)(x)|
 \le L_\Phi\,|D\mathcal U(x)|.
\]
Writing out the entries of $D\mathcal U$ yields
\[
 |D\mathcal U(x)|^2
 =\sum_{i=1}^n|\partial_i u(x)|^2
   +\sum_{i,j=1}^n|\partial_i\partial_j u(x)|^2
 =|\nabla u(x)|^2+|D^2u(x)|^2.
\]
Consequently,
\[
 |\nabla F|^2
 \le L_\Phi^2\bigl(|\nabla u|^2+|D^2u|^2\bigr)
 \quad\text{almost everywhere in }\Omega.
\]
This estimate is used only to obtain the required $H^1$ and
integrability bounds; the sharper pointwise Hessian identity is used
for the spectral comparison.
The fields $J_{1,\mathrm{reg}}$ and $J_2$ are compactly supported and belong to $W^{1,\infty}$,
the space of bounded functions with bounded first weak derivatives,
by Lemma~\ref{lem:lipschitz} and the chain rule. The remaining part of
$J_1$ is $\eta P\nabla\left|u\right|$.

We use the following distributional form of Kato's inequality:
for $u\in W^{2,1}_{\mathrm{loc}}(\Omega)$,
\[
 \Delta|u|\ge \operatorname{sign}(u)\,\Delta u
 \quad\text{in }\mathcal D'(\Omega),
 \qquad \operatorname{sign}(0)=0.
\]
Equivalently, the distribution
$\Delta|u|-\operatorname{sign}(u)\Delta u$ is nonnegative and hence
is a Radon measure. For the smooth functions considered here this is
the classical Kato inequality; see Kato~\cite{Kato}. The weighted
product rule below determines its contribution to the divergence
correction.
\begin{lemma}\label{lem:nodal}
There is a finite nonnegative Radon measure $\mu$ such that
\[
 \Delta\left|u\right|=\operatorname{sign}\left(u\right)\Delta u+\mu,
\]
and the weak product rule gives
\begin{equation}\label{eq:nodal}
 \operatorname{div}\left(P\nabla\left|u\right|\right)
=\operatorname{sign}\left(u\right)\left(\nabla P\cdot\nabla u+P\Delta u\right)\dd x+P\mu.
\end{equation}
At a regular nodal hypersurface the last measure is
$2P^2\,\dd\mathcal H^{n-1}$, where $\mathcal H^{n-1}$ denotes
surface measure on the smooth hypersurface. A regular nodal point is
a point with $u=0$ and $P>0$.
\end{lemma}
\begin{proof}
Kato's distributional inequality \cite{Kato} gives
$\mu=\Delta\left|u\right|-\operatorname{sign}\left(u\right)\Delta u\ge0$.
A positive distribution is a Radon measure, that is, a regular Borel
measure finite on compact sets. On either open set $\left\{u>0\right\}$
or $\left\{u<0\right\}$, the ordinary chain rule gives equality without
$\mu$. Hence $\mu$ is supported on the nodal set. It also vanishes
outside the compact support of $u$. Choose
$\psi\in C_c^\infty\left(\Omega\right)$ equal to $1$ in a neighbourhood
of this support. Testing the distributional identity with $\psi$ gives
\[
 \mu\left(\Omega\right)
 =-\int_\Omega\operatorname{sign}\left(u\right)\Delta u\dd x
 \le\int_\Omega\left|\Delta u\right|\dd x<+\infty,
\]
because $\left\langle\Delta\left|u\right|,\psi\right\rangle
=\int_\Omega\left|u\right|\Delta\psi\dd x=0$.
The additional step needed here is multiplication by $P$.
The function $P$ is Lipschitz: it is the norm of the smooth vector
field $\nabla u$, and $\left|\nabla P\right|\le\left|D^2u\right|$
almost everywhere. Extending it by zero and mollifying, choose smooth
compactly supported $P_j$ converging uniformly to $P$
with $\nabla P_j\to\nabla P$ in $L^1$ on a fixed compact set and with
uniformly bounded gradients. For each $j$, the product rule in distributions is
\[
 \operatorname{div}\left(P_j\nabla\left|u\right|\right)
 =\nabla P_j\cdot\nabla\left|u\right|\dd x
  +P_j\operatorname{sign}\left(u\right)\Delta u\dd x+P_j\mu.
\]
This identity passes to the limit: the measure term converges by uniform convergence,
and the gradient term by $L^1$ convergence and boundedness of
$\nabla\left|u\right|$. Using
$\nabla\left|u\right|=\operatorname{sign}\left(u\right)\nabla u$
almost everywhere proves \eqref{eq:nodal}.
At a regular nodal point, the implicit function theorem makes
$\left\{u=0\right\}$ a smooth hypersurface. Orient its unit normal
from $\left\{u<0\right\}$ to $\left\{u>0\right\}$, so that it is
$\nabla u/P$ on the interface. The two traces of the normal derivative
of $\left|u\right|$ are $-P$ and $P$. Integration by parts on the two
sides therefore gives $\mu=2P\,\dd\mathcal H^{n-1}$ there, and
multiplication by $P$ gives the asserted $2P^2$ density.
In particular $P\mu$ is singular with
respect to Lebesgue measure: where $P>0$, it is carried by the smooth
nodal hypersurfaces, and where $P=0$ its weight is zero. A countable
subcover of the regular nodal set suffices for this conclusion.
\end{proof}

\subsection{Integral comparison, scaling, and equality}\label{sec:main-proof}
The preceding weak formulations allow the local identity to be integrated
across the nodal and critical sets. Scaling then yields the spectral
comparison, and Faber--Krahn gives the volume inequality and its rigidity.

\begin{proof}[Proof of Theorem~\ref{thm:main}]
Let $u\not\equiv0$, $u\in C_c^\infty\left(\Omega\right)$ be real, and use the
comparison function $F$ and divergence fields defined above.

\smallskip\noindent\emph{Step 1: the integral comparison.}
Off $\left\{u=0\right\}\cup\left\{P=0\right\}$, apply \eqref{eq:hessian} to $u$ or $-u$.
For $u<0$, the map attached to $-u$ is still
$F=P f\left(\left|u\right|/P\right)$. Its first correction field is
$-P b\left(-u/P\right)\nabla u=P b\left(u/P\right)\nabla u$
by the odd extension of $b$. Its second correction field is $J_2$,
because both the Hessian and the gradient change sign. Thus on both
sign regions the first field splits as
$J_1=J_{1,\mathrm{reg}}+\eta P\nabla\left|u\right|$,
and the same nonnegative local identity applies.
The set $\left\{u=0,P>0\right\}$ is locally a smooth hypersurface and is null for
Lebesgue measure. On the possibly positive-measure set $Z_0=\left\{P=0\right\}$,
the Sobolev level-set property applied to every component
$\partial_i u$ gives $\partial_j\partial_i u=0$ almost everywhere
on $Z_0$, for every $i,j$. Hence $D^2u=0$ and $\Delta u=0$
almost everywhere on $Z_0$.
The function $F-c_\infty\left|u\right|$ is Lipschitz and vanishes
on $Z_0$ by the definition of $\Phi$. Its weak gradient vanishes
there almost everywhere. Since
$\nabla\left|u\right|=\operatorname{sign}\left(u\right)\nabla u=0$
there, it follows that $\nabla F=0$ almost everywhere on $Z_0$.
The regular fields vanish on $Z_0$, so their weak
derivatives vanish almost everywhere there as well. The ordinary
density in \eqref{eq:nodal} also vanishes there, and multiplication by
$P$ removes any measure supported on $Z_0$. Therefore the energy
density on this set reduces to
\[
 -u^2+\alpha\left(c_\infty\left|u\right|\right)^2
 =\chi u^2\ge0
 \quad\text{almost everywhere on }Z_0,
\]
by Lemma~\ref{lem:center} and \eqref{eq:chi}.

Consequently, in distributions,
\[
\left(\left|D^2u\right|^2-u^2-\left|\nabla F\right|^2+\alpha F^2\right)\dd x
+\operatorname{div}\left(J_{1,\mathrm{reg}}+J_2+\eta P\nabla\left|u\right|\right)
=Q_u\dd x+\eta P\mu,
\]
where, explicitly,
\[
\begin{aligned}
 Q_u={}&\left|D^2u\right|^2-u^2-\left|\nabla F\right|^2+\alpha F^2
       +\operatorname{div}\left(J_{1,\mathrm{reg}}+J_2\right)\\
 &+\eta\operatorname{sign}\left(u\right)
       \left(\nabla P\cdot\nabla u+P\Delta u\right).
\end{aligned}
\]
The divergences in this formula are the almost-everywhere weak
derivatives of Lipschitz fields. Thus $Q_u$ is bounded on the compact support;
the pointwise identity proves its sign off the exceptional sets, and
the preceding argument proves its sign on the critical set. Thus
$Q_u\ge0$ almost everywhere. This definition also establishes
integrability without requiring separate estimates for each completed
square near $P=0$. Since $\eta>0$, Lemma~\ref{lem:nodal} shows that
the measure contribution $\eta P\mu$ is nonnegative.

Pair with a cutoff equal to $1$ near the compact support. All divergence
terms integrate to zero: for each compactly supported field $J$,
$\left\langle\operatorname{div}J,\psi\right\rangle
=-\int J\cdot\nabla\psi\dd x=0$.
For each pair $i,j$, integration by parts first in $x_i$ and then in
$x_j$ gives
\[
 \int_\Omega u_{ij}^2\dd x
 =-\int_\Omega u_j u_{iij}\dd x
 =\int_\Omega u_{jj}u_{ii}\dd x.
\]
Summing these equalities yields
\[
 \int_\Omega\left|D^2u\right|^2\dd x
 =\sum_{i,j=1}^n\int_\Omega u_{ii}u_{jj}\dd x
 =\int_\Omega\left(\Delta u\right)^2\dd x.
\]
In particular, the integrated distributional identity is
\[
 \int_\Omega\left(\left(\Delta u\right)^2-u^2
                  -\left|\nabla F\right|^2+\alpha F^2\right)\dd x
 =\int_\Omega Q_u\dd x+\eta\int_\Omega P\dd\mu\ge0,
\]
which is \eqref{eq:comparison-intro}.
The Lipschitz composition $F$ is compactly supported $H^1$ and therefore
belongs to $H^1_0\left(\Omega\right)$ by interior mollification: its
support has positive distance from $\R^n\setminus\Omega$, so mollifiers
of sufficiently small radius stay inside $\Omega$ and converge in $H^1$.
It is nonzero. At any point where $u\ne0$, either $P>0$, in which case
$F=P f\left(\left|u\right|/P\right)>0$, or $P=0$, in which case
$F=c_\infty\left|u\right|>0$. Continuity then gives
$\int_\Omega F^2\dd x>0$.

\smallskip\noindent\emph{Step 2: the spectral inequalities.}
The integral comparison has a fixed reference frequency. A dilation
puts the trial function at that frequency. Set
$E=\int_\Omega\left(\Delta u\right)^2\dd x$, $M=\int_\Omega u^2\dd x$, and
$t=\left(E/M\right)^{1/4}$. Here $M>0$ and $E>0$: if $E=0$, then
$\int_\Omega\left|\nabla u\right|^2\dd x
=-\int_\Omega u\Delta u\dd x=0$, which contradicts $u\ne0$
and compact support. On $t\Omega$ take
$u_t\left(x\right)=u\left(x/t\right)$.
With $y=x/t$, the derivative and volume factors are
\[
 \nabla u_t\left(x\right)=t^{-1}\nabla u\left(y\right),\qquad
 \Delta u_t\left(x\right)=t^{-2}\Delta u\left(y\right),\qquad
 \dd x=t^n\dd y.
\]
Consequently
\[
 E\left(u_t\right)=t^{n-4}E=t^nM=M\left(u_t\right).
\]
Apply \eqref{eq:comparison-intro} to $u_t$ and write
$F_t=\Phi\left(u_t,\nabla u_t\right)$ for its comparison function.
The right side is zero, while $F_t\not\equiv0,\quad F_t\in H^1_0\left(t\Omega\right)$
by Step 1. The variational definition therefore gives
\[
 \lambda_1\left(t\Omega\right)
 \le\frac{\int_{t\Omega}\left|\nabla F_t\right|^2\dd x}
          {\int_{t\Omega}F_t^2\dd x}
 \le\alpha.
\]
For a Dirichlet test function the same substitution multiplies its
energy by $t^{n-2}$ and its squared $L^2$ norm by $t^n$.
The dilation is a bijection between the two Dirichlet test spaces, so
\[
 \lambda_1\left(t\Omega\right)=t^{-2}\lambda_1\left(\Omega\right),
 \qquad
 \lambda_1\left(\Omega\right)\le\alpha t^2
      =\alpha\sqrt{E/M}.
\]
Squaring the last inequality yields
$E/M\ge\lambda_1\left(\Omega\right)^2/\alpha^2$.
For a nonzero $u\in H^2_0\left(\Omega\right)$, choose
$u_j\in C_c^\infty\left(\Omega\right)$ converging to $u$ in $H^2$.
Then $\int u_j^2\dd x\to\int u^2\dd x>0$ and
$\int\left(\Delta u_j\right)^2\dd x
\to\int\left(\Delta u\right)^2\dd x$, so the quotient inequality
passes to the limit. Taking the infimum and using
$\alpha^{-2}=R^4/q^4$ proves the first inequality in \eqref{eq:main}.
The Faber--Krahn theorem \cite[Chapter 3]{Henrot} gives
$\lambda_1\left(\Omega\right)\ge q^2\left(\omega_n/\left|\Omega\right|\right)^{2/n}$ and proves the
second inequality, since
\[
 \frac{R^4}{q^4}
 \left[q^2\left(\frac{\omega_n}{\left|\Omega\right|}\right)^{2/n}\right]^2
 =R^4\left(\frac{\omega_n}{\left|\Omega\right|}\right)^{4/n}.
\]
Its finite-volume open-set formulation follows
from the usual rearrangement argument on $H^1_0$. Finally the ball spectral
formulas show that the first constant is attained on every ball.

\smallskip\noindent\emph{Step 3: equality in the volume comparison.}
Suppose that
$\Gamma_1\left(\Omega\right)=\Gamma_1\left(B_\Omega\right)$.
Both intervening inequalities in \eqref{eq:main} must then be equalities.
In particular
$\lambda_1\left(\Omega\right)^2
=q^4\left(\omega_n/\left|\Omega\right|\right)^{4/n}$,
so equality holds in Faber--Krahn. Its quantitative form
\cite{BrascoDePhilippisVelichkov} implies
\[
 \inf_{y\in\R^n}\left|\Omega\triangle\left(y+B_\Omega\right)\right|=0.
\]
This gives an actual ball, rather than only an approximating sequence
of balls. Indeed a sequence of centres for which the symmetric
differences tend to zero is bounded: otherwise a subsequence of
centres would escape to infinity, and its intersections with the
finite-measure set $\Omega$ would have measure tending to zero.
This contradicts convergence of the symmetric differences to zero.
Pass to a convergent subsequence
of centres. Continuity of ball translations in measure then yields
a ball $B$ with $\left|B\right|=\left|\Omega\right|$ and
$\left|\Omega\triangle B\right|=0$.
Openness implies $\Omega\subset B$: a point of $\Omega$ outside
$\overline B$ has a neighbourhood of positive measure outside $B$,
and a point of $\Omega\cap\partial B$ also has such a neighbourhood.
Both possibilities contradict $\left|\Omega\setminus B\right|=0$.

If $\Omega$ has locally Lipschitz boundary, every point of
$B\setminus\Omega$ belongs to $\partial\Omega$, because
$\left|B\setminus\Omega\right|=0$. An exterior cone at such a point
would have positive measure inside $B$, a contradiction. Hence
$\Omega=B$. Conversely, the ball spectral formulas give equality
in the volume comparison for every ball.
\end{proof}

\subsection{Equality on arbitrary open sets}
For arbitrary open sets, the almost-everywhere identification in
Theorem~\ref{thm:main} admits a more precise characterization in terms
of the clamped Sobolev space.

\begin{proposition}\label{prop:equality}
Let $n=4,5$ and let $\Omega\subset\R^n$ be open with
$0<\left|\Omega\right|<+\infty$. Equality
$\Gamma_1\left(\Omega\right)=\Gamma_1\left(B_\Omega\right)$ holds if and
only if there is a ball $B$ such that
\[
 \Omega\subset B,\qquad \left|\Omega\right|=\left|B\right|,\qquad
 H^2_0\left(\Omega\right)=H^2_0\left(B\right),
\]
where the spaces are identified by zero extension.
\end{proposition}
\begin{proof}
The equality argument in the proof of Theorem~\ref{thm:main} gives
a ball $B$ such that $\Omega\subset B$ and
$\left|\Omega\right|=\left|B\right|$.

The inclusion $H^2_0\left(\Omega\right)\subset H^2_0\left(B\right)$
is immediate from the defining smooth test classes. To see that the
infimum on $\Omega$ is attained, take a minimizing sequence with
$L^2$ norm $1$ and view its zero extensions in $H^2_0\left(B\right)$.
The standard norm equivalence on $H^2_0\left(B\right)$ bounds their
$H^2$ norms by their Laplacian energies. The compact embedding into
$L^2\left(B\right)$ therefore gives a subsequence converging weakly
in $H^2$ and strongly in $L^2$; see \cite[Chapter 9]{Brezis}.
The zero-extension image of $H^2_0\left(\Omega\right)$ is a closed
linear subspace of $H^2_0\left(B\right)$, hence is weakly closed.
Its limit remains in that subspace and has $L^2$ norm $1$.
Weak lower semicontinuity of the Laplacian energy shows that this
limit attains the infimum on $\Omega$.

Since the two first eigenvalues are equal, the same zero extension
also minimizes the ball quotient. The first ball
eigenvalue is simple, and its profile $U_B$ is smooth and strictly
positive in $B$ \cite{GGS}. Consequently $U_B\in H^2_0\left(\Omega\right)$.

Choose $u_j\in C_c^\infty\left(\Omega\right)$ converging to $U_B$ in
$H^2\left(B\right)$. For any $v\in C_c^\infty\left(B\right)$, the multiplier
$h=v/U_B$ is smooth with compact support in $B$: on a compact
neighbourhood of the support of $v$, the positive smooth function
$U_B$ is bounded away from zero. Multiplication by $h$ is continuous
on $H^2\left(B\right)$. Explicitly, for $w\in H^2\left(B\right)$,
\[
\begin{aligned}
 \partial_i\left(hw\right)&=h\partial_iw+\left(\partial_i h\right)w,\\
 \partial_{ij}\left(hw\right)&=h\partial_{ij}w
 +\left(\partial_i h\right)\partial_jw
 +\left(\partial_j h\right)\partial_iw
 +\left(\partial_{ij}h\right)w.
\end{aligned}
\]
The bounded derivatives of $h$ through order two consequently give
$\left\|hw\right\|_{H^2\left(B\right)}
\le C_h\left\|w\right\|_{H^2\left(B\right)}$ for a finite constant
$C_h$ depending only on $h$. Apply this bound to $w=u_j-U_B$. It yields
\[
 \frac v{U_B}u_j\in C_c^\infty\left(\Omega\right),\qquad
 \frac v{U_B}u_j\longrightarrow v\quad\hbox{in }H^2\left(B\right).
\]
Thus every $v\in C_c^\infty\left(B\right)$ belongs to
$H^2_0\left(\Omega\right)$. Taking the $H^2$ closure gives the reverse
inclusion $H^2_0\left(B\right)\subset H^2_0\left(\Omega\right)$.
Conversely,
equality of these spaces identifies the two variational infima, so
$\left|\Omega\right|=\left|B\right|$ gives equality in the volume comparison.
\end{proof}

The condition on the Sobolev spaces in Proposition~\ref{prop:equality}
means that the removed set is invisible to the clamped variational
space. For a compact hole $K\subset B$, this is equivalent to zero clamped
$H^2$-capacity \cite[Proposition 2.1]{FelliRomani}. Zero capacity means
that there are smooth cutoff functions compactly supported in $B$,
equal to $1$ near $K$, whose $H^2$ norms tend to zero.
It is stronger than the hole having zero Lebesgue measure. For example,
a point has zero such capacity when $n\ge4$
\cite[Proposition 2.3]{FelliRomani}, and therefore
$B\setminus\left\{x_0\right\}$ has the same first eigenvalue as $B$ for
$x_0\in B$ and $n=4,5$. On the other hand, remove a closed spherical
hypersurface lying strictly inside $B$. It has zero volume, but a
function in the resulting $H^2_0$ space has zero value trace there
(the Sobolev restriction to the smooth hypersurface), whereas $U_B>0$;
see the trace theorem in \cite[Chapter 9]{Brezis}. Thus the spaces differ and the volume inequality is
strict.

Equality in the first comparison of \eqref{eq:main} alone does not
have the same rigidity. If $\Omega$ is the disjoint union of two
separated balls of radius $\rho$, then
\[
 \lambda_1\left(\Omega\right)=q^2/\rho^2,\qquad
 \Gamma_1\left(\Omega\right)=R^4/\rho^4.
\]
Indeed the quotient of any test function on the union is a weighted
average of the quotients of its nonzero restrictions to the two
components, with weights given by their squared $L^2$ norms.
The infimum is therefore the smaller component eigenvalue; a test
function supported on either ball gives the reverse inequality.
This applies to both variational problems. The first comparison is
an equality, whereas
\[
 \Gamma_1\left(B_\Omega\right)
 =\frac{R^4}{\left(2^{1/n}\rho\right)^4}
 =2^{-4/n}\Gamma_1\left(\Omega\right)
 <\Gamma_1\left(\Omega\right),
\]
because $\left|\Omega\right|=2\omega_n\rho^n$.

\section{Dimension dependence and the six-dimensional obstruction}\label{sec:six}
The proof of Theorem~\ref{thm:main} depends on specific signs of the
reference coefficients. The centre coefficient has the opposite sign in
dimension six. This also gives a direct obstruction to the integral
comparison \eqref{eq:comparison-intro} for the same local map.

Using the same definitions of $q,R,\alpha$ and the normalized
profiles, rational enclosures give the following bounds for
$\chi=\alpha c_\infty^2-1$:
\begingroup
\setlength{\tabcolsep}{\arraycolsep}
\setlength{\LTpre}{\medskipamount}
\setlength{\LTpost}{\medskipamount}
\begin{longtable}{>{$}c<{$}|>{$}r<{$}>{$}r<{$}}
n&\text{lower bound for }\chi&\text{upper bound for }\chi\\
\hline
\endfirsthead
n&\text{lower bound for }\chi&\text{upper bound for }\chi\\
\hline
\endhead
\hline
\endfoot
\hline
\endlastfoot
4&0.10964&0.10965\\
 5&0.03538&0.03539\\
 6&-0.021442&-0.021441\\
\end{longtable}
\endgroup
The four- and five-dimensional intervals are obtained by substituting
the parameter boxes in Appendix~\ref{subsec:appA-parameters} into the
centre-amplitude identity~\eqref{eq:chi} and the formulas for
$a_0,U_c,V_c$ in Lemma~\ref{lem:center}. Thus the appendix supplies
the rational parameter inputs, while the centre-amplitude formulas
come from the main-text lemma. The resulting intervals agree with
Lemma~\ref{lem:center}.
In dimension six the limiting matrix determinant lies between $-.032220$ and $-.032219$.
Define
$a_0=-U''\left(0\right)/U\left(0\right)$ and
$\det_\infty=\lim_{s\to+\infty}\left(ae-d^2/4\right)$. The centre expansions give
\[
a_0=\frac{1+c}{n\left(1-c\right)},\qquad
G_\infty=\frac{\chi/\left(na_0^2\right)-1}{n-1}.
\]
The centre expansions used in \eqref{eq:chi} follow from the radial
equations in every dimension considered here; their derivation does
not use coefficient positivity. In particular,
$g=O\left(s^{-1}\right)$ and $G_s=O\left(s^{-3}\right)$ as
$s\to+\infty$. Thus $a\to1$, $d\to2G_\infty$, and
$e\to\left[1+\left(n-2\right)G_\infty\right]/\left(n-1\right)$.
The determinant is therefore
\[
 \det_\infty
 =\frac{1+\left(n-2\right)G_\infty-\left(n-1\right)G_\infty^2}{n-1}
 =\frac{\left(1-G_\infty\right)
             \left(1+\left(n-1\right)G_\infty\right)}{n-1}.
\]
Since $a_0>0$, the two negative signs are
linked by
\[
 1+\left(n-1\right)G_\infty=\frac{\chi}{na_0^2}<0,
 \qquad G_\infty<-\frac1{n-1}.
\]
The last inequality gives $1-G_\infty>0$, so $\chi<0$ forces
$\det_\infty<0$. Both signs express the same centre defect.
Both sign enclosures follow from the
normalized Bessel series through index $40$, using the rational brackets
\[
\begin{gathered}
5.90567823542052287967<R<5.90567823542052287969,\\
5.13562230184068255630<q<5.13562230184068255631.
\end{gathered}
\]
Here $r^2/4<9$, $\nu=2$, and the absolute tail after $j=40$ is bounded by
$\frac{100}{99}\frac{2\,9^{41}}{41!43!}<10^{-55}$; the next-term ratio
is at most $9/\left(42\cdot44\right)<1/100$.
The numerical inequality here is an exact comparison of integers
after multiplication by $99\cdot41!43!\cdot10^{55}$.
More explicitly, let $B_\mu,M_\mu$ be the normalized series defined
in Appendix~\ref{subsec:appA-parameters}; the term recurrence and
outward interval evaluation used there are applied here with
$r^2/4<9$ and truncation index $40$. The rational brackets for $R,q$
displayed above are the separate six-dimensional inputs to this
calculation; they are not entries in the four- and five-dimensional
parameter table of Appendix~\ref{subsec:appA-parameters}. In dimension
six the central values are
\[
 U\left(0\right)=\frac12\left[
       \frac1{M_2\left(R\right)}-\frac1{B_2\left(R\right)}\right],
 \qquad
 \phi\left(0\right)=\frac{6R}{q^2B_3\left(q\right)}.
\]
Consequently the quantity to enclose is
\[
 \chi+1=\frac{144}{q^2B_3\left(q\right)^2
       \left[1/M_2\left(R\right)-1/B_2\left(R\right)\right]^2}.
\]
For each normalized series, substitute the root intervals into its
polynomial through degree $40$ in $r^2/4$ and enlarge the resulting
interval by $10^{-55}$ on both sides. The denominators above are
bounded away from zero: $M_2\left(R\right)>0$,
$B_2\left(R\right)<0$, and $B_3\left(q\right)>0$.
Interval division therefore gives the stated enclosure for $\chi$.
Substitution of $c=B_2\left(R\right)/M_2\left(R\right)$ into $a_0$
and then into $G_\infty$ gives the determinant enclosure by the same
rational operations.
The first-root isolation argument in
Appendix~\ref{subsec:appA-parameters} applies here: the lower bounds
for the first zeros and the Sturm-comparison gap estimate place the
displayed brackets before the respective second zeros, while the
endpoint signs and the cross-product derivative identify the first
positive cross-product zero. In the present six-dimensional case, the
specific checks are $j_{2,1}\ge\sqrt{12}$,
$j_{3,1}\ge4$, and a gap greater than $\pi$ between consecutive
zeros.
The displayed brackets lie before the respective second zeros;
the same forty-term sum gives $B_3\left(R\right)>3/100>0$.
Consequently $q<R<j_{3,1}$, and the cross-product derivative
identifies the displayed $R$ as its first positive zero.
In particular $\chi_6<-1/50$ and $\det_\infty<-3/100$.
The negative determinant prevents this particular quadratic completion
from being positive. The negative centre coefficient has a stronger
consequence: the integral comparison itself fails for this map,
regardless of the chosen divergence correction.

The proof of the Lipschitz property of $\Phi$ in
Lemma~\ref{lem:lipschitz} uses only the expansions at $r=0$ and $r=R$ and
$s'<0$. These properties also hold in dimension six, so $\Phi$
remains globally Lipschitz independently of the coefficient signs
in \eqref{eq:positivity}.

Choose a fixed smooth function $\theta:\R\to\left[0,1\right]$ equal to $1$ on
$\left(-\infty,0\right]$ and $0$ on $\left[1,+\infty\right)$. Set
$u_L\left(x\right)=\theta\left(\left|x\right|-L\right)$ for $L>1$. This is smooth, equals $1$ on $B_L$, and belongs to
$C_c^\infty\left(B_{L+2}\right)$, since its support is contained in
$\overline B_{L+1}$. On $B_L$ its comparison
function is $F_{u_L}=c_\infty$; hence the density of the right side
of \eqref{eq:comparison-intro} minus the left side is $\chi$.
To estimate the transition shell $B_{L+1}\setminus B_L$, write
$r=\left|x\right|$ and $e_r=x/r$. Differentiating the radial cutoff gives
\[
\begin{aligned}
 \nabla u_L&=\theta'\left(r-L\right)e_r,\\
 D^2u_L&=\theta''\left(r-L\right)e_r\otimes e_r
       +\frac{\theta'\left(r-L\right)}r
                  \left(I_6-e_r\otimes e_r\right),\\
 \Delta u_L&=\theta''\left(r-L\right)
                  +\frac5r\theta'\left(r-L\right).
\end{aligned}
\]
Here $e_r\otimes e_r$ is the matrix with entries
$\left(e_r\right)_i\left(e_r\right)_j$.
Put $M_{\theta,1}=\sup_{\R}\left|\theta'\right|$ and
$M_{\theta,2}=\sup_{\R}\left|\theta''\right|$. Since $r\ge L>1$,
\[
 \left|\nabla u_L\right|\le M_{\theta,1},\qquad
 \left|D^2u_L\right|^2\le M_{\theta,2}^2+5M_{\theta,1}^2,\qquad
 \left|\Delta u_L\right|\le M_{\theta,2}+5M_{\theta,1}.
\]
Let $L_\Phi$ be a Lipschitz constant for $\Phi$ with respect to the
Euclidean norm on $\R\times\R^6$. Since $\Phi\left(0,0\right)=0$,
the composition rule gives, almost everywhere on the shell,
\[
 \left|F_{u_L}\right|^2\le L_\Phi^2\left(1+M_{\theta,1}^2\right),\qquad
 \left|\nabla F_{u_L}\right|^2
 \le L_\Phi^2\left(\left|\nabla u_L\right|^2
                         +\left|D^2u_L\right|^2\right)
 \le L_\Phi^2\left(M_{\theta,2}^2+6M_{\theta,1}^2\right).
\]
Thus the absolute value of the defect density on the shell is at most
\[
 C_\theta=\left(M_{\theta,2}+5M_{\theta,1}\right)^2+1
       +L_\Phi^2\left(M_{\theta,2}^2+6M_{\theta,1}^2\right)
       +\alpha L_\Phi^2\left(1+M_{\theta,1}^2\right),
\]
independently of $L$. Its integral is bounded in absolute value by
$C_\theta\omega_6\left(\left(L+1\right)^6-L^6\right)$.
Define the comparison defect
\[
 \mathcal D_L=\int_{\R^6}
 \left(\left(\Delta u_L\right)^2-u_L^2
       -\left|\nabla F_{u_L}\right|^2+\alpha F_{u_L}^2\right)\dd x.
\]
The interior contribution is exactly $\chi\omega_6L^6$, and the shell
estimate gives
\[
 \left|\frac{\mathcal D_L}{\omega_6L^6}-\chi\right|
 \le C_\theta\left(\left(1+1/L\right)^6-1\right)
 \longrightarrow0.
\]
For $L\ge1$, the binomial formula gives
\[
 \left(1+1/L\right)^6-1
 =\sum_{j=1}^6\binom6jL^{-j}\le\frac{63}{L}.
\]
For example, $L>126C_\theta/\left|\chi\right|$ makes the last error
smaller than $\left|\chi\right|/2$, and hence
$\mathcal D_L<\chi\omega_6L^6/2<0$.
Thus \eqref{eq:comparison-intro} fails in dimension six for this map,
even on smooth nonnegative test functions in smooth balls.

\begin{figure}[!htbp]
\centering
\begin{tikzpicture}[x=1cm,y=1cm,font=\small,
 every path/.style={line width=.6pt},
 arr/.style={->,>=latex,line width=.65pt}]
\node at (0,2.35) {A section of the support in $\R^6$};
\fill[pattern=north east lines,pattern color=black!35] (0,0) circle (1.8);
\draw (0,0) circle (1.8);
\draw[fill=white] (0,0) circle (1.4);
\node[text badly centered,text width=2.5cm] at (0,.55) {$u_L=1$\\$F_{u_L}=c_\infty$};
\fill (0,0) circle (.025);
\draw (0,0)--(1.4,0);
\node[below] at (.7,0) {$L$};
\draw[arr] (2.25,-1.5)--(1.23,-1.1);
\node[anchor=west,text ragged,text width=3.2cm,font=\footnotesize] at (2.1,-1.6)
 {Shell of width $1$};
\node at (0,-2.3) {Interior: $\chi\omega_6L^6<0$};

\begin{scope}[shift={(5.2,-.65)}]
\node at (2,3) {Radial cutoff $u_L\left(r\right)$};
\fill[black!8] (2.2,0) rectangle (3.35,2);
\draw[arr] (-.15,0)--(4.4,0) node[right] {$r$};
\draw[arr] (0,-.15)--(0,2.55);
\draw[dashed,black!55] (2.2,0)--(2.2,2);
\draw[dashed,black!55] (3.35,0)--(3.35,2);
\draw[line width=1pt] (0,2)--(2.2,2)
.. controls (2.65,2) and (2.85,0).. (3.35,0)--(4.05,0);
\node[left] at (0,2) {$1$};
\node[below left] at (0,0) {$0$};
\node[below] at (2.2,0) {$L$};
\node[below] at (3.35,0) {$L+1$};
\node[text badly centered,text width=5.5cm] at (2,-1.65) {Transition contribution: $O\left(L^5\right)$};
\end{scope}
\node[draw,rounded corners=2pt,fill=black!3,inner sep=8pt]
 at (3.35,-3.3)
 {$\mathcal D_L=\chi\omega_6L^6+O\left(L^5\right)<0
   \quad\text{for sufficiently large }L$};
\end{tikzpicture}
\caption{The six-dimensional obstruction, shown schematically.
The shaded annulus and shaded interval depict the same unit-width
transition region. The cutoff is smooth and constant on $B_L$.
As $L$ grows, the negative interior contribution dominates the shell
contribution. The integral comparison for this map therefore fails in dimension six.}
\label{fig:plateau}
\end{figure}
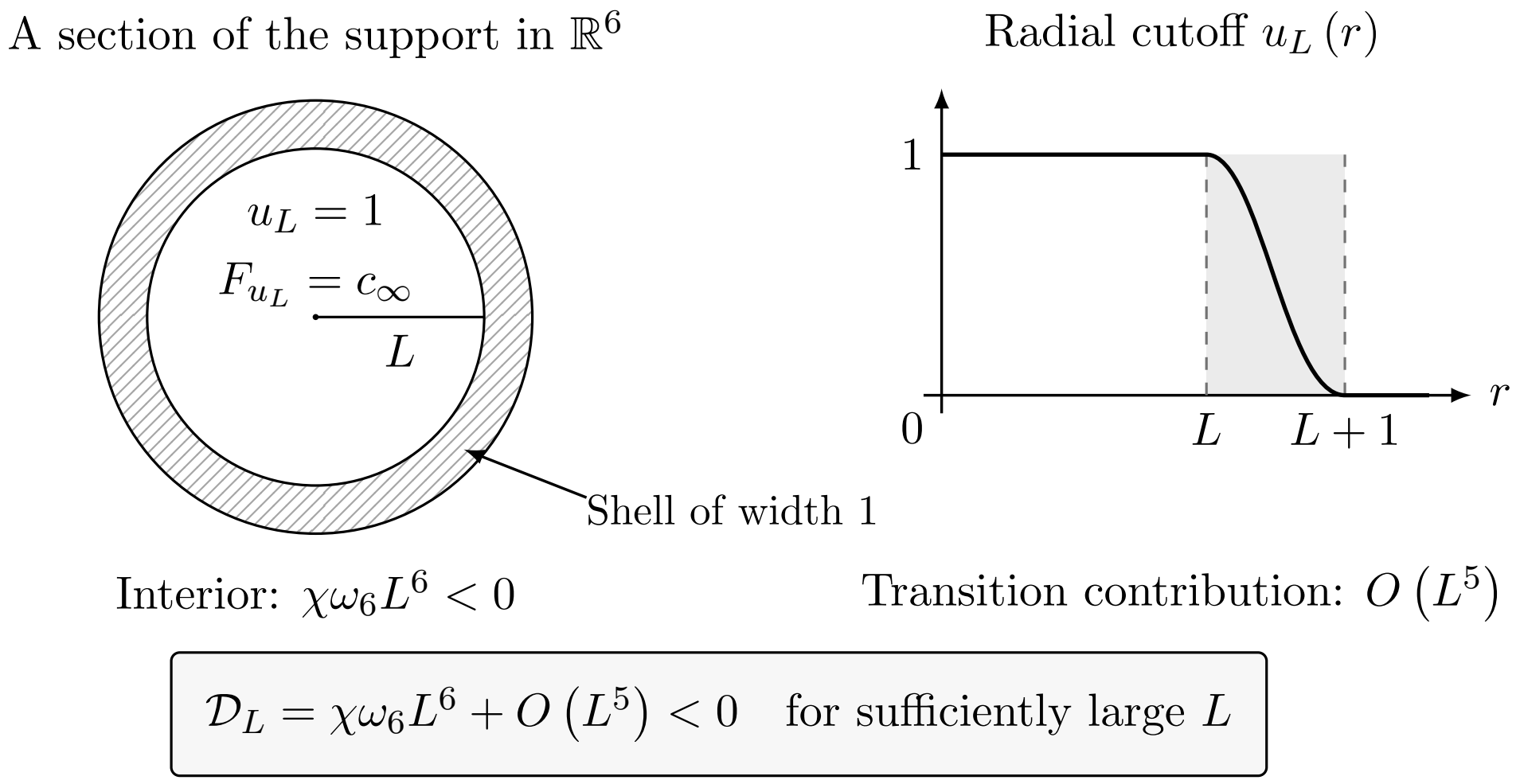

Figure~\ref{fig:plateau} separates the two contributions: the negative
interior term has order $L^6$, while the shell has fixed width and
contributes only order $L^5$.

For the same family of maps in dimension $n$, write
$\chi_n=\alpha_n\left[\phi_n\left(0\right)/U_n\left(0\right)\right]^2-1$,
where the subscripts indicate the reference data. The plateau argument
shows that $\chi_n\ge0$ is necessary for a universal integral comparison.
The proof in dimensions four and five also uses the coefficient
inequalities \eqref{eq:positivity} for the full Hessian identity.
The calculation above excludes this map in dimension six.
The signs for $n\ge7$ and other comparison maps are not determined here.

\appendix
\section{Rational coefficient certificate}\label{sec:paper-input}
This appendix specifies the exact arithmetic used in
Section~\ref{sec:certificate}. It consists of parameter enclosures,
outward-rounded profile recurrences, and finite coefficient sums.
The tables include the auxiliary coefficients and every integer
division remainder, followed by the complete finite convolution sums.
Subsection~\ref{subsec:appA-parameters} certifies the spectral
parameters; Subsection~\ref{subsec:appA-profiles} gives the profile
coefficient intervals; Subsection~\ref{subsec:appA-formulas} derives
the finite numerator coefficients; Subsection~\ref{subsec:appA-positive}
certifies their positivity; and Subsection~\ref{subsec:appA-finite-tail}
accounts for products of finite blocks whose total degree is at least
twenty-one.

\subsection{Parameter enclosures for the spectral data}\label{subsec:appA-parameters}
For $n=4,5$, the following table encloses $R,q,T,Q,Z$.
Every finite decimal in an enclosure or a bound is rational.
\begingroup
\setlength{\LTpre}{\medskipamount}
\setlength{\LTpost}{\medskipamount}
\begin{longtable}{>{$}c<{$}>{$}c<{$}|>{$}c<{$}}
n&\text{parameter}&\text{rational interval}\\
\hline
\endfirsthead
n&\text{parameter}&\text{rational interval}\\
\hline
\endhead
\hline
\endfoot
\hline
\endlastfoot
4&R&\left(4.610899879049055,4.610899879049056\right)\\
4&q&\left(3.831705970207512,3.831705970207513\right)\\
4&T&\left[21.260397694614,21.260397694615\right]\\
4&Q&\left[14.681970642123,14.681970642124\right]\\
4&Z&\left[3.218890982800,3.218890982801\right]\\
5&R&\left(5.267657530336814,5.267657530336815\right)\\
5&q&\left(4.493409457909064,4.493409457909065\right)\\
5&T&\left[27.748215856914,27.748215856915\right]\\
5&Q&\left[20.190728556426,20.190728556427\right]\\
5&Z&\left[3.501551477382,3.501551477383\right]\\
\end{longtable}
\endgroup
Squaring preserves order on positive numbers. Thus the $Q=q^2$
and $T=R^2$ intervals follow by squaring the positive endpoints of
the $q$ and $R$ intervals and rounding outwards.
The calculation of $Z$ below involves a
positive denominator, so its lower bound uses the lower numerator and
upper denominator, and its upper bound uses the opposite endpoints.
To verify these enclosures, use the normalized entire functions
\[
 B_\mu\left(t\right)=\sum_{j\ge0}\frac{\left(-1\right)^j\left(t^2/4\right)^j}{j!\left(\mu+1\right)_j},
 \qquad M_\mu\left(t\right)=\sum_{j\ge0}\frac{\left(t^2/4\right)^j}{j!\left(\mu+1\right)_j}.
\]
Here $\left(a\right)_0=1$ and $\left(a\right)_j=a\left(a+1\right)\cdots\left(a+j-1\right)$ for $j\ge1$ is
the rising factorial. These series equal
$2^\mu\mathrm{\Gamma}\left(\mu+1\right)t^{-\mu}J_\mu\left(t\right)$ and
$2^\mu\mathrm{\Gamma}\left(\mu+1\right)t^{-\mu}I_\mu\left(t\right)$, respectively, where
$\mathrm{\Gamma}$ is the Euler gamma function
\cite[Sections 10.2 and 10.25]{DLMF}.
For $\mu\ge1$ and $t^2/4<7$, the absolute tail after $j=24$ is at most
\[
 \frac{100}{99}\frac{7^{25}}{25!26!}<3\cdot10^{-31}.
\]
Indeed, the first omitted term is at most $7^{25}/\left(25!26!\right)$, and every
subsequent ratio is at most $7/\left(26\cdot27\right)<1/100$.
More explicitly, at a rational argument $t$ the positive summands can
be generated without powers or factorials by starting with $a_0=1$ and
using
\[
 a_{j+1}=\frac{t^2a_j}{4\left(j+1\right)\left(\mu+j+1\right)}.
\]
The sums $\sum_{j=0}^{24}\left(-1\right)^ja_j$ and
$\sum_{j=0}^{24}a_j$, enlarged at each end by $3\cdot10^{-31}$,
therefore enclose $B_\mu\left(t\right)$ and $M_\mu\left(t\right)$.
For a root interval the same recurrence uses interval multiplication;
for an endpoint sign it uses the single rational endpoint. In a cross
product, the enclosures are multiplied before the two products are
added, so the stated error bound is applied to each factor rather than
to the product as a whole.
At the stated left and right root endpoints, $B_\nu$ and the normalized
cross product $B_\nu M_{\nu+1}+B_{\nu+1}M_\nu$ have signs respectively
$>10^{-18}$ and $<-10^{-18}$. The same finite sums give
$B_{\nu+1}\left(R\right)>0$, $-B_\nu\left(R\right)>1/20$, $B_{\nu+1}\left(q\right)>1/10$ and
\[
 Z=\frac Tn\frac{M_{\nu+1}\left(R\right)}{M_\nu\left(R\right)}
\]
in its stated box.

The root interpretation uses the standard simplicity and interlacing of
Bessel zeros \cite[Section 10.21]{DLMF} and the ball spectral formulas
already cited. The isolation can be checked directly as follows. Before
its first zero, the radial equation gives
$B_\nu\left(t\right)\ge1-t^2/\left(2n\right)$, so that zero is at least $\sqrt{2n}$.
The Sturm comparison theorem applied to
$\sqrt{t}J_\nu\left(t\right)$ gives a gap greater than $\pi$ between consecutive
zeros for the orders used here. The enclosures lie below
$\sqrt{2n}+3$ and hence before the second zero. The same argument with
order $\nu+1$, and the verified sign $B_{\nu+1}\left(R\right)>0$, gives
$q<R<j_{\nu+1,1}$. Finally, if
$\mathcal C=J_\nu I_{\nu+1}+J_{\nu+1}I_\nu$, the Bessel recurrence
relations give $\left(t\mathcal C\right)'=2tJ_\nu I_\nu$. Thus
$t\mathcal C$ is strictly decreasing on $\left(q,j_{\nu+1,1}\right)$, where its
endpoint signs differ. This identifies the first positive zero of the cross product.

\subsection{Outward-rounded profile coefficient intervals}\label{subsec:appA-profiles}
We next construct intervals for the Taylor coefficients $u_j,v_j,W_j$.
Their recurrences are given in Section~\ref{sec:certificate}.
Nonnegativity of $u_j,v_j$ follows independently from the Bessel
products at the start of Appendix~\ref{sec:tail}.
Put $\mathsf P=10^{12}$ and $\mathsf D=10^{16}$.
For an interval $I=\left[a,b\right]$, define its outward rounding by
\[
 \mathcal R\left(I\right)=
 \left[\frac{\left\lfloor\mathsf D a\right\rfloor}{\mathsf D},
       \frac{\left\lceil\mathsf D b\right\rceil}{\mathsf D}\right].
\]
The floor and ceiling are the greatest integer below and the least
integer above a real number. Interval addition, subtraction, and
multiplication mean
\[
\begin{gathered}
 \left[a,b\right]+\left[c,d\right]=\left[a+c,b+d\right],\qquad
 -\left[a,b\right]=\left[-b,-a\right],\\
 \left[a,b\right]\left[c,d\right]=
 \left[\min\left(ac,ad,bc,bd\right),\max\left(ac,ad,bc,bd\right)\right].
\end{gathered}
\]
Division by a positive rational multiplies both endpoints by its
reciprocal. Let $\mathbf T,\mathbf Q,\mathbf Z$ denote the parameter
intervals above. For the recurrence involving $T^2$, use
\begingroup
\setlength{\tabcolsep}{\arraycolsep}
\setlength{\LTpre}{\medskipamount}
\setlength{\LTpost}{\medskipamount}
\begin{longtable}{>{$}c<{$}|>{$}r<{$}>{$}r<{$}}
n&\mathsf P\inf\mathbf S&\mathsf P\sup\mathbf S\\
\hline
\endfirsthead
n&\mathsf P\inf\mathbf S&\mathsf P\sup\mathbf S\\
\hline
\endhead
\hline
\endfoot
\hline
\endlastfoot
4&452004510133148&452004510133191\\
 5&769963483241893&769963483241950\\
\end{longtable}
\endgroup
Squaring the endpoints of $\mathbf T$ proves $T^2\in\mathbf S$.
Start from the exact intervals
$\mathbf U_0=\mathbf U_1=\mathbf V_0=\left[0,0\right]$,
$\mathbf V_1=\mathbf W_0=\left[1,1\right]$, and
$\mathbf W_1=-\mathbf Z/2$. For $0\le j\le15$ use
\[
\begin{aligned}
 \mathbf U_{j+2}
 &=\mathcal R\left(\frac{2j+n}{2\left(j+2\right)}\mathbf U_{j+1}
       +\frac{2\mathbf W_j}{\left(j+1\right)\left(j+2\right)}\right)
       \cap\left[0,+\infty\right),\\
 \mathbf V_{j+2}
 &=\mathcal R\left(\frac{2j+n}{2\left(j+2\right)}\mathbf V_{j+1}
       -\frac{\mathbf Q\mathbf V_j}{4\left(j+1\right)\left(j+2\right)}\right)
       \cap\left[0,+\infty\right),\\
 \mathbf W_{j+2}
 &=\mathcal R\left(\frac{2j+n}{2\left(j+2\right)}\mathbf W_{j+1}
       +\frac{\mathbf S\mathbf U_j}{32\left(j+1\right)\left(j+2\right)}\right).
\end{aligned}
\]
Each update uses only already constructed intervals. To see the
endpoint choices in the subtraction term, write
$\mathbf Q=\left[Q_-,Q_+\right]$ and
$\mathbf V_j=\left[v_-,v_+\right]$. Both intervals are nonnegative,
and hence
\[
 -\mathbf Q\mathbf V_j=\left[-Q_+v_+,-Q_-v_-\right].
\]
This explains why the lower recurrence for $v$ uses the upper
endpoints of both factors. The $W$ intervals may contain negative
numbers, but they are multiplied only by positive scalars in the first
term of their recurrence. The product $\mathbf S\mathbf U_j$ in the
second term has nonnegative factors. Every operation therefore
encloses the corresponding exact recurrence value; outward rounding
only enlarges that enclosure. Intersecting with $\left[0,+\infty\right)$
is justified by the independently established nonnegativity of the
profile coefficients. Starting from the exact initial values,
induction proves $u_j\in\mathbf U_j$, $v_j\in\mathbf V_j$, and
$W_j\in\mathbf W_j$ at every subsequent index.
The following integer form of the same recurrence makes every
rounding step explicit. In the tables write
\[
\begin{aligned}
 U_j^-&=\max\left(0,a_j^u\right),& U_j^+&=b_j^u,\\
 V_j^-&=\max\left(0,a_j^v\right),& V_j^+&=b_j^v,\\
 W_j^-&=a_j^W,& W_j^+&=b_j^W.
\end{aligned}
\]
Thus the assertions to be checked are
\[
 \frac{U_j^-}{\mathsf D}\le u_j\le\frac{U_j^+}{\mathsf D},\qquad
 \frac{V_j^-}{\mathsf D}\le v_j\le\frac{V_j^+}{\mathsf D},\qquad
 \frac{W_j^-}{\mathsf D}\le W_j\le\frac{W_j^+}{\mathsf D}.
\]
The lower endpoints in the $u,v$ columns are replaced by their
positive parts after the indicated integer division.
The $W$ interval is not clipped. Thus
$\mathbf U_j=\left[U_j^-/\mathsf D,U_j^+/\mathsf D\right]$,
$\mathbf V_j=\left[V_j^-/\mathsf D,V_j^+/\mathsf D\right]$, and
$\mathbf W_j=\left[W_j^-/\mathsf D,W_j^+/\mathsf D\right]$.
Write $T\in\left[t_-/\mathsf P,t_+/\mathsf P\right]$,
$Q\in\left[q_-/\mathsf P,q_+/\mathsf P\right]$, and
$S=T^2\in\left[s_-/\mathsf P,s_+/\mathsf P\right]$.
The integers $t_\pm$ and $q_\pm$ are obtained by multiplying the
endpoints in the $T$ and $Q$ rows of the parameter table by
$\mathsf P$, respectively.
The following integers give the slightly enlarged $S$ boxes and certify
their outward rounding:
\begingroup
\setlength{\tabcolsep}{\arraycolsep}
\setlength{\LTpre}{\medskipamount}
\setlength{\LTpost}{\medskipamount}
\begin{longtable}{>{$}c<{$}|>{$}r<{$}>{$}r<{$}|>{$}r<{$}>{$}r<{$}}
n&s_-&s_+&t_-^2-\mathsf P s_-&\mathsf P s_+-t_+^2\\
\hline
\endfirsthead
n&s_-&s_+&t_-^2-\mathsf P s_-&\mathsf P s_+-t_+^2\\
\hline
\endhead
\hline
\endfoot
\hline
\endlastfoot
4&452004510133148&452004510133191&286004608996&193200001775\\
5&769963483241893&769963483241950&551321603396&952246682775\\
\end{longtable}
\endgroup
Both remainder columns lie in $\left[0,\mathsf P\right)$.
Here are integer identities that certify every subsequent rounding.
For $0\le j\le15$, set
\[
 \mathfrak a_j=\left(2j+n\right)\left(j+1\right),\qquad
 \mathfrak b_j=2\left(j+1\right)\left(j+2\right).
\]
The remainder columns in the tables mean exactly
\[
\begin{aligned}
r_{u,-}={}&\mathfrak a_jU_{j+1}^-+4W_j^-
                       -\mathfrak b_ja_{j+2}^u,\\
r_{u,+}={}&\mathfrak b_jb_{j+2}^u
                       -\mathfrak a_jU_{j+1}^+-4W_j^+,\\
r_{W,-}={}&16\mathsf P\mathfrak a_jW_{j+1}^-+s_-U_j^-
                       -16\mathsf P\mathfrak b_ja_{j+2}^W,\\
r_{W,+}={}&16\mathsf P\mathfrak b_jb_{j+2}^W
                       -16\mathsf P\mathfrak a_jW_{j+1}^+-s_+U_j^+,\\
r_{v,-}={}&2\mathsf P\mathfrak a_jV_{j+1}^--q_+V_j^+
                       -2\mathsf P\mathfrak b_ja_{j+2}^v,\\
r_{v,+}={}&2\mathsf P\mathfrak b_jb_{j+2}^v
                       -2\mathsf P\mathfrak a_jV_{j+1}^++q_-V_j^-.
\end{aligned}
\]
In each row the printed integers satisfy these identities and
\[
 0\le r_{u,\pm}<\mathfrak b_j,\quad
 0\le r_{W,\pm}<16\mathsf P\mathfrak b_j,\quad
 0\le r_{v,\pm}<2\mathsf P\mathfrak b_j.
\]
These inequalities prove that the raw lower and upper integers are the
outward floor and ceiling of the appropriate recurrence bounds.
Starting with $u_0=u_1=v_0=0$, $v_1=W_0=1$, and $W_1=-Z/2$,
induction therefore proves all the displayed coefficient enclosures.
For example, in dimension four the $u_3$ row has $j=1$,
$\mathfrak a_1=12$, $\mathfrak b_1=12$, and
\[
\begin{aligned}
 12\cdot10^{16}-4\cdot16094454914005000
       -12\cdot4635181695331666&=8,\\
 12\cdot4635181695333334
       -12\cdot10^{16}+4\cdot16094454914000000&=8.
\end{aligned}
\]
For comparison, the membrane recurrence in this dimension gives
$v_2=1$ and $v_3=1-Q/24$. Substitution of the two endpoints of
$\mathbf Q$ gives
\[
 \frac{3882512232448333}{10^{16}}<v_3
 <\frac{3882512232448750}{10^{16}}.
\]
The second displayed endpoint is the exact upper value of the
interval calculation; the strict inequality for the actual coefficient
follows from the strict inclusion of $Q$ in its parameter box.
The $v_3$ row lists the same outward-rounded intervals. At indices
where a raw lower integer is negative, such as the four-dimensional
$u_{14}$ row, the certified lower bound is zero, while we use the positive
upper bound. The remaining rows have the same integer
interpretation.

Only indices up to $25$ will be needed in the coefficient formulas below.
For every $j\ge18$ it suffices to use the common interval
\[
 0\le u_j,v_j\le10^{-16}.
\]
Indeed, the series bound in Appendix~\ref{sec:tail} gives
$u_j,v_j\le20\,7^j/\left[j!\left(j+1\right)!\right]$ for $j\ge12$.
To check the bound at index $18$, compare the following integers:
\[
\begin{aligned}
 20\,7^{18}\,10^{16}&=325682719582089800000000000000000,\\
 18!19!&=778817392288148379660189696000000.
\end{aligned}
\]
The second integer is larger than the first, proving the assertion
at $18$. The successive majorant ratio
$7/\left[\left(j+1\right)\left(j+2\right)\right]$ is less than $1$.
These positive intervals also cover the high-degree coefficients; they do not
set them equal to zero.
{\small\begingroup
\setlength{\tabcolsep}{\arraycolsep}
\setlength{\LTpre}{\medskipamount}
\setlength{\LTpost}{\medskipamount}
\begin{longtable}{>{$}c<{$}|>{$}r<{$}>{$}r<{$}|>{$}r<{$}>{$}r<{$}}
\multicolumn{5}{c}{$\text{Dimension 4.}$}\\[4pt]
k&a_k^u&b_k^u&a_k^v&b_k^v\\
\hline
\endfirsthead
\multicolumn{5}{c}{$\text{Dimension 4.}$}\\[4pt]
k&a_k^u&b_k^u&a_k^v&b_k^v\\
\hline
\endhead
\hline
\endfoot
\hline
\endlastfoot
0&0&0&0&0\\
1&0&0&10000000000000000&10000000000000000\\
2&10000000000000000&10000000000000000&10000000000000000&10000000000000000\\
3&4635181695331666&4635181695333334&3882512232448333&3882512232448750\\
4&1952772542997499&1952772543000001&823768348672499&823768348673125\\
5&343327051596999&343327051600001&111231715991250&111231715992002\\
6&55093442977825&55093442981236&10443860065168&10443860066004\\
7&5099424289763&5099424293538&723022037384&723022038287\\
8&441038127104&441038131197&38484328100&38484329059\\
9&25218593774&25218598143&1625340809&1625341815\\
10&1373541325&1373545935&55824727&55825773\\
11&53368820&53373640&1590134&1591214\\
12&2000437&2005442&37799&38910\\
13&51921&57090&359&1497\\
14&-3062&2253&-426&735\\
15&-1630&754&-27&729\\
16&-133&739&-12&729\\
17&-91&754&-10&729\\
\end{longtable}
\endgroup}
{\small\begingroup
\setlength{\tabcolsep}{\arraycolsep}
\setlength{\LTpre}{\medskipamount}
\setlength{\LTpost}{\medskipamount}
\begin{longtable}{>{$}c<{$}|>{$}r<{$}>{$}r<{$}|>{$}r<{$}>{$}r<{$}}
\multicolumn{5}{c}{$\text{Dimension 4: auxiliary coefficients and remainders.}$}\\[4pt]
k&a_k^W&b_k^W&r_{W,-}\left(k-2\right)&r_{W,+}\left(k-2\right)\\
\hline
\endfirsthead
\multicolumn{5}{c}{$\text{Dimension 4: auxiliary coefficients and remainders.}$}\\[4pt]
k&a_k^W&b_k^W&r_{W,-}\left(k-2\right)&r_{W,+}\left(k-2\right)\\
\hline
\endhead
\hline
\endfoot
\hline
\endlastfoot
0&10000000000000000&10000000000000000&--&--\\
1&-16094454914005000&-16094454914000000&--&--\\
2&-16094454914005000&-16094454914000000&0&0\\
3&-16094454914005000&-16094454914000000&0&0\\
4&-4323504129287605&-4323504129281484&320000000000000&144000000000000\\
5&-1049874392449287&-1049874392441675&630183680664568&353480517911206\\
6&-130434812554448&-130434812545569&785275520996852&650220776866809\\
7&-14969503199865&-14969503189965&8731527222852&534665834266809\\
8&-1073027360167&-1073027349404&1495736361443100&843861526195924\\
9&-72609487748&-72609476243&2024976463363924&1006722539380242\\
10&-3390313237&-3390301089&4583587643392&1761991797740373\\
11&-151984225&-151971515&2823726063820552&607484764735687\\
12&-5003432&-4990227&4122260030341100&3016889193371415\\
13&-171111&-157469&908484151645360&1933774711514760\\
14&-15856&-1825&1426237224185676&1827189473174578\\
15&-12364&2016&2286170623177308&6582516496125810\\
16&-12364&2149&0&3073838669920677\\
17&-12364&2189&0&7348599359573986\\
\end{longtable}
\endgroup}
{\small\begingroup
\setlength{\tabcolsep}{\arraycolsep}
\setlength{\LTpre}{\medskipamount}
\setlength{\LTpost}{\medskipamount}
\begin{longtable}{>{$}c<{$}|>{$}r<{$}>{$}r<{$}|>{$}r<{$}>{$}r<{$}}
\multicolumn{5}{c}{$\text{Dimension 4: profile division remainders.}$}\\[4pt]
j&r_{u,-}&r_{u,+}&r_{v,-}&r_{v,+}\\
\hline
\endfirsthead
\multicolumn{5}{c}{$\text{Dimension 4: profile division remainders.}$}\\[4pt]
j&r_{u,-}&r_{u,+}&r_{v,-}&r_{v,+}\\
\hline
\endhead
\hline
\endfoot
\hline
\endlastfoot
0&0&0&0&0\\
1&8&8&8000000000000&0\\
2&8&8&32000000000000&0\\
3&0&0&61213578855000&13892430930959\\
4&20&36&51820368282500&13497661075377\\
5&60&68&63996411707752&39929949423750\\
6&16&84&94051097247504&165377585871664\\
7&60&84&69611296998412&207846191126232\\
8&152&176&93922457318684&226029212556300\\
9&108&72&108753462384940&81182446297507\\
10&164&84&519740121338148&490918531175421\\
11&92&236&402766663301464&588705041614482\\
12&84&240&204522314955160&227808301607277\\
13&156&296&701089948740372&230827460522157\\
14&416&100&728751578038860&0\\
15&48&96&176843401891604&0\\
\end{longtable}
\endgroup}

{\small\begingroup
\setlength{\tabcolsep}{\arraycolsep}
\setlength{\LTpre}{\medskipamount}
\setlength{\LTpost}{\medskipamount}
\begin{longtable}{>{$}c<{$}|>{$}r<{$}>{$}r<{$}|>{$}r<{$}>{$}r<{$}}
\multicolumn{5}{c}{$\text{Dimension 5.}$}\\[4pt]
k&a_k^u&b_k^u&a_k^v&b_k^v\\
\hline
\endfirsthead
\multicolumn{5}{c}{$\text{Dimension 5.}$}\\[4pt]
k&a_k^u&b_k^u&a_k^v&b_k^v\\
\hline
\endhead
\hline
\endfoot
\hline
\endlastfoot
0&0&0&0&0\\
1&0&0&10000000000000000&10000000000000000\\
2&10000000000000000&10000000000000000&12500000000000000&12500000000000000\\
3&5830747537695000&5830747537696667&6170529768155416&6170529768155834\\
4&2912141524299583&2912141524302501&1683843760938645&1683843760939376\\
5&650141057804437&650141057808377&294884492055173&294884492056161\\
6&126150673500376&126150673505291&36141263821769&36141263822978\\
7&14923347716842&14923347722737&3282706864956&3282706866372\\
8&1579374879579&1579374886444&230204403230&230204404844\\
9&115965738067&115965745878&12853803169&12853804972\\
10&7791260145&7791268873&585397075&585399060\\
11&390147973&390157589&22170380&22172539\\
12&18178552&18189026&708476&710802\\
13&656402&667706&18289&20775\\
14&17557&29664&-772&1868\\
15&-4810&8074&-500&1491\\
16&-773&7903&-40&1538\\
17&-152&8312&-28&1584\\
\end{longtable}
\endgroup}
{\small\begingroup
\setlength{\tabcolsep}{\arraycolsep}
\setlength{\LTpre}{\medskipamount}
\setlength{\LTpost}{\medskipamount}
\begin{longtable}{>{$}c<{$}|>{$}r<{$}>{$}r<{$}|>{$}r<{$}>{$}r<{$}}
\multicolumn{5}{c}{$\text{Dimension 5: auxiliary coefficients and remainders.}$}\\[4pt]
k&a_k^W&b_k^W&r_{W,-}\left(k-2\right)&r_{W,+}\left(k-2\right)\\
\hline
\endfirsthead
\multicolumn{5}{c}{$\text{Dimension 5: auxiliary coefficients and remainders.}$}\\[4pt]
k&a_k^W&b_k^W&r_{W,-}\left(k-2\right)&r_{W,+}\left(k-2\right)\\
\hline
\endhead
\hline
\endfoot
\hline
\endlastfoot
0&10000000000000000&10000000000000000&--&--\\
1&-17507757386915000&-17507757386910000&--&--\\
2&-21884696733643750&-21884696733637500&0&0\\
3&-25532146189251042&-25532146189243750&64000000000000&0\\
4&-8672532086816459&-8672532086806770&112000000000000&320000000000000\\
5&-2524999851704773&-2524999851691589&589820656635000&232443130419350\\
6&-399747097947814&-399747097931017&231864988030619&782624026883050\\
7&-55841478853556&-55841478833274&320269559679241&293566172184850\\
8&-5128774417138&-5128774393472&943634634451768&39256841842550\\
9&-426539319788&-426539292836&1553429956061906&357205512782850\\
10&-25622888686&-25622858550&2160141403003047&782521772874200\\
11&-1421263281&-1421230065&3424307639240831&2224314710817900\\
12&-60268144&-60231952&3943936325254485&422337837177650\\
13&-2409930&-2370861&4334844022632889&4816278884341450\\
14&-92699&-50849&2434213880478936&3160261607159300\\
15&-20580&23961&1314326945048986&4138458450533300\\
16&-19463&27684&488875277915401&1803233110795200\\
17&-20036&29213&4864000000000000&4626836304495700\\
\end{longtable}
\endgroup}
{\small\begingroup
\setlength{\tabcolsep}{\arraycolsep}
\setlength{\LTpre}{\medskipamount}
\setlength{\LTpost}{\medskipamount}
\begin{longtable}{>{$}c<{$}|>{$}r<{$}>{$}r<{$}|>{$}r<{$}>{$}r<{$}}
\multicolumn{5}{c}{$\text{Dimension 5: profile division remainders.}$}\\[4pt]
j&r_{u,-}&r_{u,+}&r_{v,-}&r_{v,+}\\
\hline
\endfirsthead
\multicolumn{5}{c}{$\text{Dimension 5: profile division remainders.}$}\\[4pt]
j&r_{u,-}&r_{u,+}&r_{v,-}&r_{v,+}\\
\hline
\endhead
\hline
\endfoot
\hline
\endlastfoot
0&0&0&0&0\\
1&0&4&16000000000000&16000000000000\\
2&8&15&4000000000000&12000000000000\\
3&4&36&46787201754882&54958493503216\\
4&9&35&112326057830448&84407606482770\\
5&20&74&11263458503253&13412235691698\\
6&94&93&44046787820394&37221123637594\\
7&136&40&11867649227156&2778008007256\\
8&11&86&182967193867612&294966711655980\\
9&138&134&193096244844956&263007314113994\\
10&103&89&154352477241380&445050752853950\\
11&300&108&395644208641847&28572815861880\\
12&230&342&311760634575546&200604742466776\\
13&218&348&537614240229075&284234568475114\\
14&244&206&683719056594364&390000000000000\\
15&368&204&359623722367343&832000000000000\\
\end{longtable}
\endgroup}

The interval recurrences and the integer remainder identities above
certify the displayed enclosures for the profile coefficients
$u_j,v_j,W_j$. These are precisely the inputs used in the finite
coefficient formulas below. For every subsequent profile index needed
in those formulas, we use
\[
 0\le u_j,v_j\le10^{-16}.
\]
The finite-index bounds follow by induction from the outward-rounded
recurrences. The common bound for $j\ge18$ follows from the
positive-series estimate proved in Appendix~\ref{subsec:appB-profiles}.

\subsection{Finite coefficient formulas and their rational inputs}\label{subsec:appA-formulas}
The profile intervals now give finite bounds for the numerator
coefficients. We first write their coefficient formulas and then
obtain determinant lower bounds by retaining selected positive terms. Put
\[
 \mathcal C=4H^2N_a,\qquad
 \mathcal E=N_e,\qquad \mathcal Z=N_d,
 \qquad
 c_j=\left[\mathcal C\right]_j,\quad e_j=\left[\mathcal E\right]_j,\quad
 z_j=\left[\mathcal Z\right]_j.
\]
The calligraphic symbols distinguish these series from the radial coefficients. Their defining relation is
\begin{equation}\label{eq:CEZ-paperoutput}
 N_{\det}=\mathcal C\mathcal E-\mathcal Z^2.
\end{equation}

All sequences below are zero at negative indices. For two sequences
write
\[
 \left(a*b\right)_j=\sum_{i=0}^j a_i b_{j-i}.
\]
For nonnegative integers $r,s$ define
\[
 \mathcal F_{r,s}\left(a\right)_j
 =\sum_{i=0}^j
 \frac{\left(i+r\right)!}{i!}
 \frac{\left(j-i+s\right)!}{\left(j-i\right)!}
 a_{i+r}a_{j-i+s},
 \qquad j\ge0.
\]
Here $a\left(y\right)=\sum_{i\ge0}a_i y^i$ denotes the series with
coefficient sequence $a$. Termwise differentiation gives
\[
 \partial_y^r a=\sum_{i\ge0}\frac{\left(i+r\right)!}{i!}a_{i+r}y^i,
\]
and multiplying the two differentiated series shows that
$\mathcal F_{r,s}\left(a\right)_j$ is the coefficient of $y^j$ in
$\left(\partial_y^r a\right)\left(\partial_y^s a\right)$.
For example,
\[
 \mathcal F_{1,2}\left(u\right)_j
 =\sum_{i=0}^j\left(i+1\right)\left(j-i+1\right)
              \left(j-i+2\right)u_{i+1}u_{j-i+2}.
\]
The convention at negative indices makes the same formulas applicable
when a factor $y$ or $y^2$ shifts the required degree below zero.
Set
\[
 h_j=u_{j+2},\qquad k_j=v_{j+1},\qquad H_j=\left(j+2\right)h_j.
\]
The formula for $H_j$ follows directly from $H=2h+yh_y$.
Write $\rho_j=\left[D\right]_j$ and $\omega_j=\left[N_g\right]_j$. Then
\[
\begin{aligned}
 \rho_j={}&\sum_{a+b=j}\left[\left(a-b\right)^2-j-4\right]h_a h_b\\
 &+\sum_{a+b=j-1}
       \left[\frac{j+3}{2}-\left(a-b\right)^2\right]h_a h_b,\\
 \omega_j={}&\sum_{a+b=j}\left(a-b-1\right)k_a h_b.
\end{aligned}
\]
For the first identity, $H_b=\left(b+2\right)h_b$ gives
\[
 \left[hH\right]_j
 =\sum_{a+b=j}\left(b+2\right)h_a h_b
 =\frac{j+4}{2}\sum_{a+b=j}h_a h_b.
\]
Since $\partial_x=-\partial_y$, the corresponding differentiated
product satisfies
\[
\begin{aligned}
 \left[y\left(h_xH-hH_x\right)\right]_j
 &=\sum_{a+b=j}\left(b-a\right)\left(b+2\right)h_a h_b\\
 &=\frac12\sum_{a+b=j}\left(a-b\right)^2h_a h_b.
\end{aligned}
\]
The second equality in each case averages the summand with the summand
obtained by interchanging $a,b$. Multiplication by $2\left(1-y\right)$
in $D$ subtracts twice the preceding coefficient, and the remaining
term $-yhH$ subtracts it once more with its $hH$ weight. The weights
at degrees $j$ and $j-1$ are therefore exactly those displayed for
$\rho_j$. For $N_g$, differentiation of $h_b y^b$ contributes $-b$,
differentiation of $k_a y^a$ contributes $-a$, and the undifferentiated
term contributes $-1$. Thus $y\left(kh_x-hk_x\right)-kh$ has the
weight $a-b-1$ in the formula for $\omega_j$.

The following expression is $\beta_j=\left[B\right]_j$:
\[
\begin{aligned}
 \beta_j={}&\mathcal F_{1,1}\left(u\right)_j
 +\left(2n+6\right)
   \left[\mathcal F_{1,2}\left(u\right)_j
              -\mathcal F_{1,2}\left(u\right)_{j-1}\right]\\
 &+4\left[\mathcal F_{2,2}\left(u\right)_j
       -2\mathcal F_{2,2}\left(u\right)_{j-1}
       +\mathcal F_{2,2}\left(u\right)_{j-2}\right]\\
 &-8\left[\mathcal F_{1,3}\left(u\right)_j
       -2\mathcal F_{1,3}\left(u\right)_{j-1}
       +\mathcal F_{1,3}\left(u\right)_{j-2}\right]\\
 &-m\left(n+2\right)\mathcal F_{0,2}\left(u\right)_j
 +2m\left[\mathcal F_{0,3}\left(u\right)_j
              -\mathcal F_{0,3}\left(u\right)_{j-1}\right].
\end{aligned}
\]
The signs can be read term by term from
$u_x=-u_y$, $u_{xx}=u_{yy}$, and $u_{xxx}=-u_{yyy}$.
For example, $-\left(2n+6\right)xu_xu_{xx}$ becomes
$\left(2n+6\right)\left(1-y\right)u_yu_{yy}$, whereas
$-8x^2u_xu_{xxx}$ becomes
$-8\left(1-y\right)^2u_yu_{yyy}$. For any series $A$,
\[
 \left[\left(1-y\right)A\right]_j=A_j-A_{j-1},\qquad
 \left[\left(1-y\right)^2A\right]_j=A_j-2A_{j-1}+A_{j-2}.
\]
These two identities give each shifted difference in $\beta_j$.
Writing $\ell_j=\left[L\right]_j$, the exact factor $N_G=y^3L$
gives $\ell_j=\left[N_G\right]_{j+3}$. Consequently,
\[
\begin{aligned}
 \ell_j=T\left\{\right.&\left.4\left[\mathcal F_{1,1}\left(v\right)_{j+3}
                  -\mathcal F_{1,1}\left(v\right)_{j+2}\right]
 -2m\mathcal F_{0,1}\left(v\right)_{j+3}\right.\\
 &\left.+Q\mathcal F_{0,0}\left(v\right)_{j+3}
 -\frac{T^2}{16}\mathcal F_{0,0}\left(u\right)_{j+3}
 -\frac14\beta_{j+3}\right\}.
\end{aligned}
\]
Finally define
\[
\begin{gathered}
 p_j=\left(\rho*\rho\right)_j
       -16\left(\omega*\omega\right)_j
       +16\left(\omega*\omega\right)_{j-1},\\
 J_j=\sum_{a+b=j}\left(-a+2b-1\right)\ell_a H_b.
\end{gathered}
\]
The three coefficient families are
\begin{equation}\label{eq:coefficients-paperoutput}
\begin{aligned}
 c_j&=\frac T4\sum_{i=0}^j\left(H*H\right)_i p_{j-i+2},\\
 e_j&=\frac{mT}{4}\left(H*H\right)_j
                    +\left(m-1\right)\ell_{j-1},\\
 z_j&=\left(\ell*\rho\right)_j-\left(h*J\right)_j
                              +\left(h*J\right)_{j-1}.
\end{aligned}
\end{equation}
These formulas follow by expanding the defining numerators. For example,
$p_j$ is the coefficient of $D^2-16\left(1-y\right)N_g^2$, and the
expression defining $J_j$ is the coefficient of
$y\left(L_xH-2LH_x\right)-LH$: differentiating $\ell_a y^a$
with respect to $x$ supplies $-a$, and differentiating $H_b y^b$
in the second term supplies $2b$. The final term supplies $-1$.
Since $K=T/16$, the exact factor in $N_a$ gives
\[
 \left[N_a\right]_j=\frac{T}{16}p_{j+2}.
\]
Multiplication by $4H^2$ yields the formula for $c_j$.
The identity $N_e=\left(mT/4\right)H^2+\left(m-1\right)yL$
gives $e_j$. Finally, $N_d=LD-\left(1-y\right)hJ$ gives the
three convolution terms in $z_j$. Thus the index shifts and signs in
\eqref{eq:coefficients-paperoutput} follow directly from multiplication
and exact division by the already proved endpoint factors.

For $j\le20$, the formulas for $c_j,e_j$ involve at most
$u_{24},v_{23}$, and that for $z_j$ involves at most $u_{25},v_{24}$.
The apparent extra indices in $\mathcal F$ have zero factors because
$u_0=u_1=v_0=0$. Thus the displayed input intervals through index $17$
and the interval $\left[0,10^{-16}\right]$ for every subsequent required coefficient
supply every input to \eqref{eq:coefficients-paperoutput}.

All numbers in the following tables are rational. The columns mean
\[
 c_i^-<c_i<c_i^+,\qquad \left|e_i\right|<e_i^+,\qquad \left|z_i\right|<d_i^+.
\]
All operations after the input rounding are exact interval operations.
In the formula for $\ell_j$, evaluate $T^2$ by multiplying
$\mathbf T$ by itself; $\mathbf S$ is used only to construct the
profile inputs.
The bounds below follow by substituting the input intervals in the
finite sums and applying interval addition and multiplication.

\subsection{Positive lower bounds for the determinant coefficients}\label{subsec:appA-positive}
For the lower coefficients, the following strict rational bounds
follow by substitution in these sums. The entries are rounded down
only after the interval calculation.
\begingroup
\setlength{\tabcolsep}{\arraycolsep}
\setlength{\LTpre}{\medskipamount}
\setlength{\LTpost}{\medskipamount}
\begin{longtable}{>{$}c<{$}|>{$}r<{$}>{$}r<{$}|>{$}r<{$}>{$}r<{$}}
&\multicolumn{2}{c|}{$n=4$}&\multicolumn{2}{c}{$n=5$}\\
 j&\left[-L\right]_j&\left[N_e\right]_j&\left[-L\right]_j&\left[N_e\right]_j\\
\hline
\endfirsthead
&\multicolumn{2}{c|}{$n=4$}&\multicolumn{2}{c}{$n=5$}\\
 j&\left[-L\right]_j&\left[N_e\right]_j&\left[-L\right]_j&\left[N_e\right]_j\\
\hline
\endhead
\hline
\endfoot
\hline
\endlastfoot
0&12.957053&63.781193&27.776916&110.992863\\
1&21.466191&62.777118&56.474191&110.820660\\
2&12.715447&37.720177&42.745647&44.771785\\
3&5.294963&20.156767&22.312929&20.922841\\
4&1.567725&8.859678&8.321330&10.670057\\
5&0.367442&2.834257&2.439149&4.557487\\
6&0.067182&0.739224&0.562468&1.662088\\
7&0.010245&0.152611&0.107647&0.487614\\
8&0.001293&0.026355&0.017143&0.117154\\
\end{longtable}
\endgroup
For $0\le j\le5$ the direct convolution
$\left[N_{\det}\right]_j=\sum_{i=0}^j
\left(c_i e_{j-i}-z_i z_{j-i}\right)$ gives
\begingroup
\setlength{\tabcolsep}{\arraycolsep}
\setlength{\LTpre}{\medskipamount}
\setlength{\LTpost}{\medskipamount}
\begin{longtable}{>{$}c<{$}|>{$}r<{$}>{$}r<{$}}
j&n=4&n=5\\
\hline
\endfirsthead
j&n=4&n=5\\
\hline
\endhead
\hline
\endfoot
\hline
\endlastfoot
0&3551.919&5440.203\\
1&22896.369&50823.663\\
2&62224.434&170414.718\\
3&109706.562&349811.383\\
4&145118.069&517714.754\\
5&155600.390&600348.054\\
\end{longtable}
\endgroup
The constant terms are
\[
\begin{gathered}
 \left[-L\right]_0=\frac{mT\left(Z-2\right)}6,\qquad
 \left[N_e\right]_0=mT,\\
 \left[N_{\det}\right]_0=
 \frac{mT^2}{9}\left[\left(12-m\right)Z^2+16mZ-16m-12Q\right].
\end{gathered}
\]
For the first two identities, $H_0=2$ gives
$\left[N_e\right]_0=\left(mT/4\right)H_0^2=mT$, and the
coefficient of $y^3$ in $N_G$ gives
$\ell_0=-mT\left(Z-2\right)/6$. These provide short checks on the
normalization and the sign of $L$. For example, the four-dimensional
parameter intervals imply
$\left[-L\right]_0>3\cdot21.26\cdot1.218/6>12$.
At a nonconstant degree, the same calculation uses the convolution
instead of evaluation at zero; for instance,
\[
 e_1=\frac{mT}{2}H_0H_1+\left(m-1\right)\ell_0
     =3mT u_3+\left(m-1\right)\ell_0.
\]
Each factor on the right already has an interval enclosure. A lower
bound uses the lower endpoint after summing both terms, including the
negative contribution from $\ell_0$.
The preceding two tables prove the first twenty-four inequalities
required in Section~\ref{sec:certificate}. For the higher determinant
coefficients and the low-block remainder, the following enclosures
keep the factors $\mathcal C,\mathcal E,\mathcal Z$ separate.

{\small
\begingroup
\setlength{\LTpre}{\medskipamount}
\setlength{\LTpost}{\medskipamount}
\begin{longtable}{r|rrrr}
\multicolumn{5}{c}{\textbf{Dimension four.}}\\[4pt]
$i$&$c_i^-$&$c_i^+$&$e_i^+$&$d_i^+$\\
\hline
\endfirsthead
\multicolumn{5}{c}{\textbf{Dimension four.}}\\[4pt]
$i$&$c_i^-$&$c_i^+$&$e_i^+$&$d_i^+$\\
\hline
\endhead
\hline
\endfoot
\hline
\endlastfoot
0&66.2&66.3&64&26\\
1&337&338&63&54.1\\
2&726&727&38&93.8\\
3&1023&1024&21&98.2\\
4&1085&1086&9&71.9\\
5&928&929&3&38.7\\
6&667&668&0.75&16\\
7&413&414&0.16&5.26\\
8&223&224&0.027&1.41\\
9&106&107&0.004&0.305\\
10&45.4&45.5&0.0005&0.0507\\
11&17.3&17.4&0.00005&0.00517\\
12&5.96&5.97&0.00001&0.000376\\
13&1.86&1.87&0.00001&0.000405\\
14&0.53&0.531&0.00001&0.00022\\
15&0.138&0.139&0.00001&0.000183\\
16&0.0332&0.0333&0.00001&0.000211\\
17&0.00738&0.00739&0.00001&0.000263\\
18&0.00151&0.00153&0.00001&0.000327\\
19&0.000288&0.000293&0.00001&0.000403\\
20&0.00005&0.000054&0.00001&0.000492\\
\end{longtable}
\endgroup
}
{\small
\begingroup
\setlength{\LTpre}{\medskipamount}
\setlength{\LTpost}{\medskipamount}
\begin{longtable}{r|rrrr}
\multicolumn{5}{c}{\textbf{Dimension five.}}\\[4pt]
$i$&$c_i^-$&$c_i^+$&$e_i^+$&$d_i^+$\\
\hline
\endfirsthead
\multicolumn{5}{c}{\textbf{Dimension five.}}\\[4pt]
$i$&$c_i^-$&$c_i^+$&$e_i^+$&$d_i^+$\\
\hline
\endhead
\hline
\endfoot
\hline
\endlastfoot
0&76.8&76.9&111&55.6\\
1&527&528&111&146\\
2&1462&1463&45&293\\
3&2610&2611&21&377\\
4&3469&3470&11&347\\
5&3698&3699&4.6&240\\
6&3299&3300&1.7&130\\
7&2528&2529&0.49&57.1\\
8&1691&1692&0.12&21\\
9&1000&1001&0.023&6.53\\
10&527&528&0.004&1.73\\
11&249&250&0.0006&0.383\\
12&106&107&0.0001&0.0691\\
13&41.6&41.7&0.0001&0.00944\\
14&14.7&14.8&0.0001&0.000975\\
15&4.81&4.82&0.0001&0.000962\\
16&1.44&1.45&0.0001&0.00109\\
17&0.401&0.402&0.0001&0.00134\\
18&0.103&0.104&0.0001&0.00169\\
19&0.0248&0.0249&0.0001&0.00214\\
20&0.00554&0.00558&0.0001&0.00269\\
\end{longtable}
\endgroup
}

The lower bounds for the first nine coefficients of $\mathcal E$ are
\begingroup
\setlength{\tabcolsep}{\arraycolsep}
\setlength{\LTpre}{\medskipamount}
\setlength{\LTpost}{\medskipamount}
\begin{longtable}{>{$}c<{$}|>{$}r<{$}>{$}r<{$}>{$}r<{$}>{$}r<{$}>{$}r<{$}>{$}r<{$}>{$}r<{$}>{$}r<{$}>{$}r<{$}}
n\backslash j&0&1&2&3&4&5&6&7&8\\
\hline
\endfirsthead
n\backslash j&0&1&2&3&4&5&6&7&8\\
\hline
\endhead
\hline
\endfoot
\hline
\endlastfoot
4&63&62&37&20&8&2.8&.7&.15&.02\\
 5&110&110&44&20&10&4.5&1.6&.48&.11\\
\end{longtable}
\endgroup
Denote these bounds by $e_{n,j}^-$. They follow directly from the
coefficient sums above, with the same rational interval inputs.
For the remaining indices needed here, the same finite sums give
\[
 e_j>-\delta_n\quad\left(9\le j\le20\right),\qquad
 \delta_4=10^{-5},\quad\delta_5=10^{-4}.
\]
The sums of the $c_i^+$ columns are
\[
 \sum_{i=0}^{20}c_i^+
 =5654.752567<6000\quad\left(n=4\right),\qquad
 \sum_{i=0}^{20}c_i^+
 =21318.20648<22000\quad\left(n=5\right).
\]
Write $U_4=6000$ and $U_5=22000$.

For an index $j$ between $6$ and $20$, choose a set $I_j$ of indices
$i$ with $0\le j-i\le8$, as specified below. Since every $c_i$ in use
is positive, \eqref{eq:CEZ-paperoutput} implies
\begin{equation}\label{eq:sparse-paperoutput}
 \left[N_{\det}\right]_j>
 \sum_{i\in I_j}c_i^-e_{n,j-i}^-
 -\mathcal Q_j-\delta_n U_n,
\end{equation}
where
\[
 \mathcal Q_j=
 \begin{cases}
  \left(d_{j/2}^+\right)^2,&j\text{ even},\\
  0,&j\text{ odd}
 \end{cases}
 +2\sum_{i=0}^{\left\lfloor\left(j-1\right)/2\right\rfloor}
              d_i^+d_{j-i}^+.
\]
Indeed, the unselected terms whose $\mathcal E$ index is at most eight
are nonnegative and may be discarded. All terms with a larger
$\mathcal E$ index together contribute more than
$-\delta_n\sum c_i>-\delta_n U_n$. Pairing symmetric indices bounds
$\left[\mathcal Z^2\right]_j$ above by $\mathcal Q_j$.
Thus \eqref{eq:sparse-paperoutput} is a sufficient lower bound that
uses only a small part of the positive convolution. More explicitly,
split $\sum_{i=0}^j c_i e_{j-i}$ according to $j-i\le8$ and
$j-i\ge9$. In the first part all summands are positive; replacing
the selected factors by $c_i^-,e_{n,j-i}^-$ decreases the sum. In the
second part, $e_{j-i}>-\delta_n$ and $0<c_i<c_i^+$ imply
$c_i e_{j-i}>-\delta_n c_i^+$. This argument does not require those
high-index $e_{j-i}$ to have either sign. Similarly,
$z_i z_{j-i}\le\left|z_i\right|\left|z_{j-i}\right|$ bounds the
subtracted square from above. If $j$ is even its central product
occurs once; every other pair occurs twice, which accounts for the
two parts of $\mathcal Q_j$.

In the next two tables, $A_j$ is the selected positive sum in
\eqref{eq:sparse-paperoutput}, and $B_j$ is a strict upper bound for
$\mathcal Q_j$. Every row requires at most eleven products to check
$\mathcal Q_j<B_j$.
{\small
\begingroup
\setlength{\LTpre}{\medskipamount}
\setlength{\LTpost}{\medskipamount}
\begin{longtable}{r|crr}
\multicolumn{4}{c}{\textbf{Dimension four.}}\\[4pt]
$j$&$I_j$&$A_j$&$B_j$\\
\hline
\endfirsthead
\multicolumn{4}{c}{\textbf{Dimension four.}}\\[4pt]
$j$&$I_j$&$A_j$&$B_j$\\
\hline
\endhead
\hline
\endfoot
\hline
\endlastfoot
6&$\left\{5\right\}$&57536&28200\\
7&$\left\{6\right\}$&41354&23400\\
8&$\left\{7\right\}$&25606&16500\\
9&$\left\{7\right\}$&15281&9870\\
10&$\left\{8\right\}$&8251&5140\\
11&$\left\{8\right\}$&4460&2340\\
12&$\left\{9\right\}$&2120&936\\
13&$\left\{9\right\}$&848&333\\
14&$\left\{10\right\}$&363.2&105\\
15&$\left\{11\right\}$&138.4&29.5\\
16&$\left\{12\right\}$&47.68&7.43\\
17&$\left\{12\right\}$&16.688&1.77\\
18&$\left\{13\right\}$&5.208&0.487\\
19&$\left\{14\right\}$&1.484&0.253\\
20&$\left\{15\right\}$&0.3864&0.245\\
\end{longtable}
\endgroup
}
{\small
\begingroup
\setlength{\LTpre}{\medskipamount}
\setlength{\LTpost}{\medskipamount}
\begin{longtable}{r|crr}
\multicolumn{4}{c}{\textbf{Dimension five.}}\\[4pt]
$j$&$I_j$&$A_j$&$B_j$\\
\hline
\endfirsthead
\multicolumn{4}{c}{\textbf{Dimension five.}}\\[4pt]
$j$&$I_j$&$A_j$&$B_j$\\
\hline
\endhead
\hline
\endfoot
\hline
\endlastfoot
6&$\left\{5,6\right\}$&769670&431000\\
7&$\left\{6,7\right\}$&640970&447000\\
8&$\left\{7,8\right\}$&464090&397000\\
9&$\left\{7,8,9\right\}$&407242&305000\\
10&$\left\{8,9,10\right\}$&242374&206000\\
11&$\left\{8,9,10\right\}$&135790&123000\\
12&$\left\{9,10,11\right\}$&70578&65000\\
13&$\left\{10,11,12\right\}$&33156&31100\\
14&$\left\{10,11,12\right\}$&14914&13400\\
15&$\left\{10,11,12\right\}$&6981.5&5260\\
16&$\left\{11,12\right\}$&2180.5&1880\\
17&$\left\{12,13\right\}$&893&613\\
18&$\left\{13,14\right\}$&334.2&185\\
19&$\left\{13\right\}$&66.56&52.5\\
20&$\left\{14\right\}$&23.52&16.4\\
\end{longtable}
\endgroup
}

Every row satisfies
$A_j-B_j-\delta_nU_n>1/20$. The smallest four-dimensional row can
be checked in full as follows:
\[
\begin{aligned}
 \mathcal Q_{20}={}&0.0507^2
 +2\left(26\cdot0.000492+54.1\cdot0.000403
             +93.8\cdot0.000327+98.2\cdot0.000263\right.\\
 &\left.\quad+71.9\cdot0.000211+38.7\cdot0.000183
             +16\cdot0.000220+5.26\cdot0.000405\right.\\
 &\left.\quad+1.41\cdot0.000376+0.305\cdot0.00517\right)\\
 ={}&0.24477811<0.245.
\end{aligned}
\]
The selected index is $15$, and therefore
\[
 \left[N_{\det}\right]_{20}>0.138\cdot2.8-0.245-0.06
 =0.0814>\frac1{20}.
\]
For dimension five the corresponding check is
\[
\begin{aligned}
 \mathcal Q_{20}={}&1.73^2
 +2\left(55.6\cdot0.00269+146\cdot0.00214
             +293\cdot0.00169+377\cdot0.00134\right.\\
 &\left.\quad+347\cdot0.00109+240\cdot0.000962
             +130\cdot0.000975+57.1\cdot0.00944\right.\\
 &\left.\quad+21\cdot0.0691+6.53\cdot0.383\right)\\
 ={}&16.371556<16.4,
\end{aligned}
\]
and
\[
 \left[N_{\det}\right]_{20}>14.7\cdot1.6-16.4-2.2
 =4.92>\frac1{20}.
\]
The remaining rows give the common bounds
\begin{proposition}\label{prop:finite-determinant-coefficients}
For $n=4,5$, the determinant coefficients satisfy
\begin{equation}\label{eq:det-positive-paperoutput}
 \left[N_{\det}\right]_j>\frac1{20},\qquad 6\le j\le20.
\end{equation}
\end{proposition}
\begin{proof}
Fix $n\in\{4,5\}$ and $6\le j\le20$. The sparse-convolution
inequality~\eqref{eq:sparse-paperoutput} gives
\[
 [N_{\det}]_j>A_j-\mathcal Q_j-\delta_nU_n,
\]
where $A_j$ is the selected positive sum listed in the corresponding
dimension table, and $B_j$ is the displayed strict upper bound for
$\mathcal Q_j$. Thus it is enough to check
$A_j-B_j-\delta_nU_n>1/20$ in every row. The error terms are
\[
 \delta_4U_4=10^{-5}\cdot6000=0.06,\qquad
 \delta_5U_5=10^{-4}\cdot22000=2.2.
\]
Substitution of the entries in the two tables gives the required
inequality in every row. The smallest margin in dimension four occurs
at $j=20$. In that row, $I_{20}=\{15\}$, so
\[
 A_{20}=c_{15}^-e_{4,5}^-=0.138\cdot2.8.
\]
The preceding $d_i^+$ table gives
\[
 \begin{aligned}
 \mathcal Q_{20}
 &=0.0507^2+2\bigl(26\cdot0.000492
       +54.1\cdot0.000403+93.8\cdot0.000327\\
 &\quad+98.2\cdot0.000263+71.9\cdot0.000211
       +38.7\cdot0.000183+16\cdot0.000220\\
 &\quad+5.26\cdot0.000405+1.41\cdot0.000376
       +0.305\cdot0.00517\bigr)\\
 &=0.24477811<0.245.
 \end{aligned}
\]
Consequently,
\[
 [N_{\det}]_{20}>
 0.138\cdot2.8-0.245-0.06
 =0.0814>\frac1{20}.
\]
In dimension five, the last row has $I_{20}=\{14\}$ and gives
\[
 [N_{\det}]_{20}>
 14.7\cdot1.6-16.4-2.2=4.92>\frac1{20}.
\]
The other rows have larger margins, as is seen directly from the
tabulated values of $A_j,B_j$ and the corresponding constant
$\delta_nU_n$. This proves the assertion for all $6\le j\le20$.
\end{proof}
Only positivity of these coefficients is needed when the series is
evaluated on $0\le y\le1$; the proposition supplies that property.

\subsection{High-degree products of the finite coefficient blocks}\label{subsec:appA-finite-tail}
To bound the part of the determinant tail formed from these finite
blocks, let $\mathcal C_0,\mathcal E_0,\mathcal Z_0$ denote the
polynomials through degree twenty. The five-dimensional upper columns dominate
the four-dimensional columns term by term. Consequently the displayed
rational inputs give
\[
\begin{aligned}
 \mathcal T_{21}
   \left[\left(\mathcal C_0\mathcal E_0-\mathcal Z_0^2\right)^\#\right]
 &\le\sum_{i=1}^{20}c_i^+
           \sum_{k=21-i}^{20}e_k^+
       +\sum_{i=1}^{20}d_i^+
           \sum_{k=21-i}^{20}d_k^+\\
 &=92.983482190329<100.
\end{aligned}
\]
This is a finite sum of rational products. For instance, its first
row is $528\cdot0.0001+146\cdot0.00269=0.44554$.
The lower limits $21-i$ account for every product of degree at least
twenty-one, even though each factor has degree at most twenty. Indeed,
the total coefficient sum over these degrees is the sum over all
ordered pairs $\left(i,k\right)$ with $0\le i,k\le20$ and
$i+k\ge21$. The index $i=0$ contributes nothing; fixing
$1\le i\le20$ leaves exactly $21-i\le k\le20$.
No factor of two is needed in this row-wise sum, because the two
ordered pairs are already included in their respective rows. For
example, the second row is
\[
 1463\left(0.0001+0.0001\right)=0.2926,\qquad
 293\left(0.00214+0.00269\right)=1.41519.
\]
Every entry below is an upper bound for the corresponding row product.
Thus each row can be checked by a suffix sum and a single multiplication.
{\small
\begingroup
\setlength{\LTpre}{\medskipamount}
\setlength{\LTpost}{\medskipamount}
\begin{longtable}{r|rr}
$i$&$c_i^+\displaystyle\sum_{k=21-i}^{20}e_k^+$&$d_i^+\displaystyle\sum_{k=21-i}^{20}d_k^+$\\
\hline
\endfirsthead
$i$&$c_i^+\displaystyle\sum_{k=21-i}^{20}e_k^+$&$d_i^+\displaystyle\sum_{k=21-i}^{20}d_k^+$\\
\hline
\endhead
\hline
\endfoot
\hline
\endlastfoot
1&0.052800&0.392740\\
2&0.292600&1.415190\\
3&0.783300&2.458040\\
4&1.388000&2.727420\\
5&1.849500&2.148000\\
6&1.980000&1.288560\\
7&1.770300&0.621648\\
8&1.353600&0.426867\\
9&0.900900&0.583959\\
10&0.792000&0.817299\\
11&1.375000&0.843530\\
12&3.049500&0.603411\\
13&6.192450&0.280675\\
14&9.449800&0.084662\\
15&11.271570&0.208593\\
16&10.060825&0.497948\\
17&7.211277&1.077136\\
18&4.049604&1.995607\\
19&2.090069&3.154002\\
20&1.087757&4.357350\\
\hline
sum&67.000852&25.982637\\
\end{longtable}
\endgroup
}

The two columns of rounded upper bounds sum to
$67.000852+25.982637=92.983489<100$, consistently with the exact
sum above. This accounts for all products of the finite blocks whose
total degree is at least twenty-one. Appendix~\ref{sec:tail} bounds
the infinite remainders of the factors.

\section{Complete-series remainder estimates}\label{sec:tail}
We now control the complete Taylor series, including every coefficient
above the finite blocks. The bounds in this appendix concern the actual
reference profiles. All auxiliary majorants have nonnegative
coefficients and are given explicitly as rational functions. Thus their
tails are bounded by positive sums rather than by cancellation between
large quantities. Subsection~\ref{subsec:appB-profiles} constructs the
profile majorants, Subsection~\ref{subsec:appB-quadratic} bounds the
tails of $-L$ and $N_e$, and Subsection~\ref{subsec:appB-determinant}
bounds the tails of the three factors entering $N_{\det}$.

\subsection{Boundary coefficients and positive profile majorants}\label{subsec:appB-profiles}
To obtain these bounds, we start with the boundary coefficients of
the profiles. Put $U_c=u\left(1\right)$, $V_c=v\left(1\right)$, and
$c=B_\nu\left(R\right)/M_\nu\left(R\right)\in\left(-1,0\right)$.
The profile formulas of Section~\ref{sec:profiles} give
\[
 U_c=-\frac{4\left(1-c\right)}{T B_\nu\left(R\right)},\qquad
 V_c=\frac{2n}{Q B_{\nu+1}\left(q\right)}.
\]
To extract boundary coefficients, write the profiles in the form
\[
\begin{aligned}
 v\left(y\right)&=V_c B_\nu\left(q\sqrt{1-y}\right),\\
 u\left(y\right)&=\frac{U_c}{1-c}
 \left[B_\nu\left(R\sqrt{1-y}\right)
                   -cM_\nu\left(R\sqrt{1-y}\right)\right].
\end{aligned}
\]
Expanding $\left(1-y\right)^k$ in the defining entire series gives,
for $a>0$,
\[
\begin{aligned}
 \left[B_\nu\left(a\sqrt{1-y}\right)\right]_j
 &=\frac{\left(a^2/4\right)^j}{j!\left(\nu+1\right)_j}
       B_{\nu+j}\left(a\right),\\
 \left[M_\nu\left(a\sqrt{1-y}\right)\right]_j
 &=\frac{\left(-1\right)^j\left(a^2/4\right)^j}
              {j!\left(\nu+1\right)_j}M_{\nu+j}\left(a\right).
\end{aligned}
\]
For example, set $k=j+\ell$ in the first coefficient sum and use
$\binom{j+\ell}{j}/\left(j+\ell\right)!=1/\left(j!\ell!\right)$
and $\left(\nu+1\right)_{j+\ell}
=\left(\nu+1\right)_j\left(\nu+j+1\right)_\ell$.
The remaining sum is precisely $B_{\nu+j}\left(a\right)$.
Absolute convergence on compact sets justifies collecting these
coefficients. The Taylor coefficients at the boundary are therefore
\[
\begin{aligned}
 v_j&=\frac{V_c\left(Q/4\right)^j B_{\nu+j}\left(q\right)}
             {j!\left(\nu+1\right)_j},\\
 u_j&=\frac{U_c\left(T/4\right)^j
       \left[B_{\nu+j}\left(R\right)
       -\left(-1\right)^jcM_{\nu+j}\left(R\right)\right]}
       {\left(1-c\right)j!\left(\nu+1\right)_j}.
\end{aligned}
\]
These coefficients are nonnegative. For $v_j$, $j\ge1$, this follows
from the increase of Bessel zeros with their order. For even $j\ge2$,
the assertion for $u_j$ follows from $c<0$ and
$R<j_{\nu+j,1}$. For odd $j\ge3$, the canonical product shows that
$I_\mu\left(R\right)/J_\mu\left(R\right)$ decreases as $\mu$ increases
while $R<j_{\mu,1}$; at $\mu=\nu+1$ this quotient equals $-1/c$.
The bracket in $u_j$ is therefore positive. The forced coefficients
$u_0,u_1,v_0$ are zero. These uses of the canonical product and order
monotonicity are standard \cite[Section 10.21]{DLMF}.

We use the amplitude bounds
\[
 U_c<2,\qquad V_c<\frac{10}{3}.
\]
Here is a direct rational verification. For $\gamma>0$ define the
finite alternating polynomial
\[
 \beta_N\left(\gamma,z\right)
 =\sum_{j=0}^{N}\frac{\left(-1\right)^jz^j}
                         {j!\left(\gamma\right)_j}.
\]
Whenever the omitted terms decrease in absolute value, even and odd
partial sums are respectively upper and lower bounds for the
corresponding Bessel series. The ratio that verifies this condition is
$z/\left[\left(j+1\right)\left(\gamma+j\right)\right]$.
Also $B_\nu$ is decreasing on $\left(0,R\right)$, since
\[
 B_\nu'\left(t\right)=-\frac{t}{n}B_{\nu+1}\left(t\right)<0
 \quad\hbox{for }0<t\le R.
\]
The parameter enclosures give $T>21$ when $n=4$ and $T>277/10$ when
$n=5$. Hence
\[
\begin{aligned}
 n=4:\quad B_\nu\left(R\right)
 &<\beta_8\left(2,21/4\right)
   =-\frac{2955988871}{26843545600}<-\frac{11}{100},\\
 n=5:\quad B_\nu\left(R\right)
 &<\beta_8\left(5/2,277/40\right)\\
 &=-\frac{16714385790512438931839}
           {225268704460800000000000}<-\frac{37}{500}.
\end{aligned}
\]
The positive series for $M_\nu$, using only its first four terms,
gives
\[
 M_\nu\left(R\right)>\frac{7093}{1024}>6\quad\left(n=4\right),
 \qquad
 M_\nu\left(R\right)>\frac{119689993}{15120000}>6
       \quad\left(n=5\right).
\]
Since
\[
 U_c=\frac{4}{T\left(-B_\nu\left(R\right)\right)}
          +\frac{4}{T M_\nu\left(R\right)},
\]
we obtain
\[
\begin{aligned}
 n=4:\quad U_c
 &<\frac{400}{231}+\frac{2}{63}
 <\frac74+\frac1{30}=\frac{107}{60}<2,\\
 n=5:\quad U_c
 &<\frac{20000}{10249}+\frac{20}{831}
 <\frac{49}{25}+\frac1{40}=\frac{397}{200}<2.
\end{aligned}
\]

For the membrane amplitude, the parameter boxes give $14<Q<15$ in
dimension four and $20<Q<41/2$ in dimension five. Odd partial sums
give
\[
\begin{aligned}
 n=4:\quad B_{\nu+1}\left(q\right)
 &>\beta_3\left(3,15/4\right)=\frac{97}{512}>\frac3{16},\\
 n=5:\quad B_{\nu+1}\left(q\right)
 &>\beta_5\left(7/2,41/8\right)
   =\frac{2567617699}{16605388800}>\frac3{20}.
\end{aligned}
\]
The monotonicity of the finite polynomials used here can be checked
without any further Bessel estimate. Indeed,
\[
 \partial_z\beta_N\left(\gamma,z\right)
 =-\frac{1}{\gamma}\beta_{N-1}\left(\gamma+1,z\right).
\]
For $\gamma>1$, the quadratic $\beta_2\left(\gamma,z\right)$ is
strictly positive on the real line: its leading coefficient is
positive and its discriminant is
$\left(1-\gamma\right)/\left[\gamma^2\left(\gamma+1\right)\right]<0$.
This proves the required monotonicity of $\beta_3\left(3,z\right)$.
For the second polynomial, $\beta_3\left(9/2,z\right)$ is decreasing,
and
\[
 \beta_3\left(9/2,41/8\right)=\frac{124675}{494208}>0.
\]
Thus $\beta_4\left(9/2,z\right)>0$ on the needed interval, which
proves that $\beta_5\left(7/2,z\right)$ decreases there. We conclude
that
\[
 V_c<\frac{64}{21}<\frac{10}{3}\quad\left(n=4\right),\qquad
 V_c<\frac{10}{3}\quad\left(n=5\right).
\]

For power series, $A\preccurlyeq B$ will mean that every coefficient
of $A$ is nonnegative and no larger than the corresponding coefficient
of $B$. The absolute-coefficient series $A^\#$ and its norm were defined in
Section~\ref{sec:certificate}.
The coefficient order is preserved by addition, multiplication,
differentiation in $y$, and the deletion of any initial coefficients.
For multiplication this follows from the triangle inequality in each
finite convolution; differentiation multiplies the coefficient of
degree $j$ by the nonnegative integer $j$. In particular,
\[
 \left(AB\right)^\#\preccurlyeq A^\# B^\#,
 \qquad
 \norm{AB}_{\mathcal A}\le\norm{A}_{\mathcal A}\norm{B}_{\mathcal A}.
\]
Define the rational polynomials
\[
\begin{aligned}
 P_u\left(y\right)={}&y^2+\frac35y^3+\frac3{10}y^4
 +\frac7{100}y^5+\frac{13}{1000}y^6+\frac3{2000}y^7\\
 &+\frac1{6000}y^8+\frac3{250000}y^9
 +\frac1{1250000}y^{10}+\frac1{25000000}y^{11},\\
 P_v\left(y\right)={}&y+\frac54y^2+\frac58y^3
 +\frac{17}{100}y^4+\frac3{100}y^5+\frac{37}{10000}y^6\\
 &+\frac1{3000}y^7+\frac1{40000}y^8
 +\frac{13}{10^7}y^9+\frac6{10^8}y^{10}
 +\frac{23}{10^{10}}y^{11}.
\end{aligned}
\]
The rational profile enclosures in the finite certificate imply
$u_j\le\left[P_u\right]_j$ and $v_j\le\left[P_v\right]_j$ for
$0\le j\le11$. Put
\[
 \mathfrak a=\frac{10\,7^{12}}{3\,12!\,13!},\qquad
 \vartheta=\frac1{25},\qquad
 U_+=P_u+\frac{\mathfrak a y^{12}}{1-\vartheta y},\qquad
 V_+=P_v+\frac{\mathfrak a y^{12}}{1-\vartheta y}.
\]
The notation $\mathfrak a$ is distinct from the Hessian coefficient
$a=1-g^2$, and $\vartheta$ is distinct from the boundary constant
$\eta$.

We claim that $u\preccurlyeq U_+$ and $v\preccurlyeq V_+$ for the
complete series. For $j\ge12$, the positive modified Bessel series
gives
\[
 M_{\nu+j}\left(R\right)
 \le\exp\left(\frac{7}{j+2}\right)
 \le\sqrt{\mathrm e}<\frac53.
\]
The last elementary bound follows, for example, by summing the
exponential series through degree three and bounding its tail by a
geometric series:
\[
 \mathrm e<\frac83+\frac5{96}=\frac{87}{32}
 <\frac{11}{4}<\frac{25}{9}.
\]
For the first inequality, each factor in
$\left(\nu+j+1\right)_k$ is at least $j+2$, so the $k$th positive
summand of $M_{\nu+j}\left(R\right)$ is at most
$\left[7/\left(j+2\right)\right]^k/k!$. Summing proves the
exponential bound. Moreover, termwise absolute values give
$\left|B_{\nu+j}\left(R\right)\right|\le M_{\nu+j}\left(R\right)$.
As the bracket in $u_j$ is nonnegative, it satisfies
\[
 0\le B_{\nu+j}\left(R\right)
       -\left(-1\right)^jcM_{\nu+j}\left(R\right)
 \le\left(1-c\right)M_{\nu+j}\left(R\right).
\]
Finally, $\left(\nu+1\right)_j\ge\left(j+1\right)!$ because
$\nu\ge1$. The coefficient formula and $U_c<2$ therefore imply
\[
 u_j<\frac{10\,7^j}{3\,j!\left(j+1\right)!},\qquad j\ge12.
\]
For $v_j$, the series for $B_{\nu+j}\left(q\right)$ is alternating
with decreasing absolute terms, and hence its value lies in
$\left(0,1\right]$. The same bound follows from $V_c<10/3$.
Finally,
\[
 \frac{7^{j+1}/\left[\left(j+1\right)!\left(j+2\right)!\right]}
      {7^j/\left[j!\left(j+1\right)!\right]}
 =\frac7{\left(j+1\right)\left(j+2\right)}
 \le\frac7{13\cdot14}<\vartheta\quad\left(j\ge12\right).
\]
Put $b_j=10\,7^j/\left[3j!\left(j+1\right)!\right]$ for this
estimate. The first tail coefficient is $b_{12}=\mathfrak a$, and the
ratio bound implies by induction that
$b_{12+k}\le\mathfrak a\vartheta^k$ for every $k\ge0$.
These are exactly the coefficients of
$\mathfrak a y^{12}/\left(1-\vartheta y\right)$. Combining this tail
estimate with the first twelve coefficient comparisons proves the
claimed rational majorants for all coefficients.

For a series $F$ with nonnegative coefficients, write
\[
 \mathcal T_kF=\sum_{j\ge k}\left[F\right]_j
             =F\left(1\right)-\sum_{j<k}\left[F\right]_j,
 \qquad
 S_kF=\sum_{j\ge0}\left[F\right]_{j+k}y^j.
\]
The operator $S_k$ deletes the first $k$ coefficients. In particular it
is defined even when these coefficients do not vanish.

\subsection{Tails of the quadratic numerators}\label{subsec:appB-quadratic}
We now apply the profile majorants to the quadratic and higher-order
numerators. In the following formulas $U=U_+$, $V=V_+$, $X_+=1+y$, and primes
mean differentiation in $y$. The common parameter bounds are
\[
 T<28,\quad Q<21,\quad K<\frac74,\quad n\le5,\quad m\le4.
\]
They show that $N_G^\#$ is bounded coefficientwise by the following
positive rational function:
\[
\begin{aligned}
 \mathcal N_+={}&112X_+V'^2+224VV'+588V^2+1372U^2+7U'^2\\
 &+112X_+U'U''+28X_+^2U''^2+56X_+^2U'U'''\\
 &+196UU''+56X_+UU'''.
\end{aligned}
\]
Each multiplier in this expression comes from an absolute bound on
the corresponding coefficient in $N_G$. For instance,
$4T<112$, $2mT<224$, $QT<588$, and $KT^2<1372$ give its first
four terms. The factor $4K<7$ multiplying $B$ gives the remaining
six terms: $7\left(2n+6\right)\le112$,
$7\cdot4=28$, $7\cdot8=56$,
$7m\left(n+2\right)\le196$, and $7\cdot2m\le56$.
Replacing each occurrence of $x=1-y$ by $X_+=1+y$ bounds its
absolute coefficients, including the coefficients of $x^2$.
If an entire numerator $F$ has the exact factor $y^k$, then
$\left(F/y^k\right)^\#=S_k\left(F^\#\right)$: division simply
shifts all coefficients by $k$. Applying this fact to $N_G=y^3L$
gives $L^\#\preccurlyeq S_3\mathcal N_+$. Also put
\[
 h_+=S_2U,\qquad H_+=2h_++yh_+'=U'/y.
\]
All coefficients in these functions are nonnegative. The following
table records rational upper bounds for the ten terms of
$\mathcal N_+$ in their displayed order:
\begingroup
\setlength{\tabcolsep}{\arraycolsep}
\setlength{\LTpre}{\medskipamount}
\setlength{\LTpost}{\medskipamount}
\begin{longtable}{>{$}c<{$}|>{$}r<{$}>{$}r<{$}}
\text{term}&\mathcal T_{11}&\mathcal T_{12}\\
\hline
\endfirsthead
\text{term}&\mathcal T_{11}&\mathcal T_{12}\\
\hline
\endhead
\hline
\endfoot
\hline
\endlastfoot
1&0.216&0.026\\
2&0.067&0.008\\
3&0.247&0.029\\
4&4.829&0.790\\
5&0.024&0.004\\
6&2.408&0.366\\
7&3.763&0.611\\
8&6.947&1.141\\
9&0.606&0.087\\
10&1.029&0.160\\\hline
\text{sum}&20.136&3.222\\
\end{longtable}
\endgroup
All entries are strict upper bounds. They can be checked by the
positive rational operations below. We first give two shorter
expansions that also illustrate the tail arithmetic explicitly.

Write $P_u=\sum_{j=2}^{11}p_jy^j$ and
$\delta=\mathfrak a/\left(1-\vartheta\right)$. Grouping the products
according to their smaller index gives the exact identity
\[
\begin{aligned}
 \mathcal T_{12}\left(U^2\right)
 ={}&2\left[p_2\left(p_{10}+p_{11}\right)
       +p_3\left(p_9+p_{10}+p_{11}\right)\right.\\
 &\left.\quad+p_4\left(p_8+p_9+p_{10}+p_{11}\right)
       +p_5\left(p_7+p_8+p_9+p_{10}+p_{11}\right)\right]\\
 &+\left(\sum_{j=6}^{11}p_j\right)^2
   +2\delta P_u\left(1\right)+\delta^2.
\end{aligned}
\]
To justify this grouping, write $U=P_u+R_u$, where
$R_u=\mathfrak a y^{12}/\left(1-\vartheta y\right)$ and
$R_u\left(1\right)=\delta$. Every coefficient of
$2P_uR_u+R_u^2$ has degree at least twelve, so its entire sum is
$2\delta P_u\left(1\right)+\delta^2$. In $P_u^2$, a smaller
index $i\in\left\{2,3,4,5\right\}$ must be paired with an index
at least $12-i$, and the pair has two orders. If both indices are at
least six, all their products belong to the tail and are counted by
the square of their sum. These cases exhaust the terms.
Substituting the displayed rational coefficients of $P_u$ and
multiplying by $1372$ gives the fourth upper bound $.790$ in the
last column.

For the second expansion, put
\[
 b_j=\left(j+2\right)p_{j+2},\qquad
 P_H=\sum_{j=0}^{9}b_jy^j,\qquad
 \delta_H=\mathfrak a\left[\frac{12}{1-\vartheta}
                 +\frac{\vartheta}{\left(1-\vartheta\right)^2}\right].
\]
Indeed, differentiating $R_u$ and dividing by $y$ gives
\[
 H_+-P_H=\mathfrak a\left[
 \frac{12y^{10}}{1-\vartheta y}
 +\frac{\vartheta y^{11}}{\left(1-\vartheta y\right)^2}\right].
\]
This remainder starts at degree ten and its value at $y=1$ is
$\delta_H$. Every product involving it is therefore in the tail
starting at degree nine. In $P_H^2$, products with smaller index at
most four are counted in pairs, and the remaining products have both
indices at least five. Thus
\[
\begin{aligned}
 \mathcal T_9\left(H_+^2\right)
 ={}&2\sum_{i=0}^{4}b_i\sum_{j=9-i}^{9}b_j
       +\left(\sum_{j=5}^{9}b_j\right)^2\\
 &+2\delta_HP_H\left(1\right)+\delta_H^2
 <\frac{3336}{10^6}.
\end{aligned}
\]
Consequently the full quadratic tails satisfy
\[
\begin{aligned}
 \mathcal T_9\left(L^\#\right)
 &\le\mathcal T_{12}\mathcal N_+<\frac{13}{4},\\
 \mathcal T_9\left(N_e^\#\right)
 &\le28\mathcal T_9\left(H_+^2\right)
               +3\mathcal T_{11}\mathcal N_+\\
 &<28\cdot.003336+3\cdot20.136
   =60.501408<61.
\end{aligned}
\]
\begin{proposition}\label{prop:quadratic-complete-tails}
The complete-series remainders of the two quadratic numerators satisfy
\[
 \mathcal T_9\left(L^\#\right)<\frac{13}{4},
 \qquad
 \mathcal T_9\left(N_e^\#\right)<61.
\]
\end{proposition}
\begin{proof}
The coefficientwise majorants $U_+,V_+$ constructed in
Appendix~\ref{subsec:appB-profiles} imply
$L^\#\preccurlyeq S_3\mathcal N_+$ and
\[
 \mathcal T_9(L^\#)\le\mathcal T_{12}(\mathcal N_+).
\]
The ten-term table for $\mathcal N_+$ gives
$\mathcal T_{12}(\mathcal N_+)<3.222<13/4$. For the second
numerator, the identity and estimate immediately preceding the
proposition give
\[
 \begin{aligned}
 \mathcal T_9(N_e^\#)
 &\le28\mathcal T_9(H_+^2)+3\mathcal T_{11}(\mathcal N_+)\\
 &<28\cdot0.003336+3\cdot20.136\\
 &=60.501408<61.
 \end{aligned}
\]
Both asserted bounds follow.
\end{proof}

For the higher-order tails, represent each positive rational function
as a finite sum
\[
 F=\sum_{j,p}c_{j,p}\frac{y^j}{\left(1-\vartheta y\right)^p},
 \qquad c_{j,p}\ge0,
\]
where $j,p$ are nonnegative integers. At addition and multiplication,
we collect terms with the same pair $\left(j,p\right)$, without
combining different denominator powers. Differentiation and deletion
of a coefficient are specified by
\[
\begin{aligned}
 \partial_y\frac{y^j}{\left(1-\vartheta y\right)^p}
 &=\frac{j y^{j-1}}{\left(1-\vartheta y\right)^p}
   +\frac{p\vartheta y^j}{\left(1-\vartheta y\right)^{p+1}},\\
 S_1\frac{y^j}{\left(1-\vartheta y\right)^p}
 &=\frac{y^{j-1}}{\left(1-\vartheta y\right)^p}
                \quad\left(j\ge1\right),\\
 S_1\frac1{\left(1-\vartheta y\right)^p}
 &=\vartheta\sum_{r=1}^{p}\frac1{\left(1-\vartheta y\right)^r}
                \quad\left(p\ge1\right),\qquad S_1\left(1\right)=0.
\end{aligned}
\]
For the last identity, put $z=\vartheta y$ and use the finite
geometric identity
\[
 \left(1-z\right)^{-p}-1
 =z\sum_{r=1}^{p}\left(1-z\right)^{-r}.
\]
Subtracting the constant coefficient and dividing by $y$ proves the
formula for $S_1$. Thus every operation preserves nonnegativity.

The finite coefficient arrays can be calculated explicitly as follows.
If $F$ has entries $c_{j,p}$ and $G$ has entries $d_{j,p}$, their
product has entries
\[
 \left(c*d\right)_{j,p}
 =\sum_{a=0}^{j}\sum_{b=0}^{p}c_{a,b}d_{j-a,p-b}.
\]
Entries not present in a finite array are zero. The entry of $F'$ at
$\left(j,p\right)$ is
\[
 \left(j+1\right)c_{j+1,p}
       +\left(p-1\right)\vartheta c_{j,p-1}\quad\left(p\ge1\right),
\]
and it is $\left(j+1\right)c_{j+1,0}$ when $p=0$.
For $S_1F$, the entry at $j\ge1$ is $c_{j+1,p}$, whereas the
entries at $j=0$ are
\[
 c_{1,0}\quad\left(p=0\right),\qquad
 c_{1,p}+\vartheta\sum_{q\ge p}c_{0,q}\quad\left(p\ge1\right).
\]
Iteration gives $S_k$. Starting with the coefficients of $P_u,P_v$
at denominator power zero and the single entry $\mathfrak a$ at
$\left(12,1\right)$ gives all subsequent arrays by these finite
positive rational operations.

There is also a direct geometric bound for the tail of each term.
If $j\ge k$, its entire coefficient sum is
$c_{j,p}\left(1-\vartheta\right)^{-p}$. If $j<k$ and $p=0$, its tail
is zero. Otherwise put $r=k-j\ge1$. The binomial series gives a first
tail term $c_{j,p}\binom{p+r-1}{r}\vartheta^r$; every subsequent ratio
is no larger than $\vartheta\left(p+r\right)/\left(r+1\right)$.
Hence
\[
 \mathcal T_k\left(\frac{c_{j,p}y^j}
                          {\left(1-\vartheta y\right)^p}\right)
 \le\frac{c_{j,p}\binom{p+r-1}{r}\vartheta^r}
              {1-\vartheta\left(p+r\right)/\left(r+1\right)}.
\]
The ratio of successive binomial-series terms, starting with index
$r$, is
\[
 \vartheta\frac{p+r}{r+1}
 =\vartheta\left(1+\frac{p-1}{r+1}\right).
\]
It is nonincreasing in $r$ because $p\ge1$. In the functions below
$p\le12$ and $r\ge1$, so this ratio is at most $13/50<1$.
Every denominator in the geometric estimate is therefore at least
$37/50$. For example, the input term with $\left(j,p\right)=\left(12,1\right)$
has tail from degree twenty-one equal to
\[
 \frac{\mathfrak a\vartheta^9}{1-\vartheta},
\]
whereas the product of two such terms has $j=24$ and its entire
coefficient sum $\mathfrak a^2/\left(1-\vartheta\right)^2$ belongs
to that tail. These examples illustrate the two different cases
$j<k$ and $j\ge k$ in the estimate.

\subsection{The complete determinant tail}\label{subsec:appB-determinant}
The actual intermediate numerators to be controlled are
\[
 \mathcal C=4H^2N_a,\qquad \mathcal E=N_e,\qquad
 \mathcal Z=N_d,\qquad N_{\det}=\mathcal C\mathcal E-\mathcal Z^2.
\]
The calligraphic symbol $\mathcal C$ is distinct from $C=-1/r$ in the
calibration identity. Set $k_+=S_1V$ and define the positive functions
\[
\begin{aligned}
 D_+={}&2X_+\left[y\left(h_+'H_++h_+H_+'\right)+h_+H_+\right]
                  +yh_+H_+,\\
 N_{g,+}={}&y\left(k_+h_+'+h_+k_+'\right)+k_+h_+,\\
 \mathcal A_+={}&S_2\left[\frac74D_+^2+28X_+N_{g,+}^2\right],
 \qquad \mathcal L_+=S_3\mathcal N_+,\\
 \mathcal C_+={}&4H_+^2\mathcal A_+,\qquad
 \mathcal E_+=28H_+^2+3y\mathcal L_+,\\
 \mathcal Z_+={}&\mathcal L_+D_+
    +X_+h_+\left[y\left(\mathcal L_+'H_+
                         +2\mathcal L_+H_+'\right)+\mathcal L_+H_+\right].
\end{aligned}
\]
Coefficient convolution, differentiation, and the already proved
endpoint factors imply
\[
 \mathcal C^\#\preccurlyeq\mathcal C_+,\qquad
 \mathcal E^\#\preccurlyeq\mathcal E_+,\qquad
 \mathcal Z^\#\preccurlyeq\mathcal Z_+.
\]
In particular, the deletions $S_2,S_3$ in the positive majorants do
not assert cancellation of their low coefficients: they simply
delete those coefficients, as prescribed above. For example,
\[
 \left(KD^2-TxN_g^2\right)^\#
 \preccurlyeq\frac74D_+^2+28X_+N_{g,+}^2.
\]
The numerator on the left has zero coefficients at degrees zero and
$1$. Shifting both coefficient sequences by two therefore bounds
$N_a^\#$ by $\mathcal A_+$. The same argument with three deleted
coefficients bounds $L^\#$ by $\mathcal L_+$. Applying convolution
and differentiation to these two bounds then gives the displayed
majorants for $\mathcal C,\mathcal E,\mathcal Z$; no sign of an
individual coefficient of the actual numerators is assumed.

For the positive representation of $F$, define
\[
\begin{aligned}
 \mathcal B_0\left(F\right)
 &=\sum_{j\ge21}c_{j,0},\\
 \mathcal B_p\left(F\right)
 &=\sum_{j\ge21}c_{j,p}\left(\frac{25}{24}\right)^p\\
 &\quad+\sum_{j<21}
 \frac{c_{j,p}\binom{p+20-j}{21-j}25^{j-21}}
      {1-\left(p+21-j\right)/\left[25\left(22-j\right)\right]}
       \quad\left(p\ge1\right).
\end{aligned}
\]
The preceding geometric estimate proves
$\mathcal T_{21}F\le\sum_p\mathcal B_p\left(F\right)$.
Only the following finite ranges occur:
\begingroup
\setlength{\tabcolsep}{\arraycolsep}
\setlength{\LTpre}{\medskipamount}
\setlength{\LTpost}{\medskipamount}
\begin{longtable}{>{$}c<{$}|>{$}c<{$}>{$}c<{$}}
F&\max j&\max p\\
\hline
\endfirsthead
F&\max j&\max p\\
\hline
\endhead
\hline
\endfoot
\hline
\endlastfoot
\mathcal C_+&66&12\\
\mathcal E_+&24&6\\
\mathcal Z_+&46&10\\
\end{longtable}
\endgroup
Substitution of $P_u,P_v,\mathfrak a,\vartheta$ in the positive
operations gives the following strict rational bounds:
\begingroup
\setlength{\tabcolsep}{\arraycolsep}
\setlength{\LTpre}{\medskipamount}
\setlength{\LTpost}{\medskipamount}
\begin{longtable}{>{$}c<{$}|>{$}c<{$}>{$}c<{$}>{$}c<{$}}
&\mathcal C_+&\mathcal E_+&\mathcal Z_+\\
\hline
\endfirsthead
&\mathcal C_+&\mathcal E_+&\mathcal Z_+\\
\hline
\endhead
\hline
\endfoot
\hline
\endlastfoot
\mathcal B_0&171/100&0&63/10000\\
\mathcal B_1&153/10000&4/10^{11}&17/10000\\
\mathcal B_2&21/100000&15/10^{11}&8/100000\\
\mathcal B_3&11/10^7&3/10^{11}&25/10^7\\
\sum_{p\ge4}\mathcal B_p&1/10^{10}&1/10^{12}&5/10^8\\
\end{longtable}
\endgroup
The zero entry is exact. In fact, the denominator-power-zero parts
of $U,V$ have degree eleven, so the corresponding part of
$\mathcal N_+$ has degree at most twenty-two. After applying $S_3$
and multiplying by $y$, its degree is at most twenty. The polynomial
part of $H_+^2$ has degree at most eighteen. Thus the polynomial part
of $\mathcal E_+$ has no coefficient of degree twenty-one or higher.
Terms that contain a positive power of the denominator still contain
such a power after these operations. They therefore do not affect this
conclusion.
The other entries are obtained by summing, for each fixed $p$, the
finite array entries with the geometric weights in
$\mathcal B_p$. Their column sums are explicitly
\[
\begin{aligned}
 \sum_p\mathcal B_p\left(\mathcal C_+\right)
 &<1.7255111001<\frac74,\\
 \sum_p\mathcal B_p\left(\mathcal E_+\right)
 &<\frac{221}{10^{12}}<10^{-9},\\
 \sum_p\mathcal B_p\left(\mathcal Z_+\right)
 &<.00808255<\frac1{100}.
\end{aligned}
\]
Combining these estimates with coefficientwise domination gives
\[
 \mathcal T_{21}\left(\mathcal C^\#\right)<\frac74,\qquad
 \mathcal T_{21}\left(\mathcal E^\#\right)<10^{-9},\qquad
 \mathcal T_{21}\left(\mathcal Z^\#\right)<\frac1{100}.
\]
\begin{proposition}\label{prop:determinant-complete-tails}
The complete-series tails of the three factors in
$N_{\det}=\mathcal C\mathcal E-\mathcal Z^2$ satisfy
\[
 \mathcal T_{21}\left(\mathcal C^\#\right)<\frac74,\qquad
 \mathcal T_{21}\left(\mathcal E^\#\right)<10^{-9},\qquad
 \mathcal T_{21}\left(\mathcal Z^\#\right)<\frac1{100}.
\]
\end{proposition}
\begin{proof}
The coefficientwise majorants
\[
 \mathcal C^\#\preccurlyeq\mathcal C_+,\qquad
 \mathcal E^\#\preccurlyeq\mathcal E_+,\qquad
 \mathcal Z^\#\preccurlyeq\mathcal Z_+
\]
were established above. For each positive rational function
$F=\sum_{j,p}c_{j,p}y^j/(1-\vartheta y)^p$, the geometric
estimate proved in this subsection gives
\[
 \mathcal T_{21}(F)\le\sum_p\mathcal B_p(F).
\]
The finite rational arrays for $\mathcal C_+,\mathcal E_+,\mathcal Z_+$
are listed in the table above. Summing their columns gives, respectively,
\[
 \begin{aligned}
 \sum_p\mathcal B_p(\mathcal C_+)
   &<1.7255111001<\frac74,\\
 \sum_p\mathcal B_p(\mathcal E_+)
   &<\frac{221}{10^{12}}<10^{-9},\\
 \sum_p\mathcal B_p(\mathcal Z_+)
   &<0.00808255<\frac1{100}.
 \end{aligned}
\]
Combining these three estimates with the coefficientwise majorants
proves all the asserted tail bounds. The geometric estimate includes
both terms whose initial degree is already at least twenty-one and
terms of lower initial degree whose denominator expansion contributes
to the tail.
\end{proof}

\par\medskip\noindent
\textbf{AI usage disclosure.}
The authors acknowledge the use of ChatGPT (GPT-6 Astra Pro) for
mathematical discussions, editorial assistance, and selected technical
tasks.
In particular, AI tools assisted with the lengthy coefficient
expansions and algebraic calculations in Section~\ref{sec:certificate},
including the polynomial divisions, coefficient convolutions, and
determinant bounds used there. They also assisted with the rational
certificates collected in Appendix~\ref{sec:paper-input} and the
complete-series remainder estimates in Appendix~\ref{sec:tail}.
The schematic drawings in Figures~\ref{fig:two-balls},
\ref{fig:radial-map}, \ref{fig:hessian}, and \ref{fig:plateau}
were prepared with AI assistance. The authors checked the
AI-assisted arguments and computations and take full responsibility
for the contents of the paper.

\par\medskip\noindent
\textbf{Data availability.}
All quantities used in the proofs are specified in this article.
\par\medskip\noindent
\textbf{Conflict of interest.}
The authors declare that they have no conflicts of interest.

\end{document}